\documentclass[]{aspm}
\articleinfo{}{}{}
\usepackage{verbatim}
\usepackage{amssymb}
\usepackage{amsbsy}
\usepackage{amscd}
\usepackage{amsmath}
\usepackage{amsthm}
\usepackage[mathscr]{eucal}

\newtheorem{theorem}{Theorem}
\numberwithin{defn}{section}

\title[Differential complexes]
      {Some differential complexes within and beyond parabolic geometry}

\author[R.L.~Bryant]{Robert~L.~Bryant}
\author[M.G.~Eastwood]{Michael~G.~Eastwood}
\author[A.R.~Gover]{A.~Rod.~Gover}
\author[K. Neusser]{Katharina~Neusser}

\address{Mathematical Sciences Research Institute, 17 Gauss Way, Berkeley,\\ 
CA 94720-5070, USA, and\\
Department of Mathematics, University of California, Berkeley,\\
CA 94720-3840, USA}
\email{bryant@msri.org}

\address{Mathematical Sciences Institute, Australian National University,\\
ACT 0200, Australia}
\email{meastwoo@member.ams.org}

\address{Department of Mathematics, University of Auckland,
Private Bag 92019,\\ Auckland 1142, New Zealand, and\\
Mathematical Sciences Institute, Australian National University,\\
ACT 0200, Australia}
\email{r.gover@auckland.ac.nz}

\address{Mathematical Sciences Institute, Australian National University,\\
ACT 0200, Australia}
\email{katharina.neusser@anu.edu.au}

\rcvdate{}
\rvsdate{}


\begin{document}

\begin{abstract}
For smooth manifolds equipped with various geometric structures, we
construct complexes that replace the de~Rham complex in providing an
alternative fine resolution of the sheaf of locally constant functions. In case
that the geometric structure is that of a parabolic geometry, our complexes
coincide with the Bernstein-Gelfand-Gelfand complex associated with the trivial
representation. However, at least in the cases we discuss, our constructions
are relatively simple and avoid most of the machinery of parabolic geometry.
Moreover, our method extends to contact and symplectic geometries (beyond the 
parabolic realm).
\end{abstract}

\dedicatory{Dedicated to Professors Reiko Miyaoka and Keizo Yamaguchi on their
sixtieth birthdays}

\maketitle

\section{Introduction}
In~\cite{CSS}, \v{C}ap, Slov\'ak, and Sou\v{c}ek construct sequences of
invariant differential operators on {\em parabolic geometries\/} of any
type~$G/P$, one for each finite-dimensional representation ${\mathbb{V}}$
of~$G$. (Here, $G$ is a semisimple Lie group and $P\subset G$ a parabolic
subgroup.) These sequences are known as {\em Bernstein-Gelfand-Gelfand\/} (BGG)
sequences since, for the homogeneous model $G/P$ of such a geometry, these
sequences are complexes, which are dual to a parallel construction due to these
authors~\cite{BGG} on the level of Verma modules. In ~\cite{CD} Calderbank and
Diemer simplify the construction of BGG sequences in~\cite{CSS}. In addition
they provide~\cite[p.~87]{CD}, for regular parabolic geometries, alternative
BGG sequences, which only coincide with the ones in ~\cite{CSS} if the geometry
is {\em torsion-free\/}. The latter sequences not only appear to be more
natural, they also have the advantage that if ${\mathbb{V}}$ is taken to be the
trivial representation, then they form complexes, providing fine resolutions of
the locally constant sheaf~${\mathbb{R}}$ (as one sees by suitably modifying
\cite[Proposition~5.5(iv)]{CD}). For the sequences of \cite{CSS} this is only
true if the geometry is torsion-free and, in this case, the two sequences are
anyway the same. In combination with the construction of canonical Cartan
connections given in~\cite{CSc}, this shows that one can find alternatives to
the de~Rham resolution for any parabolic geometry defined in terms of a regular
infinitesimal flag structure \cite[\S3.1.6]{CSl}. A hallmark of these
resolutions is that the ranks of the bundles involved are diminished as
compared to the de~Rham complex. The price one pays is that the operators may
be higher than first order. The construction of these resolutions
in~\cite{CD,CSS}, entails firstly constructing the Cartan connection as
described in~\cite{CSc} and this is not at all straightforward.

In this article we present some examples constructed by a more elementary
route. As we show, our method extends to certain non-parabolic geometries,
namely arbitrary contact and symplectic geometries. We shall use the spectral
sequence of a filtered complex~\cite{Ch} without comment and merely as a
replacement for tedious diagram chasing.

\section{The Rumin complex}\label{rumin} 
For our first example we shall construct the Rumin complex~\cite{R}. It is
defined on an arbitrary contact manifold but, for simplicity, we shall present
the $5$-dimensional case, which is typical. So let $M$ be a $5$-dimensional
smooth manifold with $H\subset TM$ a contact distribution. Equivalently, the
contact structure may be defined by $L\equiv H^\perp$, a line sub-bundle of the
bundle of $1$-forms~$\Lambda^1$. If we define a rank $4$ vector bundle
$\Lambda_H^1$ as the quotient $\Lambda^1/L$, then there are induced short exact
sequences
$$0\to \Lambda_H^{p-1}\otimes L\to\Lambda^p\to\Lambda_H^p\to 0,
\quad\mbox{for }1\leq p\leq 5$$
and the spectral sequence of the de~Rham complex filtered in this way reads, at
the $E_0$-level, 
\small$$\begin{picture}(300,45)
\put(160,25){\makebox(0,0){$\begin{array}{ccccccccccccc}
\Lambda^0&&\Lambda_H^1&&\Lambda_H^2 &&\Lambda_H^3
&&\Lambda_H^4&&0\\
&&&&\makebox[0pt][r]{${\mathcal{L}}$}\uparrow&&\uparrow&&\uparrow\\
&&0&&L&&\Lambda_H^1\otimes L&&\Lambda_H^2\otimes L&
&\Lambda_H^3\otimes L&&\Lambda_H^4\otimes L
\end{array}$}}
\put(0,0){\vector(1,0){40}}
\put(0,0){\vector(0,1){45}}
\put(30,5){$p$}
\put(3,35){$q$}
\end{picture}$$\normalsize
where ${\mathcal{L}}$ is the composition
$L\to\Lambda^1\xrightarrow{\,d\,}\Lambda^2\to\Lambda_H^2$. The Leibniz rule
shows that ${\mathcal{L}}$ is linear over the functions and is, therefore, a
homomorphism of vector bundles. It is called the {\em Levi form\/}. By
definition of contact manifold, the range of ${\mathcal{L}}$ is non-degenerate
as a skew form on~$H$, defined up to scale. Equivalently, we can choose local
co-framings $(\omega_1,\omega_2,\omega_3,\omega_4,\omega_5)$ with $\omega_1$ a
section of $L$ such that
$$d\omega_1\equiv \omega_2\wedge\omega_3+\omega_4\wedge\omega_5\bmod\omega_1.$$
Notice that
\small$$\begin{array}{cc}
d(\omega_1\wedge\omega_2)
\equiv\omega_2\wedge\omega_4\wedge\omega_5\bmod\omega_1&\quad
d(\omega_1\wedge\omega_3)
\equiv\omega_3\wedge\omega_4\wedge\omega_5\bmod\omega_1\\
d(\omega_1\wedge\omega_4)
\equiv\omega_2\wedge\omega_3\wedge\omega_4\bmod\omega_1&\quad
d(\omega_1\wedge\omega_5)
\equiv\omega_2\wedge\omega_3\wedge\omega_5\bmod\omega_1
\end{array}$$\normalsize
whence the $E_0$-differential $\Lambda_H^1\otimes L\to \Lambda_H^3$ is an
isomorphism of vector bundles. Similar reasoning shows that 
$\Lambda_H^2\otimes L\to\Lambda_H^4$ is surjective. Hence, at the $E_1$-level
we obtain 
\small$$\begin{picture}(300,45)
\put(163,25){\makebox(0,0){$\begin{array}{ccccccccccccc}
\Lambda^0&\!\!\!\stackrel{d_\perp}{\longrightarrow}\!\!\!&\Lambda_H^1&
\!\!\!\stackrel{d_\perp}{\longrightarrow}\!\!\!&
\Lambda_{H\perp}^2&&\!\!0\!\!&&0&&0\\[4pt]
&&0&&0&&\!\!0\!\!&&\Lambda_{H\perp}^2\otimes L&
\!\!\!\stackrel{d_\perp}{\longrightarrow}\!\!\!
&\Lambda_H^3\otimes L&
\!\!\!\stackrel{d_\perp}{\longrightarrow}\!\!\!&\Lambda_H^4\otimes L
\end{array}$}}
\put(0,0){\vector(1,0){40}}
\put(0,0){\vector(0,1){45}}
\put(30,5){$p$}
\put(3,35){$q$}
\end{picture}$$\normalsize
and deduce that  there is a complex
\small\begin{equation}\label{rumincomplex}
0\to{\mathbb{R}}\to\Lambda^0\stackrel{d_\perp}{\longrightarrow}\Lambda_H^1
\stackrel{d_\perp}{\longrightarrow}\Lambda_{H\perp}^2
\stackrel{d_\perp^{(2)}}{\longrightarrow}\Lambda_{H\perp}^2\otimes L
\stackrel{d_\perp}{\longrightarrow}\Lambda_H^3\otimes L
\stackrel{d_\perp}{\longrightarrow}\Lambda_H^4\otimes L\to 0
\end{equation}\normalsize
where $\Lambda_{H\perp}^p$ denotes the sub-bundle of~$\Lambda_H^p$, trace-free
with respect to the Levi-form. The operator $d_\perp^{(2)}$ is second order
and, because the spectral sequence converges to the local cohomology of the
de~Rham complex, it follows that this complex is exact on the level of sheaves.
Already, the Rumin complex goes beyond parabolic geometry. Notice that,
although a convenient co-framing was chosen to perform some calculation,
the construction itself and the resulting complex are independent of any such 
choice. This is a repeated theme in this article.

\section{The Engel complex}\label{engel}
In this section we shall be concerned with a smooth $4$-manifold $M$ equipped
with a generic distribution $H\subset TM$ of rank~$2$. Genericity entails that
$[H,H]$ has rank~$3$ and that $[H,[H,H]]=TM$. Dually, if we let $K\equiv
H^\perp$ and $L\equiv [H,H]^\perp$ then the $1$-forms are filtered $L\subset
K\subset\Lambda^1$ by the line-bundle $L$ and rank $2$ bundle~$K$. In fact,
there is a canonically defined finer filtration constructed as follows. One
easily checks that the Levi form $K\to\Lambda_H^2$, defined as the composition
$K\to\Lambda^1\xrightarrow{\,d\,}\Lambda^2\to\Lambda_H^2$, is a surjective
homomorphism of vector bundles with $L$ as kernel. It follows that the other
Levi form, defined as the composition
$$L\to\Lambda^1\xrightarrow{\,d\,}\Lambda^2\to\Lambda^2(\Lambda^1/L)$$ 
has range in the kernel of $\Lambda^2(\Lambda^1/L)\to\Lambda_H^2$. 
However, the short exact sequence
$$0\to K/L\to\Lambda^1/L\to\Lambda_H^1\to0$$
identifies this kernel as $\Lambda_H^1\otimes K/L$. In other words, we have a
canonically defined inclusion $L\otimes(K/L)^*\hookrightarrow\Lambda_H^1$ the
range of which defines a line sub-bundle $\xi$ of $\Lambda_H^1$. The result is
that we can write
$$\Lambda^1=\Lambda_H^1/\xi+\xi+K/L+L,$$
meaning that $\Lambda^1$ is filtered with composition factors being line
bundles as indicated (ordered so that $\Lambda_H^1/\xi$ is a canonical quotient
and $L$ is a canonical sub-bundle). All in all, if we write $\lambda$ for
$\Lambda_H^1/\xi$ and untangle the identifications found above, then we
conclude that
\begin{equation}\label{filtration}\Lambda^1=
\lambda+\xi+\lambda\xi+\lambda\xi^2.\end{equation}
Equivalently, we can work locally with $(\omega^1,\omega^2,\omega^3,\omega^4)$,
an {\em adapted\/} co-framing such that
\begin{equation}\label{structure}
d\omega^1\equiv\omega^2\wedge\omega^3\bmod\omega^1\quad\mbox{and}\quad
d\omega^2\equiv\omega^3\wedge\omega^4\bmod \omega^1,\omega^2,\end{equation}
noting that the freedom in such a co-framing comprises exactly the triangular 
endomorphisms of the 
filtration~(\ref{filtration}), where 
$$\begin{array}{c}
L=\lambda\xi^2={\mathrm{span}}\{\omega^1\},\enskip
K=\lambda\xi+\lambda\xi^2={\mathrm{span}}\{\omega^1,\omega^2\},\\[4pt]
\xi+K={\mathrm{span}}\{\omega^1,\omega^2\,\omega^3\}.\end{array}$$
So far, this is the structure of an {\em Engel manifold\/}. 
As with the Rumin complex, it is clear that the first order operator 
$d_H:\Lambda^0\to\Lambda_H^1$ defined as the composition 
$\Lambda^0\xrightarrow{\,d\,}\Lambda^1\to\Lambda_H^1$ has the locally constant 
functions as its kernel. We now seek differential conditions on a section of 
$\Lambda_H^1$ in order that it be in the range of the operator~$d_H$.
Starting with any $1$-form~$\omega$,
\begin{equation}\label{primary_recipe}
\begin{array}{l}
\bullet\enskip\mbox{define }f\mbox{ by }
d\omega \equiv f\,\omega^3\wedge\omega^4 \bmod \omega^1,\omega^2\\
\bullet\enskip\mbox{define }p\mbox{ by } 
d(\omega-f\omega^2)\equiv p\,\omega^2\wedge\omega^4+
                           g\,\omega^2\wedge\omega^3\bmod \omega^1.
\end{array}\end{equation}
The structure equations (\ref{structure}) show that $p$ is well-defined and 
one easily checks that the equivalence class
$$[p\,\omega^2\wedge\omega^4]\in
\frac{{\mathrm{span}}\{\omega^2\wedge\omega^4,\omega^2\wedge\omega^3,
\omega^1\wedge\omega^4,\omega^1\wedge\omega^3,\omega^1\wedge\omega^2\}}
{{\mathrm{span}}\{\omega^2\wedge\omega^3,
\omega^1\wedge\omega^4,\omega^1\wedge\omega^3,\omega^1\wedge\omega^2\}}
\cong\lambda^2\xi$$
depends only on the equivalence class $[\omega]\in\Lambda_H^1$ and is
independent of choice of co-framing. We have a well-defined second 
order differential operator
$$\Lambda_H^1\ni[\omega]\stackrel{{\mathcal{P}}}{\longmapsto}
[p\,\omega^2\wedge\omega^4]\in\lambda^2\xi,$$
giving what we shall call the {\em primary\/} obstruction to $[\omega]$ being
in the range of~$d_H$. In a chosen co-frame, one can easily proceed to find a
{\em secondary\/} obstruction $s$ as follows. Define $f,p,g$ by
(\ref{primary_recipe}) and then
$$\bullet\enskip\mbox{define }s\mbox{ by }
d(\omega-f\omega^2-g\omega^1)=p\,\omega^2\wedge\omega^4
+r\,\omega^1\wedge\omega^4+s\,\omega^1\wedge\omega^3
+t\,\omega^1\wedge\omega^2.$$
If $p$ vanishes, then 
$$0=d^2(\omega-f\omega^2-g\omega^1)=r\,\omega^2\wedge\omega^3\wedge\omega^4+
\cdots$$
so $r$ vanishes. If, in addition $s$ vanishes, then 
$$0=d^2(\omega-f\omega^2-g\omega^1)=d(t\,\omega^1\wedge\omega^2)=
-t\,\omega^1\wedge\omega^3\wedge\omega^4$$
so $t$ vanishes. Hence, if both $p$ and $s$ vanish, then
$d(\omega-f\omega^2-g\omega^1)=0$. By the Poincar\'e Lemma, it follows that
$[\omega]$ is locally in the range of~$d_H$, as required. If the primary 
obstruction vanishes, then the equivalence class
$$[s\,\omega^1\wedge\omega^3]\in
\frac{{\mathrm{span}}\{\omega^1\wedge\omega^3,\omega^1\wedge\omega^2\}}
{{\mathrm{span}}\{\omega^1\wedge\omega^2\}}
\cong\lambda\xi^3$$
is independent of choice of co-framing. Otherwise, the change
\begin{equation}\label{badchange}
\omega_4\mapsto\omega_4+h\omega_3\end{equation} 
induces severe complications with $s$ changing by $r$ and its derivatives. If
one wants to avoid these complications, it suffices to prohibit 
(\ref{badchange}) to arrive at an invariantly defined differential operator
$$({\mathcal{P}},{\mathcal{S}}):\Lambda_H^1\to
\lambda^2\xi\oplus\lambda\xi^3,$$
whose kernel is locally the range of~$d_H$. More precisely, we may eliminate 
(\ref{badchange}) by choosing a complement to the line sub-bundle
$\xi\hookrightarrow\Lambda_H^1$. In other words, we choose a splitting
$\Lambda_H^1=\lambda\oplus\xi$. An adapted co-framing yields such a splitting
and, conversely, a fixed choice of splitting restricts the choice of adapted
co-framings precisely by preventing the addition of any multiple of
$\omega^3$ to $\omega^4$. The forms on an Engel manifold endowed with this
extra structure are filtered as follows.
$$\begin{array}{rcl} 
\Lambda^1&=&(\lambda\oplus\xi)+\lambda\xi+\lambda\xi^2\\[2pt]
\Lambda^2&=&\lambda\xi+(\lambda^2\xi\oplus\lambda\xi^2)+
(\lambda^2\xi^2\oplus\lambda\xi^3)+\lambda^2\xi^3\\[2pt]
\Lambda^3&=&\lambda^2\xi^2+\lambda^2\xi^3+(\lambda^3\xi^3\oplus\lambda^2\xi^4)
\end{array}$$
and the spectral sequence of the de~Rham complex filtered in this way reads, at
the $E_0$-level, 
\small$$\begin{picture}(300,85)(0,5)
\put(170,50){\makebox(0,0){$\begin{array}{cccccccc}
\Lambda^0&\lambda\oplus\xi&\lambda\xi &0&0&0&0\\
&&\uparrow\\
&0&\lambda\xi&\lambda^2\xi\oplus\lambda\xi^2&\lambda^2\xi^2&0&0\\
&&&\uparrow&\uparrow\\
&0&0&\lambda\xi^2&\lambda^2\xi^2\oplus\lambda\xi^3&
\lambda^2\xi^3&0\\
&&&&&\uparrow\\
&0&0&0&0&\lambda^2\xi^3&\lambda^3\xi^3\oplus\lambda^2\xi^4&\Lambda^4\,.
\end{array}$}}
\put(0,0){\vector(1,0){40}}
\put(0,0){\vector(0,1){45}}
\put(30,5){$p$}
\put(3,35){$q$}
\end{picture}$$\normalsize
The $E_0$-differentials are easily computed in our adapted co-frame. For 
example
\small$$\begin{array}{rcl}
f\,\omega^1\wedge\omega^2&\!\stackrel{d}{\longrightarrow}\!&
f\,\omega^1\wedge\omega^3\wedge\omega^4
\scriptstyle\bmod\omega^1\wedge\omega^2\\
h\,\omega^1\wedge\omega^4+g\,\omega^1\wedge\omega^3+f\,\omega^1\wedge\omega^2
&\!\stackrel{d}{\longrightarrow}\!&
h\,\omega^2\wedge\omega^3\wedge\omega^4
\scriptstyle\bmod\omega^1\wedge\omega^2,\,\omega^1\wedge\omega^3\wedge\omega^4
\end{array}$$\normalsize
deals with the two rightmost differentials. Consequently, at the $E_1$-level we
obtain 
$$\begin{picture}(300,75)(0,5)
\put(165,45){\makebox(0,0){$\begin{array}{ccccccccccccccc}
\Lambda^0&\!&\lambda\oplus\xi&&0&&0&&0&&0&&0\\[4pt]
&\!&0&&0&&\lambda^2\xi&&0&&0&&0\\[4pt]
&\!&0&&0&&0&&\lambda\xi^3&&0&&0\\[4pt]
&\!&0&&0&&0&&0&&0&&\lambda^3\xi^3\oplus\lambda^2\xi^4&\!&\Lambda^4\,.
\end{array}$}}
\put(0,0){\vector(1,0){40}}
\put(0,0){\vector(0,1){45}}
\put(30,5){$p$}
\put(3,35){$q$}
\end{picture}$$
The bundles $\lambda\oplus\xi$ and $\lambda^3\xi^3\oplus\lambda^2\xi^4$ may be
identified with $\Lambda_H^1$ and $\Lambda_H^1\otimes\Lambda^2\!K$,
respectively. The line bundles $\lambda^2\xi$ and $\lambda\xi^3$ combine to
give a rank $2$ vector bundle $\lambda^2\xi+\lambda\xi^3$ but, in fact, this
bundle canonically splits as can readily be seen in our adapted co-frame:
$$\textstyle\lambda^2\xi+\lambda\xi^3=
\frac{{\mathrm{span}}\{\omega^2\wedge\omega^4,\omega^1\wedge\omega^2,
\omega^1\wedge\omega^4\}}
{{\mathrm{span}}\{\omega^1\wedge\omega^2,\omega^1\wedge\omega^4\}}\oplus
\frac{{\mathrm{span}}\{\omega^1\wedge\omega^3,\omega^1\wedge\omega^2\}}
{{\mathrm{span}}\{\omega^1\wedge\omega^2\}},$$
independent of choice of co-frame. We conclude that there is a complex of
differential operators (cf.~\cite{P})
\begin{equation}\label{fullEngel}
\Lambda^0\xrightarrow{\,d_H\,}\Lambda_H^1\to \lambda^2\xi\oplus\lambda\xi^3
\to\Lambda_H^1\otimes\Lambda^2\!K\to\Lambda^4\to 0\end{equation}
resolving the locally constant sheaf~${\mathbb{R}}$. Following through the
spectral sequence more explicitly as a diagram chase shows that 
$\Lambda_H^1\to\lambda^2\xi\oplus\lambda\xi^3$ is given by our previous 
recipe.  We shall see later in \S\ref{engelrevisited}, that (\ref{fullEngel}) 
is a BGG complex for an appropriate parabolic geometry.

\section{The Rumin complex revisited}
Since a contact manifold with no extra structure is not a parabolic geometry,
the Rumin complex lies outside the realm of parabolic geometry. Nevertheless,
there is a parabolic geometry in which the Rumin complex finds its genesis. Let
us denote by ${\mathrm{Sp}}(2n,{\mathbb{R}})$ the simple Lie group of linear
automorphisms of ${\mathbb{R}}^{2n}$ preserving a fixed non-degenerate
symplectic form. Viewing the $(2n+1)$-sphere $S^{2n+1}$ as
$$\{x\in{\mathbb{R}}^{2n+2}\mbox{ s.t. }x\not=0\}/\{x\sim\lambda x
\mbox{ for }\lambda>0\}$$
(i.e.~the space of rays emanating from the origin in ${\mathbb{R}}^{2n+2}$),
the group $G={\mathrm{Sp}}(2n+2,{\mathbb{R}})$ acts smoothly and transitively
on $S^{2n+1}$. The stabiliser subgroup $P$ of this action is parabolic.
Parabolic geometries modelled on this particular homogeneous space
$S^{2n+1}=G/P$ are known as {\em contact projective\/}~\cite[\S4.2.6]{CSl}. In
any case, when viewed in this way, the sphere $S^{2n+1}$ inherits a
$G$-invariant contact structure from the symplectic form on
${\mathbb{R}}^{2n+2}$. As in~\S\ref{rumin}, let us now consider the case $n=2$.
Adopting the notation from~\cite{BE}, this homogeneous space is written as
\begin{picture}(42,5)
\put(20,.1){\line(1,0){16}}
\put(20,3.3){\line(1,0){16}}
\put(20,1.5){\makebox(0,0){$\bullet$}}
\put(36,1.5){\makebox(0,0){$\bullet$}}
\put(28,1.5){\makebox(0,0){$\langle$}}
\put(4,1.5){\makebox(0,0){$\times$}}
\put(4,1.5){\line(1,0){15}}
\end{picture} and the Bernstein-Gelfand-Gelfand complex corresponding to the 
trivial representation of ${\mathrm{Sp}}(6,{\mathbb{R}})$ is
$$\begin{array}{l}0\to
\begin{picture}(42,5)
\put(20,.1){\line(1,0){16}}
\put(20,3.3){\line(1,0){16}}
\put(20,1.5){\makebox(0,0){$\bullet$}}
\put(36,1.5){\makebox(0,0){$\bullet$}}
\put(28,1.5){\makebox(0,0){$\langle$}}
\put(20,8){\makebox(0,0){$\scriptstyle 0$}}
\put(36,8){\makebox(0,0){$\scriptstyle 0$}}
\put(4,1.5){\makebox(0,0){$\bullet$}}
\put(4,1.5){\line(1,0){15}}
\put(4,8){\makebox(0,0){$\scriptstyle 0$}}
\end{picture}\to
\begin{picture}(42,5)
\put(20,.1){\line(1,0){16}}
\put(20,3.3){\line(1,0){16}}
\put(20,1.5){\makebox(0,0){$\bullet$}}
\put(36,1.5){\makebox(0,0){$\bullet$}}
\put(28,1.5){\makebox(0,0){$\langle$}}
\put(20,8){\makebox(0,0){$\scriptstyle 0$}}
\put(36,8){\makebox(0,0){$\scriptstyle 0$}}
\put(4,1.5){\makebox(0,0){$\times$}}
\put(4,1.5){\line(1,0){15}}
\put(4,8){\makebox(0,0){$\scriptstyle 0$}}
\end{picture}\to 
\begin{picture}(42,5)
\put(20,.1){\line(1,0){16}}
\put(20,3.3){\line(1,0){16}}
\put(20,1.5){\makebox(0,0){$\bullet$}}
\put(36,1.5){\makebox(0,0){$\bullet$}}
\put(28,1.5){\makebox(0,0){$\langle$}}
\put(20,8){\makebox(0,0){$\scriptstyle 1$}}
\put(36,8){\makebox(0,0){$\scriptstyle 0$}}
\put(4,1.5){\makebox(0,0){$\times$}}
\put(4,1.5){\line(1,0){15}}
\put(4,8){\makebox(0,0){$\scriptstyle -2$}}
\end{picture}\to
\begin{picture}(42,5)
\put(20,.1){\line(1,0){16}}
\put(20,3.3){\line(1,0){16}}
\put(20,1.5){\makebox(0,0){$\bullet$}}
\put(36,1.5){\makebox(0,0){$\bullet$}}
\put(28,1.5){\makebox(0,0){$\langle$}}
\put(20,8){\makebox(0,0){$\scriptstyle 0$}}
\put(36,8){\makebox(0,0){$\scriptstyle 1$}}
\put(4,1.5){\makebox(0,0){$\times$}}
\put(4,1.5){\line(1,0){15}}
\put(4,8){\makebox(0,0){$\scriptstyle -3$}}
\end{picture}\\[8pt] 
\mbox{ }\hspace{100pt}\to
\begin{picture}(42,5)
\put(20,.1){\line(1,0){16}}
\put(20,3.3){\line(1,0){16}}
\put(20,1.5){\makebox(0,0){$\bullet$}}
\put(36,1.5){\makebox(0,0){$\bullet$}}
\put(28,1.5){\makebox(0,0){$\langle$}}
\put(20,8){\makebox(0,0){$\scriptstyle 0$}}
\put(36,8){\makebox(0,0){$\scriptstyle 1$}}
\put(4,1.5){\makebox(0,0){$\times$}}
\put(4,1.5){\line(1,0){15}}
\put(4,8){\makebox(0,0){$\scriptstyle -5$}}
\end{picture}\to
\begin{picture}(42,5)
\put(20,.1){\line(1,0){16}}
\put(20,3.3){\line(1,0){16}}
\put(20,1.5){\makebox(0,0){$\bullet$}}
\put(36,1.5){\makebox(0,0){$\bullet$}}
\put(28,1.5){\makebox(0,0){$\langle$}}
\put(20,8){\makebox(0,0){$\scriptstyle 1$}}
\put(36,8){\makebox(0,0){$\scriptstyle 0$}}
\put(4,1.5){\makebox(0,0){$\times$}}
\put(4,1.5){\line(1,0){15}}
\put(4,8){\makebox(0,0){$\scriptstyle -6$}}
\end{picture}\to
\begin{picture}(42,5)
\put(20,.1){\line(1,0){16}}
\put(20,3.3){\line(1,0){16}}
\put(20,1.5){\makebox(0,0){$\bullet$}}
\put(36,1.5){\makebox(0,0){$\bullet$}}
\put(28,1.5){\makebox(0,0){$\langle$}}
\put(20,8){\makebox(0,0){$\scriptstyle 0$}}
\put(36,8){\makebox(0,0){$\scriptstyle 0$}}
\put(4,1.5){\makebox(0,0){$\times$}}
\put(4,1.5){\line(1,0){15}}
\put(4,8){\makebox(0,0){$\scriptstyle -6$}}
\end{picture}\to 0.\end{array}$$
This coincides with the Rumin complex (\ref{rumincomplex}). The reason for the
notation is fully explained in~\cite{BE}. Here, suffice it to say that
$$\begin{picture}(24,5)
\put(4,.1){\line(1,0){16}}
\put(4,3.3){\line(1,0){16}}
\put(4,1.5){\makebox(0,0){$\bullet$}}
\put(20,1.5){\makebox(0,0){$\bullet$}}
\put(12,1.5){\makebox(0,0){$\langle$}}
\put(4,8){\makebox(0,0){$\scriptstyle q$}}
\put(20,8){\makebox(0,0){$\scriptstyle r$}}
\end{picture}(\Lambda_H^1)\otimes L^s=
\begin{picture}(52,5)(-20,0)
\put(30,.1){\line(1,0){16}}
\put(30,3.3){\line(1,0){16}}
\put(30,1.5){\makebox(0,0){$\bullet$}}
\put(46,1.5){\makebox(0,0){$\bullet$}}
\put(38,1.5){\makebox(0,0){$\langle$}}
\put(30,8){\makebox(0,0){$\scriptstyle q$}}
\put(46,8){\makebox(0,0){$\scriptstyle r$}}
\put(4,1.5){\makebox(0,0){$\times$}}
\put(4,1.5){\line(1,0){25}}
\put(0,8){\makebox(0,0){$\scriptstyle -2s-2q-3r$}}
\end{picture}$$
where 
$\begin{picture}(24,5)
\put(4,.1){\line(1,0){16}}
\put(4,3.3){\line(1,0){16}}
\put(4,1.5){\makebox(0,0){$\bullet$}}
\put(20,1.5){\makebox(0,0){$\bullet$}}
\put(12,1.5){\makebox(0,0){$\langle$}}
\put(4,8){\makebox(0,0){$\scriptstyle q$}}
\put(20,8){\makebox(0,0){$\scriptstyle r$}}
\end{picture}(\Lambda_H^1)$ denotes
the bundle induced by the irreducible representation 
$\begin{picture}(24,6)
\put(4,.1){\line(1,0){16}}
\put(4,3.3){\line(1,0){16}}
\put(4,1.5){\makebox(0,0){$\bullet$}}
\put(20,1.5){\makebox(0,0){$\bullet$}}
\put(12,1.5){\makebox(0,0){$\langle$}}
\put(4,7){\makebox(0,0){$\scriptstyle q$}}
\put(20,7){\makebox(0,0){$\scriptstyle r$}}
\end{picture}$ 
of ${\mathrm{Sp}}(4,{\mathbb{R}})$ (meaning that its highest weight is $[q,r]$
with respect to the standard Bourbaki-ordered basis of fundamental weights).

In summary, there is a homogeneous contact geometry~$G/P$, with $G$
simple and $P$ parabolic, for which the BGG complex coincides with the Rumin
complex. 

\section{Pfaffian systems of rank three in five variables}\label{five}

Let $M$ be a $5$-manifold equipped with $H\subset TM$, a generic distribution 
of rank~$2$. Equivalently, let $I\subset\Lambda^1$ be a Pfaffian system of rank
$3$ that is generic in Cartan's sense, i.e.~the first derived system
$I^\prime$ has rank $2$ and the second derived system $I^{\prime\prime}$ is
zero. We have a filtration of the tangent bundle
$$H\subset [H,H]\subset TM\quad\mbox{by vector bundles of ranks}\quad 2,3,5$$
and a dual filtration of the cotangent bundle, which we shall write as
\begin{equation}\label{filter}\Lambda^1=\Lambda_H^1+L+I^\prime,\end{equation}
where $L$ is the line-bundle $I/I^\prime$. There are locally defined 
co-framings $(\omega^1,\omega^2,\omega^3,\omega^4,\omega^5)$ so that the
following congruences hold
\begin{equation}\label{G2structure}
\begin{array}{c}
d\omega^1\equiv\omega^3\wedge\omega^4\bmod\omega^1,\omega^2\qquad
d\omega^2\equiv\omega^3\wedge\omega^5\bmod\omega^1,\omega^2\\[6pt]
d\omega^3\equiv\omega^4\wedge\omega^5\bmod\omega^1,\omega^2,\omega^3
\end{array}\end{equation}
with $I^\prime={\mathrm{span}}\{\omega^1,\omega^2\}$ and
$I={\mathrm{span}}\{\omega^1,\omega^2,\omega^3\}$. We shall refer to such
co-framings as {\em adapted\/}. The Levi form for~$I$, defined as the
composition $I\to\Lambda^1\xrightarrow{\,d\,}\Lambda^2\to\Lambda_H^2$, has
$I^\prime$ as its kernel (by definition of $I^\prime$ or by viewing this form
in an adapted co-frame). Hence, the line-bundle $L$ may be canonically
identified with~$\Lambda_H^2$. Similarly, the Levi form for $I^\prime$, namely
$I^\prime\to\Lambda^1\xrightarrow{\,d\,}\Lambda^2\to
\Lambda_H^2+\Lambda_H^1\otimes L$, canonically identifies $I^\prime$ 
with $\Lambda_H^1\otimes L$. Therefore, we may rewrite (\ref{filter}) as
\begin{equation}\label{betterfilter}
\Lambda^1=\Lambda_H^1+\Lambda_H^2+\Lambda_H^1\otimes\Lambda_H^2.\end{equation}
To proceed, it is useful to have a more compact notation for the bundles
induced by $\Lambda_H^1$. Following a common convention for the irreducible
representations of ${\mathrm{GL}}(2,{\mathbb{C}})$, let us write
\begin{equation}\label{Clebsch}
(a,b)\in{\mathbb{Z}}^2\enskip\mbox{with $a\leq b$}\enskip
\mbox{for the bundle}
\enskip\textstyle\bigodot^{b-a}\!\Lambda_H^1\otimes(\Lambda_H^2)^a,
\end{equation}
where $\bigodot$ means symmetric tensor product. 
Then (\ref{betterfilter}) becomes
$$\Lambda^1=(0,1)+(1,1)+(1,2)$$
and the induced filtration on $2$-forms is
$$\Lambda^2=(1,1)+(1,2)+\begin{array}c(1,3)\\[-2pt] \oplus\\[-2pt] (2,2)
\end{array}+(2,3)+(3,3).$$ 
Without further ado, we may now consider the spectral sequence of the 
de~Rham complex filtered in this way. At the $E_0$-level we obtain
\footnotesize$$\begin{picture}(300,125)(0,5)
\put(151,70){\makebox(0,0){$\begin{array}{ccccccccccc}
\Lambda^0&\!(0,1)\!&\!(1,1)\!&0&0&0\\
&&\uparrow\\
&0&\!(1,1)\!&\!(1,2)\!&(2,2)&0&0\\
&&&\uparrow&\uparrow\\
&&0&\!(1,2)\!&\!(1,3)\!\oplus\!(2,2)\!&\!(2,3)\!&0&0\\
&&&&&\uparrow\\
&&&0&0&\!(2,3)\!&\!(2,4)\!\oplus\!(3,3)\!&\!(3,4)\!&0\\
&&&&&&\uparrow&\uparrow\\
&&&&0&0&(3,3)&\!(3,4)\!&\!(4,4)\!&0\\
&&&&&&&&\uparrow\\
&&&&&0&0&0&\!(4,4)\!&\!(4,5)\!&\!\Lambda^5\,.
\end{array}$}}
\put(0,0){\vector(1,0){40}}
\put(0,0){\vector(0,1){45}}
\put(30,5){$p$}
\put(3,35){$q$}
\end{picture}$$\normalsize
The $E_0$-level differentials are easily computed from the structure
equations~(\ref{G2structure}), the $E_1$-level is
\small$$\begin{picture}(300,80)(0,5)
\put(150,50){\makebox(0,0){$\begin{array}{ccccccccccccc}
\Lambda^0&\!\to\!&\!(0,1)\!&\;0\;&\;0\;&\;0\;&\;0\\[2pt]
&&0&0&0&0&0&0\\[2pt]
&&&0&0&\!(1,3)\!&0&0&0\\[2pt]
&&&&0&0&0&\!(2,4)\!&0&0\\[2pt]
&&&&&0&0&0&0&0&0\\[2pt]
&&&&&&\;0\;&\;0\;&\;0\;&\;0\;&\!(4,5)\!&\!\to\!&\!(5,5)\,,
\end{array}$}}
\put(0,0){\vector(1,0){40}}
\put(0,0){\vector(0,1){45}}
\put(30,5){$p$}
\put(3,35){$q$}
\end{picture}$$\normalsize
and we have shown that there is a differential complex
$$\Lambda^0\xrightarrow{\,\nabla\,}\Lambda_H^1
\xrightarrow{\,\nabla^{3}\,}(1,3)\xrightarrow{\,\nabla^{2}\,}(2,4)
\xrightarrow{\,\nabla^{3}\,}(4,5)\xrightarrow{\,\nabla\,}(5,5)$$
resolving the constant sheaf~${\mathbb{R}}$, where $\nabla^k$ simply indicates
a differential operator of order~$k$. If necessary, the structure equations
(\ref{G2structure}) can be used to compute the operators precisely. To compare
with the usual BGG complex, we follow Cartan~\cite{Ca} in realising the flat
model for this geometry as a homogeneous space $G/P$ where $G$ is the
exceptional non-compact Lie group $G_2$ and $P$ is a parabolic subgroup.
Specifically, following the notation of~\cite{BE}, the homogeneous space is
$\begin{picture}(24,5)
\put(5.6,0){\line(1,0){14.4}}
\put(4,1.6){\line(1,0){16}}
\put(5.6,3.2){\line(1,0){14.4}}
\put(4,1.5){\makebox(0,0){$\times$}}
\put(20,1.4){\makebox(0,0){$\bullet$}}
\put(12,1.5){\makebox(0,0){$\langle$}}
\end{picture}$.
The Levi factor of the parabolic subgroup is ${\mathrm{GL}}(2,{\mathbb{R}})$
but it is useful to identify its Lie algebra with ${\mathfrak{g}}_0$ where we
have graded the Lie algebra of $G_2$
$${\mathfrak{g}}={\mathfrak{g}}_{-3}\oplus
{\mathfrak{g}}_{-2}\oplus{\mathfrak{g}}_{-1}\oplus
\underbrace{{\mathfrak{g}}_0\oplus{\mathfrak{g}}_1
\oplus{\mathfrak{g}}_2\oplus{\mathfrak{g}}_3}_{\mathfrak{p}}$$
in accordance with the parabolic subalgebra~${\mathfrak{p}}$ (see~\cite{CSl}).
With these conventions, the cotangent bundle is
\begin{equation}\label{cotangent}\Lambda^1=
\begin{picture}(24,5)
\put(5.6,0){\line(1,0){14.4}}
\put(4,1.6){\line(1,0){16}}
\put(5.6,3.2){\line(1,0){14.4}}
\put(4,1.5){\makebox(0,0){$\times$}}
\put(20,1.4){\makebox(0,0){$\bullet$}}
\put(12,1.5){\makebox(0,0){$\langle$}}
\put(4,8){\makebox(0,0){$\scriptstyle -2$}}
\put(20,8){\makebox(0,0){$\scriptstyle 1$}}
\end{picture}+
\begin{picture}(24,5)
\put(5.6,0){\line(1,0){14.4}}
\put(4,1.6){\line(1,0){16}}
\put(5.6,3.2){\line(1,0){14.4}}
\put(4,1.5){\makebox(0,0){$\times$}}
\put(20,1.4){\makebox(0,0){$\bullet$}}
\put(12,1.5){\makebox(0,0){$\langle$}}
\put(4,8){\makebox(0,0){$\scriptstyle -1$}}
\put(20,8){\makebox(0,0){$\scriptstyle 0$}}
\end{picture}+
\begin{picture}(24,5)
\put(5.6,0){\line(1,0){14.4}}
\put(4,1.6){\line(1,0){16}}
\put(5.6,3.2){\line(1,0){14.4}}
\put(4,1.5){\makebox(0,0){$\times$}}
\put(20,1.4){\makebox(0,0){$\bullet$}}
\put(12,1.5){\makebox(0,0){$\langle$}}
\put(4,8){\makebox(0,0){$\scriptstyle -3$}}
\put(20,8){\makebox(0,0){$\scriptstyle 1$}}
\end{picture}\end{equation}
and the BGG complex is
\begin{equation}\label{G2BGG}\begin{picture}(24,5)
\put(5.6,0){\line(1,0){14.4}}
\put(4,1.6){\line(1,0){16}}
\put(5.6,3.2){\line(1,0){14.4}}
\put(4,1.5){\makebox(0,0){$\times$}}
\put(20,1.4){\makebox(0,0){$\bullet$}}
\put(12,1.5){\makebox(0,0){$\langle$}}
\put(4,8){\makebox(0,0){$\scriptstyle 0$}}
\put(20,8){\makebox(0,0){$\scriptstyle 0$}}
\end{picture}\xrightarrow{\,\nabla\,}
\begin{picture}(24,5)
\put(5.6,0){\line(1,0){14.4}}
\put(4,1.6){\line(1,0){16}}
\put(5.6,3.2){\line(1,0){14.4}}
\put(4,1.5){\makebox(0,0){$\times$}}
\put(20,1.4){\makebox(0,0){$\bullet$}}
\put(12,1.5){\makebox(0,0){$\langle$}}
\put(4,8){\makebox(0,0){$\scriptstyle -2$}}
\put(20,8){\makebox(0,0){$\scriptstyle 1$}}
\end{picture}\xrightarrow{\,\nabla^3\,}
\begin{picture}(24,5)
\put(5.6,0){\line(1,0){14.4}}
\put(4,1.6){\line(1,0){16}}
\put(5.6,3.2){\line(1,0){14.4}}
\put(4,1.5){\makebox(0,0){$\times$}}
\put(20,1.4){\makebox(0,0){$\bullet$}}
\put(12,1.5){\makebox(0,0){$\langle$}}
\put(4,8){\makebox(0,0){$\scriptstyle -5$}}
\put(20,8){\makebox(0,0){$\scriptstyle 2$}}
\end{picture}\xrightarrow{\,\nabla^2\,}
\begin{picture}(24,5)
\put(5.6,0){\line(1,0){14.4}}
\put(4,1.6){\line(1,0){16}}
\put(5.6,3.2){\line(1,0){14.4}}
\put(4,1.5){\makebox(0,0){$\times$}}
\put(20,1.4){\makebox(0,0){$\bullet$}}
\put(12,1.5){\makebox(0,0){$\langle$}}
\put(4,8){\makebox(0,0){$\scriptstyle -6$}}
\put(20,8){\makebox(0,0){$\scriptstyle 2$}}
\end{picture}\xrightarrow{\,\nabla^3\,}
\begin{picture}(24,5)
\put(5.6,0){\line(1,0){14.4}}
\put(4,1.6){\line(1,0){16}}
\put(5.6,3.2){\line(1,0){14.4}}
\put(4,1.5){\makebox(0,0){$\times$}}
\put(20,1.4){\makebox(0,0){$\bullet$}}
\put(12,1.5){\makebox(0,0){$\langle$}}
\put(4,8){\makebox(0,0){$\scriptstyle -6$}}
\put(20,8){\makebox(0,0){$\scriptstyle 1$}}
\end{picture}\xrightarrow{\,\nabla\,}
\begin{picture}(24,5)
\put(5.6,0){\line(1,0){14.4}}
\put(4,1.6){\line(1,0){16}}
\put(5.6,3.2){\line(1,0){14.4}}
\put(4,1.5){\makebox(0,0){$\times$}}
\put(20,1.4){\makebox(0,0){$\bullet$}}
\put(12,1.5){\makebox(0,0){$\langle$}}
\put(4,8){\makebox(0,0){$\scriptstyle -5$}}
\put(20,8){\makebox(0,0){$\scriptstyle 0$}}
\end{picture}.\end{equation}
More generally, in Dynkin diagram notation the bundle $(a,b)$ is written as 
\enskip$\begin{picture}(34,12)
\put(5.6,0){\line(1,0){24.4}}
\put(4,1.6){\line(1,0){26}}
\put(5.6,3.2){\line(1,0){24.4}}
\put(4,1.5){\makebox(0,0){$\times$}}
\put(30,1.4){\makebox(0,0){$\bullet$}}
\put(17,1.5){\makebox(0,0){$\langle$}}
\put(4,8){\makebox(0,0){$\scriptstyle a-2b$}}
\put(30,8){\makebox(0,0){$\scriptstyle b-a$}}
\end{picture}$.

In fact, there are several other complexes that can be created from the de~Rham
complex by choosing to carry out only some of the diagram chasing involved in
creating the BGG complex. We now explain two of these complexes and their
motivation. Keeping the Dynkin diagram notation, the filtration of the
$2$-forms induced from (\ref{cotangent}) is
$$\Lambda^2=
\begin{picture}(24,5)
\put(5.6,0){\line(1,0){14.4}}
\put(4,1.6){\line(1,0){16}}
\put(5.6,3.2){\line(1,0){14.4}}
\put(4,1.5){\makebox(0,0){$\times$}}
\put(20,1.4){\makebox(0,0){$\bullet$}}
\put(12,1.5){\makebox(0,0){$\langle$}}
\put(4,8){\makebox(0,0){$\scriptstyle -1$}}
\put(20,8){\makebox(0,0){$\scriptstyle 0$}}
\end{picture}+
\begin{picture}(24,5)
\put(5.6,0){\line(1,0){14.4}}
\put(4,1.6){\line(1,0){16}}
\put(5.6,3.2){\line(1,0){14.4}}
\put(4,1.5){\makebox(0,0){$\times$}}
\put(20,1.4){\makebox(0,0){$\bullet$}}
\put(12,1.5){\makebox(0,0){$\langle$}}
\put(4,8){\makebox(0,0){$\scriptstyle -3$}}
\put(20,8){\makebox(0,0){$\scriptstyle 1$}}
\end{picture}+
\begin{array}c\begin{picture}(24,5)
\put(5.6,0){\line(1,0){14.4}}
\put(4,1.6){\line(1,0){16}}
\put(5.6,3.2){\line(1,0){14.4}}
\put(4,1.5){\makebox(0,0){$\times$}}
\put(20,1.4){\makebox(0,0){$\bullet$}}
\put(12,1.5){\makebox(0,0){$\langle$}}
\put(4,8){\makebox(0,0){$\scriptstyle -5$}}
\put(20,8){\makebox(0,0){$\scriptstyle 2$}}
\end{picture}\\[0pt] \oplus\\[0pt]
\begin{picture}(24,5)
\put(5.6,0){\line(1,0){14.4}}
\put(4,1.6){\line(1,0){16}}
\put(5.6,3.2){\line(1,0){14.4}}
\put(4,1.5){\makebox(0,0){$\times$}}
\put(20,1.4){\makebox(0,0){$\bullet$}}
\put(12,1.5){\makebox(0,0){$\langle$}}
\put(4,8){\makebox(0,0){$\scriptstyle -2$}}
\put(20,8){\makebox(0,0){$\scriptstyle 0$}}
\end{picture}\end{array}+
\begin{picture}(24,5)
\put(5.6,0){\line(1,0){14.4}}
\put(4,1.6){\line(1,0){16}}
\put(5.6,3.2){\line(1,0){14.4}}
\put(4,1.5){\makebox(0,0){$\times$}}
\put(20,1.4){\makebox(0,0){$\bullet$}}
\put(12,1.5){\makebox(0,0){$\langle$}}
\put(4,8){\makebox(0,0){$\scriptstyle -4$}}
\put(20,8){\makebox(0,0){$\scriptstyle 1$}}
\end{picture}+
\begin{picture}(24,5)
\put(5.6,0){\line(1,0){14.4}}
\put(4,1.6){\line(1,0){16}}
\put(5.6,3.2){\line(1,0){14.4}}
\put(4,1.5){\makebox(0,0){$\times$}}
\put(20,1.4){\makebox(0,0){$\bullet$}}
\put(12,1.5){\makebox(0,0){$\langle$}}
\put(4,8){\makebox(0,0){$\scriptstyle -3$}}
\put(20,8){\makebox(0,0){$\scriptstyle 0$}}
\end{picture},$$
which suggests that one might cancel 
$\begin{picture}(24,5)
\put(5.6,0){\line(1,0){14.4}}
\put(4,1.6){\line(1,0){16}}
\put(5.6,3.2){\line(1,0){14.4}}
\put(4,1.5){\makebox(0,0){$\times$}}
\put(20,1.4){\makebox(0,0){$\bullet$}}
\put(12,1.5){\makebox(0,0){$\langle$}}
\put(4,8){\makebox(0,0){$\scriptstyle -1$}}
\put(20,8){\makebox(0,0){$\scriptstyle 0$}}
\end{picture}+
\begin{picture}(24,5)
\put(5.6,0){\line(1,0){14.4}}
\put(4,1.6){\line(1,0){16}}
\put(5.6,3.2){\line(1,0){14.4}}
\put(4,1.5){\makebox(0,0){$\times$}}
\put(20,1.4){\makebox(0,0){$\bullet$}}
\put(12,1.5){\makebox(0,0){$\langle$}}
\put(4,8){\makebox(0,0){$\scriptstyle -3$}}
\put(20,8){\makebox(0,0){$\scriptstyle 1$}}
\end{picture}$
from the exterior derivative $\Lambda^1\xrightarrow{\,d\,}\Lambda^2$. But, as a
sub-bundle of~$\Lambda^1$, this is precisely the original Pfaffian system~$I$.
So, we are trying to cancel from $\Lambda^1\to\Lambda^2$, the homomorphism
defined as the composition
$$I\hookrightarrow\Lambda^1\xrightarrow{\,d\,}\Lambda^2\to
\frac{\Lambda^2}{\Lambda^1\wedge I^\prime}.$$
We can accomplish this explicitly by means of an adapted co-frame.
Specifically, we define a differential operator
\begin{equation}\label{euler}{\mathcal{E}}:\Lambda_H^1=\Lambda^1/I\to
\Lambda^1\wedge I^\prime\subset\Lambda^2\end{equation}
by the following steps. Starting with any $1$-form~$\omega$,
\begin{itemize}
\item define $f$ by 
$d\omega \equiv f\,\omega^4\wedge\omega^5\bmod \omega^1,\omega^2,\omega^3$,
\item define $g,h$ by 
$d(\omega-f\omega^3)\equiv g\,\omega^3\wedge\omega^4+h\,\omega^3\wedge\omega^5
\bmod \omega^1,\omega^2$.
\end{itemize}
This is possible according to the structure equations (\ref{G2structure}),
which also imply that
$$d(\omega-f\omega^3-g\omega^1-h\omega^2)\equiv 0\bmod\omega^1,\omega^2,$$
in other words that
$${\mathcal{E}}\omega\equiv d(\omega-f\omega^3-g\omega^1-h\omega^2)
\in\Lambda^1\wedge I^\prime\subset\Lambda^2.$$
One checks easily that this definition of ${\mathcal{E}}\omega$ is independent
of choice of adapted co-framing. Moreover, if $\omega$ is actually a section of
$I$, say $\omega=F\omega^3+G\omega^1+H\omega^2$, then $f=F$, $g=G$, and $h=H$,
whence ${\mathcal{E}}\omega=0$. In other words, the differential operator
${\mathcal{E}}$ descends to $\Lambda_H^1$, as claimed in~(\ref{euler}).
\begin{theorem}\label{thm1} The sequence 
$$0\to{\mathbb{R}}\to\Lambda^0\xrightarrow{\,d_H\,}\Lambda_H^1
\xrightarrow{\,{\mathcal{E}}\,}\Lambda^1\wedge I^\prime\xrightarrow{\,d\,}
\Lambda^3\xrightarrow{\,d\,}\Lambda^4\xrightarrow{\,d\,}\Lambda^5\to 0$$
is a locally exact complex.  
\end{theorem}
\begin{proof}This is just a matter of unravelling definitions, bearing in mind
that the de~Rham complex is itself locally exact. Suppose, for example, that
$\omega$ is a $1$-form representing a section of $\Lambda_H^1$ that is
annihilated by~${\mathcal{E}}$. Locally, we need to find a smooth function
$\phi$ such that $\omega-d\phi$ is a section of~$I$. By construction of
${\mathcal{E}}$ we know $d(\omega-f\omega^3-g\omega^1-h\omega^2)=0$ for
some smooth functions $f,g,h$. Thus, by exactness of the de~Rham complex,
locally we can write $\omega-f\omega^3-g\omega^1-h\omega^2 =d\phi$ and then
$\omega-d\phi=f\omega^3+g\omega^1+h\omega^2$ is a section of $I$, as required.
The remaining verifications are similarly straightforward.
\end{proof}
The operator ${\mathcal{E}}:\Lambda_H^1\to\Lambda^1\wedge I^\prime$ has a
geometric meaning: a section $\phi$ of $\Lambda_H^1$ can be regarded as a
Lagrangian for an integral curve of~$I$. {From} this point of view
${\mathcal{E}}\phi$ are the Euler-Lagrange equations associated to this
Lagrangian. {From} its construction, one can easily verify that ${\mathcal{E}}$
is third order. More specifically, by construction, its symbol
$$\textstyle\bigodot^3\!\Lambda^1\otimes\Lambda_H^1\to
\Lambda^1\wedge I^\prime
=\begin{array}c\begin{picture}(24,5)
\put(5.6,0){\line(1,0){14.4}}
\put(4,1.6){\line(1,0){16}}
\put(5.6,3.2){\line(1,0){14.4}}
\put(4,1.5){\makebox(0,0){$\times$}}
\put(20,1.4){\makebox(0,0){$\bullet$}}
\put(12,1.5){\makebox(0,0){$\langle$}}
\put(4,8){\makebox(0,0){$\scriptstyle -5$}}
\put(20,8){\makebox(0,0){$\scriptstyle 2$}}
\end{picture}\\[0pt] \oplus\\[0pt]
\begin{picture}(24,5)
\put(5.6,0){\line(1,0){14.4}}
\put(4,1.6){\line(1,0){16}}
\put(5.6,3.2){\line(1,0){14.4}}
\put(4,1.5){\makebox(0,0){$\times$}}
\put(20,1.4){\makebox(0,0){$\bullet$}}
\put(12,1.5){\makebox(0,0){$\langle$}}
\put(4,8){\makebox(0,0){$\scriptstyle -2$}}
\put(20,8){\makebox(0,0){$\scriptstyle 0$}}
\end{picture}\end{array}+
\begin{picture}(24,5)
\put(5.6,0){\line(1,0){14.4}}
\put(4,1.6){\line(1,0){16}}
\put(5.6,3.2){\line(1,0){14.4}}
\put(4,1.5){\makebox(0,0){$\times$}}
\put(20,1.4){\makebox(0,0){$\bullet$}}
\put(12,1.5){\makebox(0,0){$\langle$}}
\put(4,8){\makebox(0,0){$\scriptstyle -4$}}
\put(20,8){\makebox(0,0){$\scriptstyle 1$}}
\end{picture}+
\begin{picture}(24,5)
\put(5.6,0){\line(1,0){14.4}}
\put(4,1.6){\line(1,0){16}}
\put(5.6,3.2){\line(1,0){14.4}}
\put(4,1.5){\makebox(0,0){$\times$}}
\put(20,1.4){\makebox(0,0){$\bullet$}}
\put(12,1.5){\makebox(0,0){$\langle$}}
\put(4,8){\makebox(0,0){$\scriptstyle -3$}}
\put(20,8){\makebox(0,0){$\scriptstyle 0$}}
\end{picture}$$
composes with the projection to 
$\begin{picture}(24,5)
\put(5.6,0){\line(1,0){14.4}}
\put(4,1.6){\line(1,0){16}}
\put(5.6,3.2){\line(1,0){14.4}}
\put(4,1.5){\makebox(0,0){$\times$}}
\put(20,1.4){\makebox(0,0){$\bullet$}}
\put(12,1.5){\makebox(0,0){$\langle$}}
\put(4,8){\makebox(0,0){$\scriptstyle -5$}}
\put(20,8){\makebox(0,0){$\scriptstyle 2$}}
\end{picture}$
as the homomorphism
$$\begin{array}{l}\bigodot^3\!\Lambda^1\otimes\Lambda_H^1=
\textstyle\bigodot^3\!\big(\,\begin{picture}(24,5)
\put(5.6,0){\line(1,0){14.4}}
\put(4,1.6){\line(1,0){16}}
\put(5.6,3.2){\line(1,0){14.4}}
\put(4,1.5){\makebox(0,0){$\times$}}
\put(20,1.4){\makebox(0,0){$\bullet$}}
\put(12,1.5){\makebox(0,0){$\langle$}}
\put(4,8){\makebox(0,0){$\scriptstyle -2$}}
\put(20,8){\makebox(0,0){$\scriptstyle 1$}}
\end{picture}+\begin{picture}(24,5)
\put(5.6,0){\line(1,0){14.4}}
\put(4,1.6){\line(1,0){16}}
\put(5.6,3.2){\line(1,0){14.4}}
\put(4,1.5){\makebox(0,0){$\times$}}
\put(20,1.4){\makebox(0,0){$\bullet$}}
\put(12,1.5){\makebox(0,0){$\langle$}}
\put(4,8){\makebox(0,0){$\scriptstyle -1$}}
\put(20,8){\makebox(0,0){$\scriptstyle 0$}}
\end{picture}+\begin{picture}(24,5)
\put(5.6,0){\line(1,0){14.4}}
\put(4,1.6){\line(1,0){16}}
\put(5.6,3.2){\line(1,0){14.4}}
\put(4,1.5){\makebox(0,0){$\times$}}
\put(20,1.4){\makebox(0,0){$\bullet$}}
\put(12,1.5){\makebox(0,0){$\langle$}}
\put(4,8){\makebox(0,0){$\scriptstyle -3$}}
\put(20,8){\makebox(0,0){$\scriptstyle 1$}}
\end{picture}\big)\otimes
\begin{picture}(24,5)
\put(5.6,0){\line(1,0){14.4}}
\put(4,1.6){\line(1,0){16}}
\put(5.6,3.2){\line(1,0){14.4}}
\put(4,1.5){\makebox(0,0){$\times$}}
\put(20,1.4){\makebox(0,0){$\bullet$}}
\put(12,1.5){\makebox(0,0){$\langle$}}
\put(4,8){\makebox(0,0){$\scriptstyle -2$}}
\put(20,8){\makebox(0,0){$\scriptstyle 1$}}
\end{picture}\\[5pt]
\qquad\to\bigodot^3\!\big(\,\begin{picture}(24,5)
\put(5.6,0){\line(1,0){14.4}}
\put(4,1.6){\line(1,0){16}}
\put(5.6,3.2){\line(1,0){14.4}}
\put(4,1.5){\makebox(0,0){$\times$}}
\put(20,1.4){\makebox(0,0){$\bullet$}}
\put(12,1.5){\makebox(0,0){$\langle$}}
\put(4,8){\makebox(0,0){$\scriptstyle -2$}}
\put(20,8){\makebox(0,0){$\scriptstyle 1$}}
\end{picture}\big)\otimes
\begin{picture}(24,5)
\put(5.6,0){\line(1,0){14.4}}
\put(4,1.6){\line(1,0){16}}
\put(5.6,3.2){\line(1,0){14.4}}
\put(4,1.5){\makebox(0,0){$\times$}}
\put(20,1.4){\makebox(0,0){$\bullet$}}
\put(12,1.5){\makebox(0,0){$\langle$}}
\put(4,8){\makebox(0,0){$\scriptstyle -2$}}
\put(20,8){\makebox(0,0){$\scriptstyle 1$}}
\end{picture}\to
\bigodot^2\!\big(\,\begin{picture}(24,5)
\put(5.6,0){\line(1,0){14.4}}
\put(4,1.6){\line(1,0){16}}
\put(5.6,3.2){\line(1,0){14.4}}
\put(4,1.5){\makebox(0,0){$\times$}}
\put(20,1.4){\makebox(0,0){$\bullet$}}
\put(12,1.5){\makebox(0,0){$\langle$}}
\put(4,8){\makebox(0,0){$\scriptstyle -2$}}
\put(20,8){\makebox(0,0){$\scriptstyle 1$}}
\end{picture}\big)\otimes
\Lambda^2\big(\,\begin{picture}(24,5)
\put(5.6,0){\line(1,0){14.4}}
\put(4,1.6){\line(1,0){16}}
\put(5.6,3.2){\line(1,0){14.4}}
\put(4,1.5){\makebox(0,0){$\times$}}
\put(20,1.4){\makebox(0,0){$\bullet$}}
\put(12,1.5){\makebox(0,0){$\langle$}}
\put(4,8){\makebox(0,0){$\scriptstyle -2$}}
\put(20,8){\makebox(0,0){$\scriptstyle 1$}}
\end{picture}\big)\\[5pt]
\mbox{ }\hspace{145pt}{}=\begin{picture}(24,5)
\put(5.6,0){\line(1,0){14.4}}
\put(4,1.6){\line(1,0){16}}
\put(5.6,3.2){\line(1,0){14.4}}
\put(4,1.5){\makebox(0,0){$\times$}}
\put(20,1.4){\makebox(0,0){$\bullet$}}
\put(12,1.5){\makebox(0,0){$\langle$}}
\put(4,8){\makebox(0,0){$\scriptstyle -4$}}
\put(20,8){\makebox(0,0){$\scriptstyle 2$}}
\end{picture}\otimes
\begin{picture}(24,5)
\put(5.6,0){\line(1,0){14.4}}
\put(4,1.6){\line(1,0){16}}
\put(5.6,3.2){\line(1,0){14.4}}
\put(4,1.5){\makebox(0,0){$\times$}}
\put(20,1.4){\makebox(0,0){$\bullet$}}
\put(12,1.5){\makebox(0,0){$\langle$}}
\put(4,8){\makebox(0,0){$\scriptstyle -1$}}
\put(20,8){\makebox(0,0){$\scriptstyle 0$}}
\end{picture}=
\begin{picture}(24,5)
\put(5.6,0){\line(1,0){14.4}}
\put(4,1.6){\line(1,0){16}}
\put(5.6,3.2){\line(1,0){14.4}}
\put(4,1.5){\makebox(0,0){$\times$}}
\put(20,1.4){\makebox(0,0){$\bullet$}}
\put(12,1.5){\makebox(0,0){$\langle$}}
\put(4,8){\makebox(0,0){$\scriptstyle -5$}}
\put(20,8){\makebox(0,0){$\scriptstyle 2$}}
\end{picture}.\end{array}$$
Furthermore, not only does the symbol have no component in 
$\begin{picture}(24,5)
\put(5.6,0){\line(1,0){14.4}}
\put(4,1.6){\line(1,0){16}}
\put(5.6,3.2){\line(1,0){14.4}}
\put(4,1.5){\makebox(0,0){$\times$}}
\put(20,1.4){\makebox(0,0){$\bullet$}}
\put(12,1.5){\makebox(0,0){$\langle$}}
\put(4,8){\makebox(0,0){$\scriptstyle -2$}}
\put(20,8){\makebox(0,0){$\scriptstyle 0$}}
\end{picture}$ but, in fact, the range of the operator ${\mathcal{E}}$ is 
entirely contained in the sub-bundle where this component vanishes. This is
easily seen in an adapted co-frame: since
$$d(\omega^1\wedge\omega^5-\omega^2\wedge\omega^4)=
2\,\omega^3\wedge\omega^4\wedge\omega^5\bmod\omega^1,\omega^2$$
and since $\omega^1\wedge\omega^5-\omega^2\wedge\omega^4$ spans 
$\begin{picture}(24,5)
\put(5.6,0){\line(1,0){14.4}}
\put(4,1.6){\line(1,0){16}}
\put(5.6,3.2){\line(1,0){14.4}}
\put(4,1.5){\makebox(0,0){$\times$}}
\put(20,1.4){\makebox(0,0){$\bullet$}}
\put(12,1.5){\makebox(0,0){$\langle$}}
\put(4,7){\makebox(0,0){$\scriptstyle -2$}}
\put(20,7){\makebox(0,0){$\scriptstyle 0$}}
\end{picture}$,
any exact $2$-form in $\Lambda^1\wedge I^\prime$ has vanishing component in
$\begin{picture}(24,5)
\put(5.6,0){\line(1,0){14.4}}
\put(4,1.6){\line(1,0){16}}
\put(5.6,3.2){\line(1,0){14.4}}
\put(4,1.5){\makebox(0,0){$\times$}}
\put(20,1.4){\makebox(0,0){$\bullet$}}
\put(12,1.5){\makebox(0,0){$\langle$}}
\put(4,7){\makebox(0,0){$\scriptstyle -2$}}
\put(20,7){\makebox(0,0){$\scriptstyle 0$}}
\end{picture}$. If we denote by $B^2$ the rank $6$ sub-bundle of 
$\Lambda^1\wedge I^\prime$ defined as the kernel of the natural projection
$\Lambda^1\wedge I^\prime\to
\begin{picture}(24,5)
\put(5.6,0){\line(1,0){14.4}}
\put(4,1.6){\line(1,0){16}}
\put(5.6,3.2){\line(1,0){14.4}}
\put(4,1.5){\makebox(0,0){$\times$}}
\put(20,1.4){\makebox(0,0){$\bullet$}}
\put(12,1.5){\makebox(0,0){$\langle$}}
\put(4,7){\makebox(0,0){$\scriptstyle -2$}}
\put(20,7){\makebox(0,0){$\scriptstyle 0$}}
\end{picture}$, and by $B^3$ the rank $9$ sub-bundle of $\Lambda^3$ generated 
by $\omega^1,\omega^2$, then we have cancelled 
$\begin{picture}(24,5)
\put(5.6,0){\line(1,0){14.4}}
\put(4,1.6){\line(1,0){16}}
\put(5.6,3.2){\line(1,0){14.4}}
\put(4,1.5){\makebox(0,0){$\times$}}
\put(20,1.4){\makebox(0,0){$\bullet$}}
\put(12,1.5){\makebox(0,0){$\langle$}}
\put(4,7){\makebox(0,0){$\scriptstyle -2$}}
\put(20,7){\makebox(0,0){$\scriptstyle 0$}}
\end{picture}$ from the complex of Theorem~\ref{thm1} and demonstrated the 
following improvement.
\begin{theorem}
The sequence 
\begin{equation}\label{basiccomplex}
0\to{\mathbb{R}}\to\Lambda^0\xrightarrow{\,d_H\,}\Lambda_H^1
\xrightarrow{\,{\mathcal{E}}\,}B^2\xrightarrow{\,d\,}
B^3\xrightarrow{\,d\,}\Lambda^4\xrightarrow{\,d\,}\Lambda^5\to 0\end{equation}
is a locally exact complex.  
\end{theorem}
The ranks of the bundles and the orders of the differential operators in
(\ref{basiccomplex}) are
$$1\xrightarrow{\,\nabla\,}2\xrightarrow{\,\nabla^3\,}6
\xrightarrow{\,\nabla\,}9\xrightarrow{\,\nabla\,}5\xrightarrow{\,\nabla\,}1$$
but if we consider $d:B^2\to B^3$ in more detail
$$B^2=\begin{picture}(24,5)
\put(5.6,0){\line(1,0){14.4}}
\put(4,1.6){\line(1,0){16}}
\put(5.6,3.2){\line(1,0){14.4}}
\put(4,1.5){\makebox(0,0){$\times$}}
\put(20,1.4){\makebox(0,0){$\bullet$}}
\put(12,1.5){\makebox(0,0){$\langle$}}
\put(4,8){\makebox(0,0){$\scriptstyle -5$}}
\put(20,8){\makebox(0,0){$\scriptstyle 2$}}
\end{picture}+\begin{picture}(24,5)
\put(5.6,0){\line(1,0){14.4}}
\put(4,1.6){\line(1,0){16}}
\put(5.6,3.2){\line(1,0){14.4}}
\put(4,1.5){\makebox(0,0){$\times$}}
\put(20,1.4){\makebox(0,0){$\bullet$}}
\put(12,1.5){\makebox(0,0){$\langle$}}
\put(4,8){\makebox(0,0){$\scriptstyle -4$}}
\put(20,8){\makebox(0,0){$\scriptstyle 1$}}
\end{picture}+\begin{picture}(24,5)
\put(5.6,0){\line(1,0){14.4}}
\put(4,1.6){\line(1,0){16}}
\put(5.6,3.2){\line(1,0){14.4}}
\put(4,1.5){\makebox(0,0){$\times$}}
\put(20,1.4){\makebox(0,0){$\bullet$}}
\put(12,1.5){\makebox(0,0){$\langle$}}
\put(4,8){\makebox(0,0){$\scriptstyle -3$}}
\put(20,8){\makebox(0,0){$\scriptstyle 0$}}
\end{picture}\rightarrow
\begin{picture}(24,5)
\put(5.6,0){\line(1,0){14.4}}
\put(4,1.6){\line(1,0){16}}
\put(5.6,3.2){\line(1,0){14.4}}
\put(4,1.5){\makebox(0,0){$\times$}}
\put(20,1.4){\makebox(0,0){$\bullet$}}
\put(12,1.5){\makebox(0,0){$\langle$}}
\put(4,8){\makebox(0,0){$\scriptstyle -4$}}
\put(20,8){\makebox(0,0){$\scriptstyle 1$}}
\end{picture}+\!\!
\begin{array}c\begin{picture}(24,5)
\put(5.6,0){\line(1,0){14.4}}
\put(4,1.6){\line(1,0){16}}
\put(5.6,3.2){\line(1,0){14.4}}
\put(4,1.5){\makebox(0,0){$\times$}}
\put(20,1.4){\makebox(0,0){$\bullet$}}
\put(12,1.5){\makebox(0,0){$\langle$}}
\put(4,8){\makebox(0,0){$\scriptstyle -6$}}
\put(20,8){\makebox(0,0){$\scriptstyle 2$}}
\end{picture}\\[0pt] \oplus\\[0pt]
\begin{picture}(24,5)
\put(5.6,0){\line(1,0){14.4}}
\put(4,1.6){\line(1,0){16}}
\put(5.6,3.2){\line(1,0){14.4}}
\put(4,1.5){\makebox(0,0){$\times$}}
\put(20,1.4){\makebox(0,0){$\bullet$}}
\put(12,1.5){\makebox(0,0){$\langle$}}
\put(4,8){\makebox(0,0){$\scriptstyle -3$}}
\put(20,8){\makebox(0,0){$\scriptstyle 0$}}
\end{picture}\end{array}
\!\!+\begin{picture}(24,5)
\put(5.6,0){\line(1,0){14.4}}
\put(4,1.6){\line(1,0){16}}
\put(5.6,3.2){\line(1,0){14.4}}
\put(4,1.5){\makebox(0,0){$\times$}}
\put(20,1.4){\makebox(0,0){$\bullet$}}
\put(12,1.5){\makebox(0,0){$\langle$}}
\put(4,8){\makebox(0,0){$\scriptstyle -5$}}
\put(20,8){\makebox(0,0){$\scriptstyle 1$}}
\end{picture}+
\begin{picture}(24,5)
\put(5.6,0){\line(1,0){14.4}}
\put(4,1.6){\line(1,0){16}}
\put(5.6,3.2){\line(1,0){14.4}}
\put(4,1.5){\makebox(0,0){$\times$}}
\put(20,1.4){\makebox(0,0){$\bullet$}}
\put(12,1.5){\makebox(0,0){$\langle$}}
\put(4,8){\makebox(0,0){$\scriptstyle -4$}}
\put(20,8){\makebox(0,0){$\scriptstyle 0$}}
\end{picture}=B^3,$$
then it suggests that we should be able to eliminate
$\begin{picture}(24,5)
\put(5.6,0){\line(1,0){14.4}}
\put(4,1.6){\line(1,0){16}}
\put(5.6,3.2){\line(1,0){14.4}}
\put(4,1.5){\makebox(0,0){$\times$}}
\put(20,1.4){\makebox(0,0){$\bullet$}}
\put(12,1.5){\makebox(0,0){$\langle$}}
\put(4,8){\makebox(0,0){$\scriptstyle -4$}}
\put(20,8){\makebox(0,0){$\scriptstyle 1$}}
\end{picture}+\begin{picture}(24,5)
\put(5.6,0){\line(1,0){14.4}}
\put(4,1.6){\line(1,0){16}}
\put(5.6,3.2){\line(1,0){14.4}}
\put(4,1.5){\makebox(0,0){$\times$}}
\put(20,1.4){\makebox(0,0){$\bullet$}}
\put(12,1.5){\makebox(0,0){$\langle$}}
\put(4,8){\makebox(0,0){$\scriptstyle -3$}}
\put(20,8){\makebox(0,0){$\scriptstyle 0$}}
\end{picture}$ from both bundles. This is, indeed, the case as can be seen in 
an adapted co-frame: writing a general section of $B^2$ as 
$$\omega=\mu\wedge\omega^1+\nu\wedge\omega^2\mbox{ s.t.\ }
d\omega\equiv 0\bmod\omega^1,\omega^2$$
we may define a differential operator ${\mathcal{F}}:B^2\to B^3$ by the
following familiar steps.
\begin{center}\begin{tabular}{ll}
$\bullet\!\!$& Define $f,g$ by $d\omega\equiv\begin{array}[t]{l} 
f\,\omega^1\wedge\omega^4\wedge\omega^5+
g\,\omega^2\wedge\omega^4\wedge\omega^5\\
\quad{}\bmod\omega^1\wedge\omega^2,\omega^1\wedge\omega^3,\omega^2\wedge\omega^3,
\end{array}$\\
$\bullet\!\!$& Define $h$ by 
$d(\omega-f\,\omega^1\wedge\omega^3-g\,\omega^2\wedge\omega^3)\equiv$\\
&\mbox{ }\hspace{120pt}$\begin{array}[t]{l}
h\,
(\omega^2\wedge\omega^3\wedge\omega^4-\omega^1\wedge\omega^3\wedge\omega^5)\\
{}+p\,\omega^1\wedge\omega^3\wedge\omega^4
+q\,\omega^2\wedge\omega^3\wedge\omega^5\\
\enskip{}+r\,
(\omega^2\wedge\omega^3\wedge\omega^4+\omega^1\wedge\omega^3\wedge\omega^5)\\
\quad{}\bmod\omega^1\wedge\omega^2.\end{array}$
\end{tabular}\end{center}
This is possible according to the structure equations (\ref{G2structure}),
which also imply that
$${\mathcal{F}}\omega\equiv
d(\omega-f\,\omega^1\wedge\omega^3-g\,\omega^2\wedge\omega^3
-h\,\omega^1\wedge\omega^2)$$
lies in the sub-bundle
$$C^3\equiv\begin{picture}(24,5)
\put(5.6,0){\line(1,0){14.4}}
\put(4,1.6){\line(1,0){16}}
\put(5.6,3.2){\line(1,0){14.4}}
\put(4,1.5){\makebox(0,0){$\times$}}
\put(20,1.4){\makebox(0,0){$\bullet$}}
\put(12,1.5){\makebox(0,0){$\langle$}}
\put(4,8){\makebox(0,0){$\scriptstyle -6$}}
\put(20,8){\makebox(0,0){$\scriptstyle 2$}}
\end{picture}+\begin{picture}(24,5)
\put(5.6,0){\line(1,0){14.4}}
\put(4,1.6){\line(1,0){16}}
\put(5.6,3.2){\line(1,0){14.4}}
\put(4,1.5){\makebox(0,0){$\times$}}
\put(20,1.4){\makebox(0,0){$\bullet$}}
\put(12,1.5){\makebox(0,0){$\langle$}}
\put(4,8){\makebox(0,0){$\scriptstyle -5$}}
\put(20,8){\makebox(0,0){$\scriptstyle 1$}}
\end{picture}+\begin{picture}(24,5)
\put(5.6,0){\line(1,0){14.4}}
\put(4,1.6){\line(1,0){16}}
\put(5.6,3.2){\line(1,0){14.4}}
\put(4,1.5){\makebox(0,0){$\times$}}
\put(20,1.4){\makebox(0,0){$\bullet$}}
\put(12,1.5){\makebox(0,0){$\langle$}}
\put(4,8){\makebox(0,0){$\scriptstyle -4$}}
\put(20,8){\makebox(0,0){$\scriptstyle 0$}}
\end{picture}$$
of $B^3\subset\Lambda^3$ and that it descends to the quotient
$$B^2=\begin{picture}(24,5)
\put(5.6,0){\line(1,0){14.4}}
\put(4,1.6){\line(1,0){16}}
\put(5.6,3.2){\line(1,0){14.4}}
\put(4,1.5){\makebox(0,0){$\times$}}
\put(20,1.4){\makebox(0,0){$\bullet$}}
\put(12,1.5){\makebox(0,0){$\langle$}}
\put(4,8){\makebox(0,0){$\scriptstyle -5$}}
\put(20,8){\makebox(0,0){$\scriptstyle 2$}}
\end{picture}+\begin{picture}(24,5)
\put(5.6,0){\line(1,0){14.4}}
\put(4,1.6){\line(1,0){16}}
\put(5.6,3.2){\line(1,0){14.4}}
\put(4,1.5){\makebox(0,0){$\times$}}
\put(20,1.4){\makebox(0,0){$\bullet$}}
\put(12,1.5){\makebox(0,0){$\langle$}}
\put(4,8){\makebox(0,0){$\scriptstyle -4$}}
\put(20,8){\makebox(0,0){$\scriptstyle 1$}}
\end{picture}+\begin{picture}(24,5)
\put(5.6,0){\line(1,0){14.4}}
\put(4,1.6){\line(1,0){16}}
\put(5.6,3.2){\line(1,0){14.4}}
\put(4,1.5){\makebox(0,0){$\times$}}
\put(20,1.4){\makebox(0,0){$\bullet$}}
\put(12,1.5){\makebox(0,0){$\langle$}}
\put(4,8){\makebox(0,0){$\scriptstyle -3$}}
\put(20,8){\makebox(0,0){$\scriptstyle 0$}}
\end{picture}\twoheadrightarrow
\begin{picture}(24,5)
\put(5.6,0){\line(1,0){14.4}}
\put(4,1.6){\line(1,0){16}}
\put(5.6,3.2){\line(1,0){14.4}}
\put(4,1.5){\makebox(0,0){$\times$}}
\put(20,1.4){\makebox(0,0){$\bullet$}}
\put(12,1.5){\makebox(0,0){$\langle$}}
\put(4,8){\makebox(0,0){$\scriptstyle -5$}}
\put(20,8){\makebox(0,0){$\scriptstyle 2$}}
\end{picture}\equiv C^2$$
of $B^2$. It is easily verified that this definition of ${\mathcal{F}}$ is
independent of choice of co-framing and that, if we denote by
$\bar{\mathcal{E}}$ the composition
$$\Lambda_H^1\xrightarrow{\,{\mathcal{E}}\,}B^2\to C^2,$$
then the expected theorem follows:
\begin{theorem}
The sequence 
\begin{equation}\label{secondcomplex}
0\to{\mathbb{R}}\to\Lambda^0\xrightarrow{\,d_H\,}\Lambda_H^1
\xrightarrow{\,\bar{\mathcal{E}}\,}C^2\xrightarrow{\,{\mathcal{F}}\,}
C^3\xrightarrow{\,d\,}\Lambda^4\xrightarrow{\,d\,}\Lambda^5\to 0\end{equation}
is a locally exact complex.  
\end{theorem}
The ranks of the bundles and the orders of the differential operators in
(\ref{secondcomplex}) are
$$1\xrightarrow{\,\nabla\,}2\xrightarrow{\,\nabla^3\,}3
\xrightarrow{\,\nabla^3\,}6\xrightarrow{\,\nabla\,}5
\xrightarrow{\,\nabla\,}1.$$
Writing (\ref{secondcomplex}) as
$$\begin{array}{l}\begin{picture}(24,5)
\put(5.6,0){\line(1,0){14.4}}
\put(4,1.6){\line(1,0){16}}
\put(5.6,3.2){\line(1,0){14.4}}
\put(4,1.5){\makebox(0,0){$\times$}}
\put(20,1.4){\makebox(0,0){$\bullet$}}
\put(12,1.5){\makebox(0,0){$\langle$}}
\put(4,8){\makebox(0,0){$\scriptstyle 0$}}
\put(20,8){\makebox(0,0){$\scriptstyle 0$}}
\end{picture}\to
\begin{picture}(24,5)
\put(5.6,0){\line(1,0){14.4}}
\put(4,1.6){\line(1,0){16}}
\put(5.6,3.2){\line(1,0){14.4}}
\put(4,1.5){\makebox(0,0){$\times$}}
\put(20,1.4){\makebox(0,0){$\bullet$}}
\put(12,1.5){\makebox(0,0){$\langle$}}
\put(4,8){\makebox(0,0){$\scriptstyle -2$}}
\put(20,8){\makebox(0,0){$\scriptstyle 1$}}
\end{picture}\to
\begin{picture}(24,5)
\put(5.6,0){\line(1,0){14.4}}
\put(4,1.6){\line(1,0){16}}
\put(5.6,3.2){\line(1,0){14.4}}
\put(4,1.5){\makebox(0,0){$\times$}}
\put(20,1.4){\makebox(0,0){$\bullet$}}
\put(12,1.5){\makebox(0,0){$\langle$}}
\put(4,8){\makebox(0,0){$\scriptstyle -5$}}
\put(20,8){\makebox(0,0){$\scriptstyle 2$}}
\end{picture}\to
\begin{picture}(24,5)
\put(5.6,0){\line(1,0){14.4}}
\put(4,1.6){\line(1,0){16}}
\put(5.6,3.2){\line(1,0){14.4}}
\put(4,1.5){\makebox(0,0){$\times$}}
\put(20,1.4){\makebox(0,0){$\bullet$}}
\put(12,1.5){\makebox(0,0){$\langle$}}
\put(4,8){\makebox(0,0){$\scriptstyle -6$}}
\put(20,8){\makebox(0,0){$\scriptstyle 2$}}
\end{picture}+\begin{picture}(24,5)
\put(5.6,0){\line(1,0){14.4}}
\put(4,1.6){\line(1,0){16}}
\put(5.6,3.2){\line(1,0){14.4}}
\put(4,1.5){\makebox(0,0){$\times$}}
\put(20,1.4){\makebox(0,0){$\bullet$}}
\put(12,1.5){\makebox(0,0){$\langle$}}
\put(4,8){\makebox(0,0){$\scriptstyle -5$}}
\put(20,8){\makebox(0,0){$\scriptstyle 1$}}
\end{picture}+\begin{picture}(24,5)
\put(5.6,0){\line(1,0){14.4}}
\put(4,1.6){\line(1,0){16}}
\put(5.6,3.2){\line(1,0){14.4}}
\put(4,1.5){\makebox(0,0){$\times$}}
\put(20,1.4){\makebox(0,0){$\bullet$}}
\put(12,1.5){\makebox(0,0){$\langle$}}
\put(4,8){\makebox(0,0){$\scriptstyle -4$}}
\put(20,8){\makebox(0,0){$\scriptstyle 0$}}
\end{picture}\to\\[5pt]
\mbox{ }\hspace{152pt}\begin{picture}(24,5)
\put(5.6,0){\line(1,0){14.4}}
\put(4,1.6){\line(1,0){16}}
\put(5.6,3.2){\line(1,0){14.4}}
\put(4,1.5){\makebox(0,0){$\times$}}
\put(20,1.4){\makebox(0,0){$\bullet$}}
\put(12,1.5){\makebox(0,0){$\langle$}}
\put(4,8){\makebox(0,0){$\scriptstyle -5$}}
\put(20,8){\makebox(0,0){$\scriptstyle 1$}}
\end{picture}+\begin{picture}(24,5)
\put(5.6,0){\line(1,0){14.4}}
\put(4,1.6){\line(1,0){16}}
\put(5.6,3.2){\line(1,0){14.4}}
\put(4,1.5){\makebox(0,0){$\times$}}
\put(20,1.4){\makebox(0,0){$\bullet$}}
\put(12,1.5){\makebox(0,0){$\langle$}}
\put(4,8){\makebox(0,0){$\scriptstyle -4$}}
\put(20,8){\makebox(0,0){$\scriptstyle 0$}}
\end{picture}+\begin{picture}(24,5)
\put(5.6,0){\line(1,0){14.4}}
\put(4,1.6){\line(1,0){16}}
\put(5.6,3.2){\line(1,0){14.4}}
\put(4,1.5){\makebox(0,0){$\times$}}
\put(20,1.4){\makebox(0,0){$\bullet$}}
\put(12,1.5){\makebox(0,0){$\langle$}}
\put(4,8){\makebox(0,0){$\scriptstyle -6$}}
\put(20,8){\makebox(0,0){$\scriptstyle 1$}}
\end{picture}\to
\begin{picture}(24,5)
\put(5.6,0){\line(1,0){14.4}}
\put(4,1.6){\line(1,0){16}}
\put(5.6,3.2){\line(1,0){14.4}}
\put(4,1.5){\makebox(0,0){$\times$}}
\put(20,1.4){\makebox(0,0){$\bullet$}}
\put(12,1.5){\makebox(0,0){$\langle$}}
\put(4,8){\makebox(0,0){$\scriptstyle -5$}}
\put(20,8){\makebox(0,0){$\scriptstyle 0$}}
\end{picture},\end{array}$$
suggests one final cancellation, specifically of 
$\begin{picture}(24,5)
\put(5.6,0){\line(1,0){14.4}}
\put(4,1.6){\line(1,0){16}}
\put(5.6,3.2){\line(1,0){14.4}}
\put(4,1.5){\makebox(0,0){$\times$}}
\put(20,1.4){\makebox(0,0){$\bullet$}}
\put(12,1.5){\makebox(0,0){$\langle$}}
\put(4,8){\makebox(0,0){$\scriptstyle -5$}}
\put(20,8){\makebox(0,0){$\scriptstyle 1$}}
\end{picture}+\begin{picture}(24,5)
\put(5.6,0){\line(1,0){14.4}}
\put(4,1.6){\line(1,0){16}}
\put(5.6,3.2){\line(1,0){14.4}}
\put(4,1.5){\makebox(0,0){$\times$}}
\put(20,1.4){\makebox(0,0){$\bullet$}}
\put(12,1.5){\makebox(0,0){$\langle$}}
\put(4,8){\makebox(0,0){$\scriptstyle -4$}}
\put(20,8){\makebox(0,0){$\scriptstyle 0$}}
\end{picture}$ 
from $C^3$ and $\Lambda^4$. The reader can readily verify that this gives the
BGG complex~(\ref{G2BGG}). It is interesting to note that the ranks of the
bundles and orders of differential operators in the BGG complex are
$$1\xrightarrow{\,\nabla\,}2\xrightarrow{\,\nabla^3\,}3
\xrightarrow{\,\nabla^2\,}3\xrightarrow{\,\nabla^3\,}2
\xrightarrow{\,\nabla\,}1.$$
In particular, the order of the differential operator in the middle has gone
down from $3$ to~$2$. Since our filtering on the de~Rham complex is, by
construction, compatible with the tautological Hodge isomorphisms
$\Lambda^p=\Lambda^5\otimes(\Lambda^{5-p})^*$, and since we have run the
spectral sequence to its end, it follows that the BGG complex is formally
self-adjoint.

\section{Pfaffian systems of rank three in six variables}\label{threeinsix}
Let $M$ be a $6$-manifold equipped with $H\subset TM$, a generic distribution of
rank~$3$. Equivalently, let $I\subset\Lambda^1$ be a Pfaffian system of rank
$3$ that is generic in Cartan's sense, i.e.~the first derived system
$I^\prime$ is zero. Locally there are co-framings 
$(\omega^1,\omega^2,\omega^3,\omega^4,\omega^5,\omega^6)$ so that 
$\omega^1,\omega^2,\omega^3$ span $I$ and the following congruences hold.
\begin{equation}\label{Spin7structure}
\begin{array}{c}
d\omega^1\equiv\omega^5\wedge\omega^6\bmod\omega^1,\omega^2,\omega^3\quad
d\omega^2\equiv\omega^6\wedge\omega^4\bmod\omega^1,\omega^2,\omega^3\\[3pt]
d\omega^3\equiv\omega^4\wedge\omega^5\bmod\omega^1,\omega^2,\omega^3.
\end{array}\end{equation}

In the terminology of~\cite{B}, these co-framings are {\em $1$-adapted\/}. As
usual, let us write $\Lambda_H^1$ for $\Lambda^1/I$. Then the Levi form
${\mathcal{L}}:I\to\Lambda_H^2$ defined as the composition
$I\hookrightarrow\Lambda^1\xrightarrow{\,d\,}
\Lambda^2\twoheadrightarrow\Lambda_H^2$ is an isomorphism and we can
canonically identify $I$ with $\Lambda_H^2$ as vector bundles. Indeed, this
isomorphism is apparent in our $1$-adapted co-framing (\ref{Spin7structure}).
We may mimic (\ref{Clebsch}) to write, up to isomorphism, the general
Schur-irreducible bundle induced by $\Lambda_H^1$ as
$(a,b,c)\in{\mathbb{Z}}^3$ with $a\leq b\leq c$ for the bundle
$$\textstyle
(\bigodot^{c-b}\!\Lambda_H^1\otimes\bigodot^{b-a}(\Lambda_H^1)^*)_\circ
\otimes(\Lambda_H^3)^b,$$
where $\circ$ as a subscript means to take the trace-free part.
These observations mean that we may write the filtration
$$\Lambda^1=\Lambda_H^1+I\quad\mbox{as}\quad\Lambda^1=(0,0,1)+(0,1,1)$$
and decompose the induced filtrations on the higher forms as
\begin{equation}\label{higherforms}\begin{array}{l}
\Lambda^2=
(0,1,1)+\begin{array}c(0,1,2)\\[-2pt] \oplus\\[-2pt] (1,1,1)\end{array}
+(1,1,2)\\[20pt]
\Lambda^3=
(1,1,1)+\begin{array}c(0,2,2)\\[-2pt] \oplus\\[-2pt] (1,1,2)\end{array}+
\begin{array}c(1,1,3)\\[-2pt] \oplus\\[-2pt] (1,2,2)\end{array}+(2,2,2)\\[20pt]
\Lambda^4=
(1,2,2)+\begin{array}c(1,2,3)\\[-2pt] \oplus\\[-2pt] (2,2,2)\end{array}+
(2,2,3)\enskip\quad
\Lambda^5=(2,2,3)+(2,3,3)\\[12pt]
\mbox{ }\hspace{140pt}\Lambda^6=(3,3,3).
\end{array}\end{equation}
{From} the structure equations (\ref{Spin7structure}) for a $1$-adapted
co-framing it is easily verified that all expected cancellations at the
$E_0$-level of the associated spectral sequence actually take place and we have
proved the following result.
\begin{theorem}\label{preBGG}
There is a canonically defined locally exact differential complex
$$\begin{array}{r}
0\to{\mathbb{R}}\to(0,0,0)\to(0,0,1)\to(0,1,2)\to
\big[(0,2,2)+(1,1,3)\big]\to\quad\\
(1,2,3)\to(2,3,3)\to(3,3,3)\to 0
\end{array}$$
on any smooth $6$-manifold equipped with a generic $3$-distribution.
\end{theorem}
For the moment, the bundle $\big[(0,2,2)+(1,1,3)\big]$ is a canonically defined
sub-quotient of $\Lambda^3$ but, in fact, one can improve matters as the
following theorem shows (the analogous step was not necessary in~\S\ref{five}).
\begin{theorem}\label{preferred} A splitting of the short exact sequence 
$$0\to(0,1,1)\to\Lambda^1\to(0,0,1)\to 0$$
gives rise to a homomorphism of vector bundles defined as the composition
\begin{equation}\label{hom}
(1,1,1)\to\Lambda^2\xrightarrow{\,d\,}\Lambda^3\to(0,2,2)\end{equation}
and there is a preferred class of splittings characterised by requiring that 
this induced homomorphism vanish. This preference canonically splits 
the bundle $\big[(0,2,2)+(1,1,3)\big]$.
\end{theorem}

\begin{proof} Certainly, a splitting of the $1$-forms splits all the other
forms and so, from~(\ref{higherforms}), one may consider the composition
(\ref{hom}) obtained by splitting the $2$-forms and $3$-forms. To see that it
is a homomorphism, rather than the differential operator it might appear to be,
notice that if $\Omega$ is a $2$-form in $(1,1,1)$, then 
\begin{equation}\label{Leibniz}d(f\Omega)=fd\Omega+df\wedge\Omega\end{equation}
and it is clear that $df\wedge\Omega$ has components only in 
$$\Lambda^1\otimes (1,1,1)=\big((0,0,1)\oplus(0,1,1)\big)\otimes(1,1,1)
=(1,1,2)\oplus(1,2,2)$$
inside
$$\Lambda^3=(1,1,1)\oplus
\begin{array}c(0,2,2)\\[-2pt] \oplus\\[-2pt] (1,1,2)\end{array}
\oplus\begin{array}c(1,1,3)\\[-2pt] \oplus\\[-2pt] (1,2,2)\end{array}
\oplus(2,2,2).$$
In particular, if we project to $(0,2,2)$ as in~(\ref{hom}), then
$df\wedge\Omega$ does not contribute and, from~(\ref{Leibniz}), the result is
linear over the functions. Now suppose we change the splitting of the
$1$-forms. The freedom in doing so lies in
\begin{equation}\label{basicfreedom}
{\mathrm{Hom}}(\big(0,0,1),(0,1,1)\big)=(-1,1,1)\oplus(0,0,1).\end{equation}
This same freedom shows up in splitting the first part of~$\Lambda^3$:
$$\Lambda^3=(1,1,1)+
\begin{array}c(0,2,2)\\[-2pt] \oplus\\[-2pt] (1,1,2)\end{array}+\cdots{}$$
and
$${\mathrm{Hom}}\!\left(\!(1,1,1),
\begin{array}c(0,2,2)\\[-2pt] \oplus\\[-2pt] (1,1,2)\end{array}\right)=
\begin{array}c(-1,1,1)\\[-2pt] \oplus\\[-2pt] (0,0,1)\end{array}.$$
Bearing in mind that the composition
$(1,1,1)\to\Lambda^2\to\Lambda^3\to(1,1,1)$ is an isomorphism, independent of
choice of splitting (it is responsible for one of the cancellations occurring
at the $E_0$-level of the spectral sequence), we conclude that we can spend the
$(-1,1,1)$-freedom in splitting precisely in setting the homomorphism
(\ref{hom}) to zero. Now let us consider how this impacts on the sub-quotient
$\big[(0,2,2)+(1,1,3)\big]$ of $\Lambda^3$. The freedom in splitting this
sub-quotient lies in
\begin{equation}\label{morefreedom}
{\mathrm{Hom}}\big((0,2,2),(1,1,3)\big)=(-1,-1,3)\oplus(-1,0,2)\oplus(-1,1,1)
\end{equation}
and one sees that the only way that (\ref{basicfreedom}) can enter is through
$(-1,1,1)$. Having eliminated this freedom by a preferred choice of splittings
for~$\Lambda^1$, it is thereby eliminated from (\ref{morefreedom}) and we have
obtained our canonical splitting.
\end{proof}
The preferred splittings of $\Lambda^1$ afforded by Theorem~\ref{preferred} can
be conveniently expressed in terms of our $1$-adapted co-framings
satisfying~(\ref{Spin7structure}). If such a co-framing is used to split
$\Lambda^1$, then the resulting sub-bundle $(1,1,1)$ of $\Lambda^2$ is spanned
by $\Omega\equiv\omega^1\wedge\omega^4+\omega^2\wedge\omega^5
+\omega^3\wedge\omega^6$ and, following through its proof, the preferred
splittings of Theorem~\ref{preferred} are characterised by requiring that
$$d\Omega\equiv3\,\omega^4\wedge\omega^5\wedge\omega^6+
\omega\wedge\Omega{\scriptstyle\bmod\omega^1\wedge\omega^2,
\omega^2\wedge\omega^3,\omega^3\wedge\omega^1},
\enskip\mbox{ for some $1$-form }\omega.$$
In the terminology of~\cite{B}, co-framings satisfying this extra congruence
are called {\em $2$-adapted\/}.

The Lie algebra ${\mathfrak{so}}(4,3)$ admits a grading of the form
\begin{equation}\label{grading}
\begin{array}{l}{\mathfrak{g}}_{-2}\oplus{\mathfrak{g}}_{-1}\oplus
\underbrace{{\mathfrak{g}}_0
\oplus{\mathfrak{g}}_1\oplus{\mathfrak{g}}_2}_{\mathfrak{p}}=\\
\enskip\begin{picture}(42,5)
\put(20,.1){\line(1,0){14.3}}
\put(20,3.2){\line(1,0){14.3}}
\put(20,1.5){\makebox(0,0){$\bullet$}}
\put(36,1.6){\makebox(0,0){$\times$}}
\put(28,1.5){\makebox(0,0){$\rangle$}}
\put(4,1.5){\makebox(0,0){$\bullet$}}
\put(4,1.5){\line(1,0){15}}
\put(4,8){\makebox(0,0){$\scriptstyle 0$}}
\put(20,8){\makebox(0,0){$\scriptstyle 1$}}
\put(36,8){\makebox(0,0){$\scriptstyle 0$}}
\end{picture}\oplus
\begin{picture}(42,5)
\put(20,.1){\line(1,0){14.3}}
\put(20,3.2){\line(1,0){14.3}}
\put(20,1.5){\makebox(0,0){$\bullet$}}
\put(36,1.6){\makebox(0,0){$\times$}}
\put(28,1.5){\makebox(0,0){$\rangle$}}
\put(4,1.5){\makebox(0,0){$\bullet$}}
\put(4,1.5){\line(1,0){15}}
\put(4,8){\makebox(0,0){$\scriptstyle 1$}}
\put(20,8){\makebox(0,0){$\scriptstyle 0$}}
\put(36,8){\makebox(0,0){$\scriptstyle 0$}}
\end{picture}\oplus\!
\begin{array}{c}
\begin{picture}(42,5)
\put(20,.1){\line(1,0){14.3}}
\put(20,3.2){\line(1,0){14.3}}
\put(20,1.5){\makebox(0,0){$\bullet$}}
\put(36,1.6){\makebox(0,0){$\times$}}
\put(28,1.5){\makebox(0,0){$\rangle$}}
\put(4,1.5){\makebox(0,0){$\bullet$}}
\put(4,1.5){\line(1,0){15}}
\put(4,8){\makebox(0,0){$\scriptstyle 1$}}
\put(20,8){\makebox(0,0){$\scriptstyle 1$}}
\put(36,8){\makebox(0,0){$\scriptstyle -2$}}
\end{picture}\\[0pt]
\oplus\\[2pt]
\begin{picture}(42,5)
\put(20,.1){\line(1,0){14.3}}
\put(20,3.2){\line(1,0){14.3}}
\put(20,1.5){\makebox(0,0){$\bullet$}}
\put(36,1.6){\makebox(0,0){$\times$}}
\put(28,1.5){\makebox(0,0){$\rangle$}}
\put(4,1.5){\makebox(0,0){$\bullet$}}
\put(4,1.5){\line(1,0){15}}
\put(4,8){\makebox(0,0){$\scriptstyle 0$}}
\put(20,8){\makebox(0,0){$\scriptstyle 0$}}
\put(36,8){\makebox(0,0){$\scriptstyle 0$}}
\end{picture}
\end{array}\!\oplus
\begin{picture}(42,5)
\put(20,.1){\line(1,0){14.3}}
\put(20,3.2){\line(1,0){14.3}}
\put(20,1.5){\makebox(0,0){$\bullet$}}
\put(36,1.6){\makebox(0,0){$\times$}}
\put(28,1.5){\makebox(0,0){$\rangle$}}
\put(4,1.5){\makebox(0,0){$\bullet$}}
\put(4,1.5){\line(1,0){15}}
\put(4,8){\makebox(0,0){$\scriptstyle 0$}}
\put(20,8){\makebox(0,0){$\scriptstyle 1$}}
\put(36,8){\makebox(0,0){$\scriptstyle -2$}}
\end{picture}\oplus
\begin{picture}(42,5)
\put(20,.1){\line(1,0){14.3}}
\put(20,3.2){\line(1,0){14.3}}
\put(20,1.5){\makebox(0,0){$\bullet$}}
\put(36,1.6){\makebox(0,0){$\times$}}
\put(28,1.5){\makebox(0,0){$\rangle$}}
\put(4,1.5){\makebox(0,0){$\bullet$}}
\put(4,1.5){\line(1,0){15}}
\put(4,8){\makebox(0,0){$\scriptstyle 1$}}
\put(20,8){\makebox(0,0){$\scriptstyle 0$}}
\put(36,8){\makebox(0,0){$\scriptstyle -2$}}
\end{picture}\end{array}\end{equation}
and one can see from this grading that the corresponding $6$-dimensional
homogeneous space $G/P$ is equipped with a canonical $3$-dimensional
distribution. The corresponding infinitesimal flag 
structure~\cite[\S3.1.6]{CSl} is exactly the geometry of such $3$-distributions
and the irreducible bundles are related in our two notations by
$$(a,b,c)=\enskip
\begin{picture}(54,5)
\put(26,.1){\line(1,0){20.3}}
\put(26,3.2){\line(1,0){20.3}}
\put(26,1.5){\makebox(0,0){$\bullet$}}
\put(48,1.6){\makebox(0,0){$\times$}}
\put(37,1.5){\makebox(0,0){$\rangle$}}
\put(4,1.5){\makebox(0,0){$\bullet$}}
\put(4,1.5){\line(1,0){21}}
\put(4,8){\makebox(0,0){$\scriptstyle b-a$}}
\put(26,8){\makebox(0,0){$\scriptstyle c-b$}}
\put(48,8){\makebox(0,0){$\scriptstyle -2c$}}
\end{picture}.$$
Recall in the proof of Theorem~\ref{preferred} that we reduced the freedom in 
splitting $\Lambda^1$ to $(0,0,1)$ in (\ref{basicfreedom}). In the Dynkin 
diagram notation this remaining freedom lies in
\begin{picture}(42,5)
\put(20,.1){\line(1,0){14.3}}
\put(20,3.2){\line(1,0){14.3}}
\put(20,1.5){\makebox(0,0){$\bullet$}}
\put(36,1.6){\makebox(0,0){$\times$}}
\put(28,1.5){\makebox(0,0){$\rangle$}}
\put(4,1.5){\makebox(0,0){$\bullet$}}
\put(4,1.5){\line(1,0){15}}
\put(4,8){\makebox(0,0){$\scriptstyle 0$}}
\put(20,8){\makebox(0,0){$\scriptstyle 1$}}
\put(36,8){\makebox(0,0){$\scriptstyle -2$}}
\end{picture}, which is exactly the action of ${\mathfrak{g}}_1$ on 
$${\mathfrak{g}}/{\mathfrak{p}}={\mathfrak{g}}_{-2}\oplus{\mathfrak{g}}_{-1}=
\begin{picture}(42,5)
\put(20,.1){\line(1,0){14.3}}
\put(20,3.2){\line(1,0){14.3}}
\put(20,1.5){\makebox(0,0){$\bullet$}}
\put(36,1.6){\makebox(0,0){$\times$}}
\put(28,1.5){\makebox(0,0){$\rangle$}}
\put(4,1.5){\makebox(0,0){$\bullet$}}
\put(4,1.5){\line(1,0){15}}
\put(4,8){\makebox(0,0){$\scriptstyle 0$}}
\put(20,8){\makebox(0,0){$\scriptstyle 1$}}
\put(36,8){\makebox(0,0){$\scriptstyle 0$}}
\end{picture}\oplus
\begin{picture}(42,5)
\put(20,.1){\line(1,0){14.3}}
\put(20,3.2){\line(1,0){14.3}}
\put(20,1.5){\makebox(0,0){$\bullet$}}
\put(36,1.6){\makebox(0,0){$\times$}}
\put(28,1.5){\makebox(0,0){$\rangle$}}
\put(4,1.5){\makebox(0,0){$\bullet$}}
\put(4,1.5){\line(1,0){15}}
\put(4,8){\makebox(0,0){$\scriptstyle 1$}}
\put(20,8){\makebox(0,0){$\scriptstyle 0$}}
\put(36,8){\makebox(0,0){$\scriptstyle 0$}}
\end{picture},$$ 
as can be seen in~(\ref{grading}). The geometric import of this observation is
that Theorem~\ref{preferred} reduces the structure group of the tangent bundle
from general $H$-preserving and Levi-form-preserving automorphisms to the
subgroup of ${\mathrm{Aut}}({\mathfrak{g}}/{\mathfrak{p}})$ defined by the
Adjoint action of~$P$, namely the group $P/\exp({\mathfrak{g}}_2)$ with Lie
algebra~${\mathfrak{g}}_0\oplus{\mathfrak{g}}_1$. Dually, the $2$-adapted
co-framings are preserved by exactly this group.

Finally, we can take the complex of Theorem~\ref{preBGG}, use the 
splitting of $\big[(0,2,2)+(1,1,3)\big]$ afforded by Theorem~\ref{preferred},
and write the result in Dynkin diagram notation to obtain the following.
\begin{theorem}
On any smooth $6$-manifold equipped with a generic $3$-distribution,
there is a canonically defined locally exact differential complex
$$\begin{array}{r}0\to{\mathbb{R}}\to\!
\begin{picture}(42,5)
\put(20,.1){\line(1,0){14.3}}
\put(20,3.2){\line(1,0){14.3}}
\put(20,1.5){\makebox(0,0){$\bullet$}}
\put(36,1.6){\makebox(0,0){$\times$}}
\put(28,1.5){\makebox(0,0){$\rangle$}}
\put(4,1.5){\makebox(0,0){$\bullet$}}
\put(4,1.5){\line(1,0){15}}
\put(4,8){\makebox(0,0){$\scriptstyle 0$}}
\put(20,8){\makebox(0,0){$\scriptstyle 0$}}
\put(36,8){\makebox(0,0){$\scriptstyle 0$}}
\end{picture}
\!\to\!
\begin{picture}(42,5)
\put(20,.1){\line(1,0){14.3}}
\put(20,3.2){\line(1,0){14.3}}
\put(20,1.5){\makebox(0,0){$\bullet$}}
\put(36,1.6){\makebox(0,0){$\times$}}
\put(28,1.5){\makebox(0,0){$\rangle$}}
\put(4,1.5){\makebox(0,0){$\bullet$}}
\put(4,1.5){\line(1,0){15}}
\put(4,8){\makebox(0,0){$\scriptstyle 0$}}
\put(20,8){\makebox(0,0){$\scriptstyle 1$}}
\put(36,8){\makebox(0,0){$\scriptstyle -2$}}
\end{picture}
\!\to\!
\begin{picture}(42,5)
\put(20,.1){\line(1,0){14.3}}
\put(20,3.2){\line(1,0){14.3}}
\put(20,1.5){\makebox(0,0){$\bullet$}}
\put(36,1.6){\makebox(0,0){$\times$}}
\put(28,1.5){\makebox(0,0){$\rangle$}}
\put(4,1.5){\makebox(0,0){$\bullet$}}
\put(4,1.5){\line(1,0){15}}
\put(4,8){\makebox(0,0){$\scriptstyle 1$}}
\put(20,8){\makebox(0,0){$\scriptstyle 1$}}
\put(36,8){\makebox(0,0){$\scriptstyle -4$}}
\end{picture}
\begin{array}{c}\nearrow\\ \searrow\end{array}\!\!
\begin{array}{c}
\begin{picture}(42,5)
\put(20,.1){\line(1,0){14.3}}
\put(20,3.2){\line(1,0){14.3}}
\put(20,1.5){\makebox(0,0){$\bullet$}}
\put(36,1.6){\makebox(0,0){$\times$}}
\put(28,1.5){\makebox(0,0){$\rangle$}}
\put(4,1.5){\makebox(0,0){$\bullet$}}
\put(4,1.5){\line(1,0){15}}
\put(4,8){\makebox(0,0){$\scriptstyle 2$}}
\put(20,8){\makebox(0,0){$\scriptstyle 0$}}
\put(36,8){\makebox(0,0){$\scriptstyle -4$}}
\end{picture}\\[2pt]
\oplus\\[4pt]
\begin{picture}(42,5)
\put(20,.1){\line(1,0){14.3}}
\put(20,3.2){\line(1,0){14.3}}
\put(20,1.5){\makebox(0,0){$\bullet$}}
\put(36,1.6){\makebox(0,0){$\times$}}
\put(28,1.5){\makebox(0,0){$\rangle$}}
\put(4,1.5){\makebox(0,0){$\bullet$}}
\put(4,1.5){\line(1,0){15}}
\put(4,8){\makebox(0,0){$\scriptstyle 0$}}
\put(20,8){\makebox(0,0){$\scriptstyle 2$}}
\put(36,8){\makebox(0,0){$\scriptstyle -6$}}
\end{picture}
\end{array}
\!\!\begin{array}{c}\searrow\\ \nearrow\end{array}\qquad\\[20pt]
\begin{picture}(42,5)
\put(20,.1){\line(1,0){14.3}}
\put(20,3.2){\line(1,0){14.3}}
\put(20,1.5){\makebox(0,0){$\bullet$}}
\put(36,1.6){\makebox(0,0){$\times$}}
\put(28,1.5){\makebox(0,0){$\rangle$}}
\put(4,1.5){\makebox(0,0){$\bullet$}}
\put(4,1.5){\line(1,0){15}}
\put(4,8){\makebox(0,0){$\scriptstyle 1$}}
\put(20,8){\makebox(0,0){$\scriptstyle 1$}}
\put(36,8){\makebox(0,0){$\scriptstyle -6$}}
\end{picture}
\!\to\!
\begin{picture}(42,5)
\put(20,.1){\line(1,0){14.3}}
\put(20,3.2){\line(1,0){14.3}}
\put(20,1.5){\makebox(0,0){$\bullet$}}
\put(36,1.6){\makebox(0,0){$\times$}}
\put(28,1.5){\makebox(0,0){$\rangle$}}
\put(4,1.5){\makebox(0,0){$\bullet$}}
\put(4,1.5){\line(1,0){15}}
\put(4,8){\makebox(0,0){$\scriptstyle 1$}}
\put(20,8){\makebox(0,0){$\scriptstyle 0$}}
\put(36,8){\makebox(0,0){$\scriptstyle -6$}}
\end{picture}
\!\to\!
\begin{picture}(42,5)
\put(20,.1){\line(1,0){14.3}}
\put(20,3.2){\line(1,0){14.3}}
\put(20,1.5){\makebox(0,0){$\bullet$}}
\put(36,1.6){\makebox(0,0){$\times$}}
\put(28,1.5){\makebox(0,0){$\rangle$}}
\put(4,1.5){\makebox(0,0){$\bullet$}}
\put(4,1.5){\line(1,0){15}}
\put(4,8){\makebox(0,0){$\scriptstyle 0$}}
\put(20,8){\makebox(0,0){$\scriptstyle 0$}}
\put(36,8){\makebox(0,0){$\scriptstyle -6$}}
\end{picture}
\!\to 0.\end{array}$$
\end{theorem}
\noindent This is the BGG complex in standard notation.

\section{The Engel complex revisited}\label{engelrevisited}
Although the complex $\Lambda^0\xrightarrow{\,d_H\,}\Lambda_H^1
\xrightarrow{\,{\mathcal{P}}\,}\lambda\xi^2$ constructed in \S\ref{engel} used
nothing beyond an Engel structure, for the full-blown resolution
(\ref{fullEngel}) it was necessary to choose some extra structure, namely a
complement to $\xi\subset\Lambda_H^1$ (equivalently, a complement to
$(\xi+K)^\perp\subset H$, the {\em Engel line field\/}~\cite{M}). As pointed
out to us by Boris Doubrov, there is a unique homogeneous space of the form
$G/P$, for $G$ semisimple and $P$ parabolic, that carries a $G$-invariant Engel
structure. Specifically, if $G={\mathrm{Sp}}(4,{\mathbb{R}})$ and $P$ is its 
Borel subgroup, then 
$${\mathfrak{g}}={\mathfrak{sp}}(4,{\mathbb{R}})={\mathfrak{g}}_{-3}\oplus
{\mathfrak{g}}_{-2}\oplus{\mathfrak{g}}_{-1}\oplus
\underbrace{{\mathfrak{g}}_0\oplus{\mathfrak{g}}_1
\oplus{\mathfrak{g}}_2\oplus{\mathfrak{g}}_3}_{\mathfrak{p}},$$
which, in Dynkin diagram notation, reads
$$\begin{picture}(24,5)
\put(4,0){\line(1,0){16}}
\put(4,3.2){\line(1,0){16}}
\put(4,1.3){\makebox(0,0){$\bullet$}}
\put(20,1.4){\makebox(0,0){$\bullet$}}
\put(12,1.5){\makebox(0,0){$\langle$}}
\put(4,8){\makebox(0,0){$\scriptstyle 2$}}
\put(20,8){\makebox(0,0){$\scriptstyle 0$}}
\end{picture}=
\begin{picture}(24,5)
\put(5.6,0){\line(1,0){12.8}}
\put(5.6,3.2){\line(1,0){12.8}}
\put(4,1.5){\makebox(0,0){$\times$}}
\put(20,1.5){\makebox(0,0){$\times$}}
\put(12,1.5){\makebox(0,0){$\langle$}}
\put(4,8){\makebox(0,0){$\scriptstyle 2$}}
\put(20,8){\makebox(0,0){$\scriptstyle 0$}}
\end{picture}\oplus
\begin{picture}(24,5)
\put(5.6,0){\line(1,0){12.8}}
\put(5.6,3.2){\line(1,0){12.8}}
\put(4,1.5){\makebox(0,0){$\times$}}
\put(20,1.5){\makebox(0,0){$\times$}}
\put(12,1.5){\makebox(0,0){$\langle$}}
\put(4,8){\makebox(0,0){$\scriptstyle 0$}}
\put(20,8){\makebox(0,0){$\scriptstyle 1$}}
\end{picture}\oplus\!\!
\begin{array}c\begin{picture}(24,5)
\put(5.6,0){\line(1,0){12.8}}
\put(5.6,3.2){\line(1,0){12.8}}
\put(4,1.5){\makebox(0,0){$\times$}}
\put(20,1.5){\makebox(0,0){$\times$}}
\put(12,1.5){\makebox(0,0){$\langle$}}
\put(4,8){\makebox(0,0){$\scriptstyle 2$}}
\put(20,8){\makebox(0,0){$\scriptstyle -1$}}
\end{picture}\\ \oplus\\
\begin{picture}(24,5)
\put(5.6,0){\line(1,0){12.8}}
\put(5.6,3.2){\line(1,0){12.8}}
\put(4,1.5){\makebox(0,0){$\times$}}
\put(20,1.5){\makebox(0,0){$\times$}}
\put(12,1.5){\makebox(0,0){$\langle$}}
\put(4,8){\makebox(0,0){$\scriptstyle -2$}}
\put(20,8){\makebox(0,0){$\scriptstyle 2$}}
\end{picture}\end{array}\!\!\oplus\!\!
\begin{array}c\begin{picture}(24,5)
\put(5.6,0){\line(1,0){12.8}}
\put(5.6,3.2){\line(1,0){12.8}}
\put(4,1.5){\makebox(0,0){$\times$}}
\put(20,1.5){\makebox(0,0){$\times$}}
\put(12,1.5){\makebox(0,0){$\langle$}}
\put(4,8){\makebox(0,0){$\scriptstyle 0$}}
\put(20,8){\makebox(0,0){$\scriptstyle 0$}}
\end{picture}\\ \oplus\\
\begin{picture}(24,5)
\put(5.6,0){\line(1,0){12.8}}
\put(5.6,3.2){\line(1,0){12.8}}
\put(4,1.5){\makebox(0,0){$\times$}}
\put(20,1.5){\makebox(0,0){$\times$}}
\put(12,1.5){\makebox(0,0){$\langle$}}
\put(4,8){\makebox(0,0){$\scriptstyle 0$}}
\put(20,8){\makebox(0,0){$\scriptstyle 0$}}
\end{picture}\end{array}\!\!\oplus\!\!
\begin{array}c\begin{picture}(24,5)
\put(5.6,0){\line(1,0){12.8}}
\put(5.6,3.2){\line(1,0){12.8}}
\put(4,1.5){\makebox(0,0){$\times$}}
\put(20,1.5){\makebox(0,0){$\times$}}
\put(12,1.5){\makebox(0,0){$\langle$}}
\put(4,8){\makebox(0,0){$\scriptstyle 2$}}
\put(20,8){\makebox(0,0){$\scriptstyle -2$}}
\end{picture}\\ \oplus\\
\begin{picture}(24,5)
\put(5.6,0){\line(1,0){12.8}}
\put(5.6,3.2){\line(1,0){12.8}}
\put(4,1.5){\makebox(0,0){$\times$}}
\put(20,1.5){\makebox(0,0){$\times$}}
\put(12,1.5){\makebox(0,0){$\langle$}}
\put(4,8){\makebox(0,0){$\scriptstyle -2$}}
\put(20,8){\makebox(0,0){$\scriptstyle 1$}}
\end{picture}\end{array}\!\!\oplus
\begin{picture}(24,5)
\put(5.6,0){\line(1,0){12.8}}
\put(5.6,3.2){\line(1,0){12.8}}
\put(4,1.5){\makebox(0,0){$\times$}}
\put(20,1.5){\makebox(0,0){$\times$}}
\put(12,1.5){\makebox(0,0){$\langle$}}
\put(4,8){\makebox(0,0){$\scriptstyle 0$}}
\put(20,8){\makebox(0,0){$\scriptstyle -1$}}
\end{picture}\oplus
\begin{picture}(24,5)
\put(5.6,0){\line(1,0){12.8}}
\put(5.6,3.2){\line(1,0){12.8}}
\put(4,1.5){\makebox(0,0){$\times$}}
\put(20,1.5){\makebox(0,0){$\times$}}
\put(12,1.5){\makebox(0,0){$\langle$}}
\put(4,8){\makebox(0,0){$\scriptstyle -2$}}
\put(20,8){\makebox(0,0){$\scriptstyle 0$}}
\end{picture}.$$
The $1$-forms on this homogeneous space $G/P$ are filtered
$$\Lambda^1=\begin{array}c\begin{picture}(24,5)
\put(5.6,0){\line(1,0){12.8}}
\put(5.6,3.2){\line(1,0){12.8}}
\put(4,1.5){\makebox(0,0){$\times$}}
\put(20,1.5){\makebox(0,0){$\times$}}
\put(12,1.5){\makebox(0,0){$\langle$}}
\put(4,8){\makebox(0,0){$\scriptstyle 2$}}
\put(20,8){\makebox(0,0){$\scriptstyle -2$}}
\end{picture}\\ \oplus\\
\begin{picture}(24,5)
\put(5.6,0){\line(1,0){12.8}}
\put(5.6,3.2){\line(1,0){12.8}}
\put(4,1.5){\makebox(0,0){$\times$}}
\put(20,1.5){\makebox(0,0){$\times$}}
\put(12,1.5){\makebox(0,0){$\langle$}}
\put(4,8){\makebox(0,0){$\scriptstyle -2$}}
\put(20,8){\makebox(0,0){$\scriptstyle 1$}}
\end{picture}\end{array}\!\!+
\begin{picture}(24,5)
\put(5.6,0){\line(1,0){12.8}}
\put(5.6,3.2){\line(1,0){12.8}}
\put(4,1.5){\makebox(0,0){$\times$}}
\put(20,1.5){\makebox(0,0){$\times$}}
\put(12,1.5){\makebox(0,0){$\langle$}}
\put(4,8){\makebox(0,0){$\scriptstyle 0$}}
\put(20,8){\makebox(0,0){$\scriptstyle -1$}}
\end{picture}+
\begin{picture}(24,5)
\put(5.6,0){\line(1,0){12.8}}
\put(5.6,3.2){\line(1,0){12.8}}
\put(4,1.5){\makebox(0,0){$\times$}}
\put(20,1.5){\makebox(0,0){$\times$}}
\put(12,1.5){\makebox(0,0){$\langle$}}
\put(4,8){\makebox(0,0){$\scriptstyle -2$}}
\put(20,8){\makebox(0,0){$\scriptstyle 0$}}
\end{picture}=
\begin{array}c\lambda\\ \oplus\\ \xi\end{array}\!\!+
\lambda\xi+\lambda\xi^2$$
and the corresponding regular infinitesimal flag structure is exactly that of 
an Engel manifold equipped with a choice of splitting 
$\Lambda_H^1=\lambda\oplus\xi$ as discussed in~\S\ref{engel}. The BGG complex 
(\ref{fullEngel}) in Dynkin diagram notation reads
$$0\to{\mathbb{R}}\to\begin{picture}(24,5)
\put(5.6,0){\line(1,0){12.8}}
\put(5.6,3.2){\line(1,0){12.8}}
\put(4,1.5){\makebox(0,0){$\times$}}
\put(20,1.5){\makebox(0,0){$\times$}}
\put(12,1.5){\makebox(0,0){$\langle$}}
\put(4,8){\makebox(0,0){$\scriptstyle 0$}}
\put(20,8){\makebox(0,0){$\scriptstyle 0$}}
\end{picture}\begin{array}{c}\nearrow\\ \searrow\end{array}\!\!
\begin{array}c\begin{picture}(24,5)
\put(5.6,0){\line(1,0){12.8}}
\put(5.6,3.2){\line(1,0){12.8}}
\put(4,1.5){\makebox(0,0){$\times$}}
\put(20,1.5){\makebox(0,0){$\times$}}
\put(12,1.5){\makebox(0,0){$\langle$}}
\put(4,8){\makebox(0,0){$\scriptstyle 2$}}
\put(20,8){\makebox(0,0){$\scriptstyle -2$}}
\end{picture}\\[2pt] \oplus\\[4pt]
\begin{picture}(24,5)
\put(5.6,0){\line(1,0){12.8}}
\put(5.6,3.2){\line(1,0){12.8}}
\put(4,1.5){\makebox(0,0){$\times$}}
\put(20,1.5){\makebox(0,0){$\times$}}
\put(12,1.5){\makebox(0,0){$\langle$}}
\put(4,8){\makebox(0,0){$\scriptstyle -2$}}
\put(20,8){\makebox(0,0){$\scriptstyle 1$}}
\end{picture}\end{array}\!\!
\begin{array}c\longrightarrow\\ \mbox{\Large\begin{picture}(0,0)
\put(0,0){\makebox(0,0){$\nearrow$}}
\put(0,0){\makebox(0,0){$\searrow$}}
\end{picture}}\\[4pt] \longrightarrow\end{array}\!\!
\begin{array}c\begin{picture}(24,5)
\put(5.6,0){\line(1,0){12.8}}
\put(5.6,3.2){\line(1,0){12.8}}
\put(4,1.5){\makebox(0,0){$\times$}}
\put(20,1.5){\makebox(0,0){$\times$}}
\put(12,1.5){\makebox(0,0){$\langle$}}
\put(4,8){\makebox(0,0){$\scriptstyle 2$}}
\put(20,8){\makebox(0,0){$\scriptstyle -3$}}
\end{picture}\\[2pt] \oplus\\[4pt]
\begin{picture}(24,5)
\put(5.6,0){\line(1,0){12.8}}
\put(5.6,3.2){\line(1,0){12.8}}
\put(4,1.5){\makebox(0,0){$\times$}}
\put(20,1.5){\makebox(0,0){$\times$}}
\put(12,1.5){\makebox(0,0){$\langle$}}
\put(4,8){\makebox(0,0){$\scriptstyle -4$}}
\put(20,8){\makebox(0,0){$\scriptstyle 1$}}
\end{picture}\end{array}\!\!
\begin{array}c\longrightarrow\\ \mbox{\Large\begin{picture}(0,0)
\put(0,0){\makebox(0,0){$\nearrow$}}
\put(0,0){\makebox(0,0){$\searrow$}}
\end{picture}}\\[4pt] \longrightarrow\end{array}\!\!
\begin{array}c\begin{picture}(24,5)
\put(5.6,0){\line(1,0){12.8}}
\put(5.6,3.2){\line(1,0){12.8}}
\put(4,1.5){\makebox(0,0){$\times$}}
\put(20,1.5){\makebox(0,0){$\times$}}
\put(12,1.5){\makebox(0,0){$\langle$}}
\put(4,8){\makebox(0,0){$\scriptstyle 0$}}
\put(20,8){\makebox(0,0){$\scriptstyle -3$}}
\end{picture}\\[2pt] \oplus\\[4pt]
\begin{picture}(24,5)
\put(5.6,0){\line(1,0){12.8}}
\put(5.6,3.2){\line(1,0){12.8}}
\put(4,1.5){\makebox(0,0){$\times$}}
\put(20,1.5){\makebox(0,0){$\times$}}
\put(12,1.5){\makebox(0,0){$\langle$}}
\put(4,8){\makebox(0,0){$\scriptstyle -4$}}
\put(20,8){\makebox(0,0){$\scriptstyle 0$}}
\end{picture}\end{array}\!\!
\begin{array}{c}\searrow\\ \nearrow\end{array}
\begin{picture}(24,5)
\put(5.6,0){\line(1,0){12.8}}
\put(5.6,3.2){\line(1,0){12.8}}
\put(4,1.5){\makebox(0,0){$\times$}}
\put(20,1.5){\makebox(0,0){$\times$}}
\put(12,1.5){\makebox(0,0){$\langle$}}
\put(4,8){\makebox(0,0){$\scriptstyle -2$}}
\put(20,8){\makebox(0,0){$\scriptstyle -2$}}
\end{picture}\to 0.$$

\section{Another geometry in five variables}
Recall that an Engel manifold is a $4$-dimensional manifold equipped with a
generic $2$-dimensional distribution. The geometry considered in \S\ref{engel}
and \S\ref{engelrevisited} was defined on an Engel manifold by a choice of
splitting of~$\Lambda_H^1$, the bundle of $1$-forms along~$H$ (rather than the
filtration that is canonically present). The geometry to be considered in this
section will very much resemble this case. 

Let us consider a Pfaffian system $I\subset\Lambda^1$ of rank $2$ on a smooth
$5$-manifold. As usual, we define the {\em Levi form\/} ${\mathcal{L}}$ as the
composition
$$I\to\Lambda^1\xrightarrow{\,d\,}\Lambda^2\to\Lambda_H^2,$$
where~$H\equiv I^\perp$. Notice that $\Lambda_H^2$ has rank $3$ and we shall
suppose that ${\mathcal{L}}$ is injective, as is generically the case. Under
the canonical identification $\Lambda_H^2=\Lambda_H^3\otimes H$ we see that the
rank $2$ sub-bundle ${\mathcal{L}}(I)\subset\Lambda_H^2$ gives rise to a rank
$2$ sub-distribution $D\subset H$. In~\cite{M2} it is observed that $D$ is the
unique rank $2$ sub-bundle of $H$ such that $[D,D]\subseteq H$. It is not
necessarily the case, however, that $[D,D]=H$ (in which case we would be back
in the five variables geometry of~\S\ref{five}). To proceed further, let us
write $L$ for the line sub-bundle $D^\perp\subset\Lambda_H^1$ and choose a
complementary rank $2$ sub-bundle~$Q$ so that we now have a splitting
$\Lambda_H^1=Q\oplus L$. This completes the definition of the structure to be
considered in this section. Equivalently, we are considering a $5$-manifold $M$
equipped with a pair of transverse distributions $D$ and $\ell$ of ranks $2$
and~$1$, respectively, such that
\begin{equation}\label{anothergeometry}
[D,D]\subseteq \ell\oplus D\qquad\mbox{and}\qquad[\ell\oplus D,\ell\oplus D]
=TM.\end{equation}
This is precisely the regular infinitesimal flag structure associated with the
grading
$$\begin{array}{l}
{\mathfrak{g}}_{-2}\oplus{\mathfrak{g}}_{-1}\oplus
\underbrace{{\mathfrak{g}}_0
\oplus{\mathfrak{g}}_1\oplus{\mathfrak{g}}_2}_{\mathfrak{p}}=\\
\quad\begin{picture}(42,5)
\put(20,1.5){\makebox(0,0){$\times$}}
\put(36,1.5){\makebox(0,0){$\bullet$}}
\put(4,1.5){\makebox(0,0){$\times$}}
\put(4,1.5){\line(1,0){32}}
\put(4,8){\makebox(0,0){$\scriptstyle 1$}}
\put(20,8){\makebox(0,0){$\scriptstyle 0$}}
\put(36,8){\makebox(0,0){$\scriptstyle 1$}}
\end{picture}\oplus
\begin{array}{c}
\begin{picture}(42,5)
\put(20,1.5){\makebox(0,0){$\times$}}
\put(36,1.5){\makebox(0,0){$\bullet$}}
\put(4,1.5){\makebox(0,0){$\times$}}
\put(4,1.5){\line(1,0){32}}
\put(4,8){\makebox(0,0){$\scriptstyle 2$}}
\put(20,8){\makebox(0,0){$\scriptstyle -1$}}
\put(36,8){\makebox(0,0){$\scriptstyle 0$}}
\end{picture}\\[-2pt]
\oplus\\[1pt]
\begin{picture}(42,5)
\put(20,1.5){\makebox(0,0){$\times$}}
\put(36,1.5){\makebox(0,0){$\bullet$}}
\put(4,1.5){\makebox(0,0){$\times$}}
\put(4,1.5){\line(1,0){32}}
\put(4,8){\makebox(0,0){$\scriptstyle -1$}}
\put(20,8){\makebox(0,0){$\scriptstyle 1$}}
\put(36,8){\makebox(0,0){$\scriptstyle 1$}}
\end{picture}
\end{array}\oplus
\begin{array}{c}
\begin{picture}(42,5)
\put(20,1.5){\makebox(0,0){$\times$}}
\put(36,1.5){\makebox(0,0){$\bullet$}}
\put(4,1.5){\makebox(0,0){$\times$}}
\put(4,1.5){\line(1,0){32}}
\put(4,8){\makebox(0,0){$\scriptstyle 0$}}
\put(20,8){\makebox(0,0){$\scriptstyle 0$}}
\put(36,8){\makebox(0,0){$\scriptstyle 0$}}
\end{picture}\\[-2pt]
\oplus\\[1pt]
\begin{picture}(42,5)
\put(20,1.5){\makebox(0,0){$\times$}}
\put(36,1.5){\makebox(0,0){$\bullet$}}
\put(4,1.5){\makebox(0,0){$\times$}}
\put(4,1.5){\line(1,0){32}}
\put(4,8){\makebox(0,0){$\scriptstyle 0$}}
\put(20,8){\makebox(0,0){$\scriptstyle -1$}}
\put(36,8){\makebox(0,0){$\scriptstyle 2$}}
\end{picture}\\[-2pt]
\oplus\\[1pt]
\begin{picture}(42,5)
\put(20,1.5){\makebox(0,0){$\times$}}
\put(36,1.5){\makebox(0,0){$\bullet$}}
\put(4,1.5){\makebox(0,0){$\times$}}
\put(4,1.5){\line(1,0){32}}
\put(4,8){\makebox(0,0){$\scriptstyle 0$}}
\put(20,8){\makebox(0,0){$\scriptstyle 0$}}
\put(36,8){\makebox(0,0){$\scriptstyle 0$}}
\end{picture}
\end{array}\oplus
\begin{array}{c}
\begin{picture}(42,5)
\put(20,1.5){\makebox(0,0){$\times$}}
\put(36,1.5){\makebox(0,0){$\bullet$}}
\put(4,1.5){\makebox(0,0){$\times$}}
\put(4,1.5){\line(1,0){32}}
\put(4,8){\makebox(0,0){$\scriptstyle 1$}}
\put(20,8){\makebox(0,0){$\scriptstyle -2$}}
\put(36,8){\makebox(0,0){$\scriptstyle 1$}}
\end{picture}\\[-2pt]
\oplus\\[1pt]
\begin{picture}(42,5)
\put(20,1.5){\makebox(0,0){$\times$}}
\put(36,1.5){\makebox(0,0){$\bullet$}}
\put(4,1.5){\makebox(0,0){$\times$}}
\put(4,1.5){\line(1,0){32}}
\put(4,8){\makebox(0,0){$\scriptstyle -2$}}
\put(20,8){\makebox(0,0){$\scriptstyle 1$}}
\put(36,8){\makebox(0,0){$\scriptstyle 0$}}
\end{picture}
\end{array}\oplus
\begin{picture}(42,5)
\put(20,1.5){\makebox(0,0){$\times$}}
\put(36,1.5){\makebox(0,0){$\bullet$}}
\put(4,1.5){\makebox(0,0){$\times$}}
\put(4,1.5){\line(1,0){32}}
\put(4,8){\makebox(0,0){$\scriptstyle -1$}}
\put(20,8){\makebox(0,0){$\scriptstyle -1$}}
\put(36,8){\makebox(0,0){$\scriptstyle 1$}}
\end{picture}\end{array}$$
of ${\mathfrak{sl}}(4,{\mathbb{R}})$. The bundles
$$\textstyle\bigodot^c\!Q\otimes(\Lambda^2\!Q)^b\otimes L^a
\enskip\mbox{for}\enskip
c\in{\mathbb{Z}}_{\geq 0}\enskip\mbox{become}\quad
\begin{picture}(78,5)
\put(50,1.5){\makebox(0,0){$\times$}}
\put(74,1.5){\makebox(0,0){$\bullet$}}
\put(14,1.5){\makebox(0,0){$\times$}}
\put(14,1.5){\line(1,0){60}}
\put(14,8){\makebox(0,0){$\scriptstyle c+2b-2a$}}
\put(50,8){\makebox(0,0){$\scriptstyle a-3b-2c$}}
\put(74,8){\makebox(0,0){$\scriptstyle c$}}
\end{picture}$$
and the $1$-forms are filtered
\begin{equation}\label{oneforms}\Lambda^1=
\raisebox{2pt}{$\begin{array}{c}
Q\\[-2pt] \oplus\\[-2pt] L\end{array}$}+\,I=
\begin{array}{c}
\begin{picture}(42,5)
\put(20,1.5){\makebox(0,0){$\times$}}
\put(36,1.5){\makebox(0,0){$\bullet$}}
\put(4,1.5){\makebox(0,0){$\times$}}
\put(4,1.5){\line(1,0){32}}
\put(4,8){\makebox(0,0){$\scriptstyle 1$}}
\put(20,8){\makebox(0,0){$\scriptstyle -2$}}
\put(36,8){\makebox(0,0){$\scriptstyle 1$}}
\end{picture}\\[-2pt]
\oplus\\[1pt]
\begin{picture}(42,5)
\put(20,1.5){\makebox(0,0){$\times$}}
\put(36,1.5){\makebox(0,0){$\bullet$}}
\put(4,1.5){\makebox(0,0){$\times$}}
\put(4,1.5){\line(1,0){32}}
\put(4,8){\makebox(0,0){$\scriptstyle -2$}}
\put(20,8){\makebox(0,0){$\scriptstyle 1$}}
\put(36,8){\makebox(0,0){$\scriptstyle 0$}}
\end{picture}
\end{array}+
\begin{picture}(42,5)
\put(20,1.5){\makebox(0,0){$\times$}}
\put(36,1.5){\makebox(0,0){$\bullet$}}
\put(4,1.5){\makebox(0,0){$\times$}}
\put(4,1.5){\line(1,0){32}}
\put(4,8){\makebox(0,0){$\scriptstyle -1$}}
\put(20,8){\makebox(0,0){$\scriptstyle -1$}}
\put(36,8){\makebox(0,0){$\scriptstyle 1$}}
\end{picture},\end{equation}
as expected. This induces filtrations on the higher forms as follows.
$$\Lambda^2=
\begin{array}{c}
\begin{picture}(42,5)
\put(20,1.5){\makebox(0,0){$\times$}}
\put(36,1.5){\makebox(0,0){$\bullet$}}
\put(4,1.5){\makebox(0,0){$\times$}}
\put(4,1.5){\line(1,0){32}}
\put(4,8){\makebox(0,0){$\scriptstyle 2$}}
\put(20,8){\makebox(0,0){$\scriptstyle -3$}}
\put(36,8){\makebox(0,0){$\scriptstyle 0$}}
\end{picture}\\[-2pt]
\oplus\\[1pt]
\begin{picture}(42,5)
\put(20,1.5){\makebox(0,0){$\times$}}
\put(36,1.5){\makebox(0,0){$\bullet$}}
\put(4,1.5){\makebox(0,0){$\times$}}
\put(4,1.5){\line(1,0){32}}
\put(4,8){\makebox(0,0){$\scriptstyle -1$}}
\put(20,8){\makebox(0,0){$\scriptstyle -1$}}
\put(36,8){\makebox(0,0){$\scriptstyle 1$}}
\end{picture}
\end{array}+
\begin{array}{c}
\begin{picture}(42,5)
\put(20,1.5){\makebox(0,0){$\times$}}
\put(36,1.5){\makebox(0,0){$\bullet$}}
\put(4,1.5){\makebox(0,0){$\times$}}
\put(4,1.5){\line(1,0){32}}
\put(4,8){\makebox(0,0){$\scriptstyle 0$}}
\put(20,8){\makebox(0,0){$\scriptstyle -3$}}
\put(36,8){\makebox(0,0){$\scriptstyle 2$}}
\end{picture}\\[-2pt]
\oplus\\[1pt]
\begin{picture}(42,5)
\put(20,1.5){\makebox(0,0){$\times$}}
\put(36,1.5){\makebox(0,0){$\bullet$}}
\put(4,1.5){\makebox(0,0){$\times$}}
\put(4,1.5){\line(1,0){32}}
\put(4,8){\makebox(0,0){$\scriptstyle 0$}}
\put(20,8){\makebox(0,0){$\scriptstyle -2$}}
\put(36,8){\makebox(0,0){$\scriptstyle 0$}}
\end{picture}\\[-2pt]
\oplus\\[1pt]
\begin{picture}(42,5)
\put(20,1.5){\makebox(0,0){$\times$}}
\put(36,1.5){\makebox(0,0){$\bullet$}}
\put(4,1.5){\makebox(0,0){$\times$}}
\put(4,1.5){\line(1,0){32}}
\put(4,8){\makebox(0,0){$\scriptstyle -3$}}
\put(20,8){\makebox(0,0){$\scriptstyle 0$}}
\put(36,8){\makebox(0,0){$\scriptstyle 1$}}
\end{picture}
\end{array}+
\begin{picture}(42,5)
\put(20,1.5){\makebox(0,0){$\times$}}
\put(36,1.5){\makebox(0,0){$\bullet$}}
\put(4,1.5){\makebox(0,0){$\times$}}
\put(4,1.5){\line(1,0){32}}
\put(4,8){\makebox(0,0){$\scriptstyle -2$}}
\put(20,8){\makebox(0,0){$\scriptstyle -1$}}
\put(36,8){\makebox(0,0){$\scriptstyle 0$}}
\end{picture}$$
and
$$\begin{array}{r}
\Lambda^3=
\begin{picture}(42,5)
\put(20,1.5){\makebox(0,0){$\times$}}
\put(36,1.5){\makebox(0,0){$\bullet$}}
\put(4,1.5){\makebox(0,0){$\times$}}
\put(4,1.5){\line(1,0){32}}
\put(4,8){\makebox(0,0){$\scriptstyle 0$}}
\put(20,8){\makebox(0,0){$\scriptstyle -2$}}
\put(36,8){\makebox(0,0){$\scriptstyle 0$}}
\end{picture}+
\begin{array}{c}
\begin{picture}(42,5)
\put(20,1.5){\makebox(0,0){$\times$}}
\put(36,1.5){\makebox(0,0){$\bullet$}}
\put(4,1.5){\makebox(0,0){$\times$}}
\put(4,1.5){\line(1,0){32}}
\put(4,8){\makebox(0,0){$\scriptstyle 1$}}
\put(20,8){\makebox(0,0){$\scriptstyle -4$}}
\put(36,8){\makebox(0,0){$\scriptstyle 1$}}
\end{picture}\\[-2pt]
\oplus\\[1pt]
\begin{picture}(42,5)
\put(20,1.5){\makebox(0,0){$\times$}}
\put(36,1.5){\makebox(0,0){$\bullet$}}
\put(4,1.5){\makebox(0,0){$\times$}}
\put(4,1.5){\line(1,0){32}}
\put(4,8){\makebox(0,0){$\scriptstyle -2$}}
\put(20,8){\makebox(0,0){$\scriptstyle -1$}}
\put(36,8){\makebox(0,0){$\scriptstyle 0$}}
\end{picture}\\[-2pt]
\oplus\\[1pt]
\begin{picture}(42,5)
\put(20,1.5){\makebox(0,0){$\times$}}
\put(36,1.5){\makebox(0,0){$\bullet$}}
\put(4,1.5){\makebox(0,0){$\times$}}
\put(4,1.5){\line(1,0){32}}
\put(4,8){\makebox(0,0){$\scriptstyle -2$}}
\put(20,8){\makebox(0,0){$\scriptstyle -2$}}
\put(36,8){\makebox(0,0){$\scriptstyle 2$}}
\end{picture}
\end{array}+
\begin{array}{c}
\begin{picture}(42,5)
\put(20,1.5){\makebox(0,0){$\times$}}
\put(36,1.5){\makebox(0,0){$\bullet$}}
\put(4,1.5){\makebox(0,0){$\times$}}
\put(4,1.5){\line(1,0){32}}
\put(4,8){\makebox(0,0){$\scriptstyle -1$}}
\put(20,8){\makebox(0,0){$\scriptstyle -3$}}
\put(36,8){\makebox(0,0){$\scriptstyle 1$}}
\end{picture}\\[-2pt] \oplus\\[1pt]
\begin{picture}(42,5)
\put(20,1.5){\makebox(0,0){$\times$}}
\put(36,1.5){\makebox(0,0){$\bullet$}}
\put(4,1.5){\makebox(0,0){$\times$}}
\put(4,1.5){\line(1,0){32}}
\put(4,8){\makebox(0,0){$\scriptstyle -4$}}
\put(20,8){\makebox(0,0){$\scriptstyle 0$}}
\put(36,8){\makebox(0,0){$\scriptstyle 0$}}
\end{picture}
\end{array}\\
\Lambda^4=
\begin{picture}(42,5)
\put(20,1.5){\makebox(0,0){$\times$}}
\put(36,1.5){\makebox(0,0){$\bullet$}}
\put(4,1.5){\makebox(0,0){$\times$}}
\put(4,1.5){\line(1,0){32}}
\put(4,8){\makebox(0,0){$\scriptstyle -1$}}
\put(20,8){\makebox(0,0){$\scriptstyle -3$}}
\put(36,8){\makebox(0,0){$\scriptstyle 1$}}
\end{picture}+
\begin{array}{c}
\begin{picture}(42,5)
\put(20,1.5){\makebox(0,0){$\times$}}
\put(36,1.5){\makebox(0,0){$\bullet$}}
\put(4,1.5){\makebox(0,0){$\times$}}
\put(4,1.5){\line(1,0){32}}
\put(4,8){\makebox(0,0){$\scriptstyle 0$}}
\put(20,8){\makebox(0,0){$\scriptstyle -4$}}
\put(36,8){\makebox(0,0){$\scriptstyle 0$}}
\end{picture}\\[-2pt]
\oplus\\[1pt]
\begin{picture}(42,5)
\put(20,1.5){\makebox(0,0){$\times$}}
\put(36,1.5){\makebox(0,0){$\bullet$}}
\put(4,1.5){\makebox(0,0){$\times$}}
\put(4,1.5){\line(1,0){32}}
\put(4,8){\makebox(0,0){$\scriptstyle -3$}}
\put(20,8){\makebox(0,0){$\scriptstyle -2$}}
\put(36,8){\makebox(0,0){$\scriptstyle 1$}}
\end{picture}
\end{array}\qquad\qquad\qquad
\Lambda^5=\begin{picture}(42,5)
\put(20,1.5){\makebox(0,0){$\times$}}
\put(36,1.5){\makebox(0,0){$\bullet$}}
\put(4,1.5){\makebox(0,0){$\times$}}
\put(4,1.5){\line(1,0){32}}
\put(4,8){\makebox(0,0){$\scriptstyle -2$}}
\put(20,8){\makebox(0,0){$\scriptstyle -3$}}
\put(36,8){\makebox(0,0){$\scriptstyle 0$}}
\end{picture}
\end{array}$$
One can readily verify using an adapted co-frame 
$\{\omega^1,\omega^2,\omega^3,\omega^4,\omega^5\}$ with
$$\begin{array}{c}I={\mathrm{span}}\{\omega^1,\omega^2\},\quad
L+I={\mathrm{span}}\{\omega^1,\omega^2,\omega^3\},\\[4pt]
Q+I={\mathrm{span}}\{\omega^1,\omega^2,\omega^4,\omega^5\}
\end{array}$$ 
and such that
\begin{equation}\label{yetanothercoframe}
d\omega^1\equiv\omega^3\wedge\omega^4\bmod\omega^1,\omega^2
\quad\mbox{and}\quad
d\omega^2\equiv\omega^3\wedge\omega^5\bmod\omega^1,\omega^2,
\end{equation}
that the expected cancellations in the $E_0$-level of the associated spectral
sequence actually take place and we have found a differential complex as
follows.
\begin{theorem}\label{protoBGGcomplex}
On any $5$-dimensional manifold equipped with a geometric structure defined by
transverse distributions $D$ and $\ell$ of ranks $2$ and $1$, respectively, and
satisfying {\rm(\ref{anothergeometry})}, there is a canonically defined locally
exact differential complex
$$\begin{array}{r}0\to{\mathbb{R}}\to
\begin{array}{c}
\begin{picture}(42,5)
\put(20,1.5){\makebox(0,0){$\times$}}
\put(36,1.5){\makebox(0,0){$\bullet$}}
\put(4,1.5){\makebox(0,0){$\times$}}
\put(4,1.5){\line(1,0){32}}
\put(4,8){\makebox(0,0){$\scriptstyle 1$}}
\put(20,8){\makebox(0,0){$\scriptstyle -2$}}
\put(36,8){\makebox(0,0){$\scriptstyle 1$}}
\end{picture}\\
\oplus\\[3pt]
\begin{picture}(42,5)
\put(20,1.5){\makebox(0,0){$\times$}}
\put(36,1.5){\makebox(0,0){$\bullet$}}
\put(4,1.5){\makebox(0,0){$\times$}}
\put(4,1.5){\line(1,0){32}}
\put(4,8){\makebox(0,0){$\scriptstyle -2$}}
\put(20,8){\makebox(0,0){$\scriptstyle 1$}}
\put(36,8){\makebox(0,0){$\scriptstyle 0$}}
\end{picture}
\end{array}\!\to\!
\begin{array}{c}
\begin{picture}(42,5)
\put(20,1.5){\makebox(0,0){$\times$}}
\put(36,1.5){\makebox(0,0){$\bullet$}}
\put(4,1.5){\makebox(0,0){$\times$}}
\put(4,1.5){\line(1,0){32}}
\put(4,8){\makebox(0,0){$\scriptstyle 2$}}
\put(20,8){\makebox(0,0){$\scriptstyle -3$}}
\put(36,8){\makebox(0,0){$\scriptstyle 0$}}
\end{picture}\\
+\\[3pt]
\begin{picture}(42,5)
\put(20,1.5){\makebox(0,0){$\times$}}
\put(36,1.5){\makebox(0,0){$\bullet$}}
\put(4,1.5){\makebox(0,0){$\times$}}
\put(4,1.5){\line(1,0){32}}
\put(4,8){\makebox(0,0){$\scriptstyle 0$}}
\put(20,8){\makebox(0,0){$\scriptstyle -3$}}
\put(36,8){\makebox(0,0){$\scriptstyle 2$}}
\end{picture}\oplus
\begin{picture}(42,5)
\put(20,1.5){\makebox(0,0){$\times$}}
\put(36,1.5){\makebox(0,0){$\bullet$}}
\put(4,1.5){\makebox(0,0){$\times$}}
\put(4,1.5){\line(1,0){32}}
\put(4,8){\makebox(0,0){$\scriptstyle -3$}}
\put(20,8){\makebox(0,0){$\scriptstyle 0$}}
\put(36,8){\makebox(0,0){$\scriptstyle 1$}}
\end{picture}
\end{array}\!\!\to\hspace{77pt}\\[30pt]
\begin{array}{c}
\begin{picture}(42,5)
\put(20,1.5){\makebox(0,0){$\times$}}
\put(36,1.5){\makebox(0,0){$\bullet$}}
\put(4,1.5){\makebox(0,0){$\times$}}
\put(4,1.5){\line(1,0){32}}
\put(4,8){\makebox(0,0){$\scriptstyle -2$}}
\put(20,8){\makebox(0,0){$\scriptstyle -2$}}
\put(36,8){\makebox(0,0){$\scriptstyle 2$}}
\end{picture}\oplus
\begin{picture}(42,5)
\put(20,1.5){\makebox(0,0){$\times$}}
\put(36,1.5){\makebox(0,0){$\bullet$}}
\put(4,1.5){\makebox(0,0){$\times$}}
\put(4,1.5){\line(1,0){32}}
\put(4,8){\makebox(0,0){$\scriptstyle 1$}}
\put(20,8){\makebox(0,0){$\scriptstyle -4$}}
\put(36,8){\makebox(0,0){$\scriptstyle 1$}}
\end{picture}\\
+\\[3pt]
\begin{picture}(42,5)
\put(20,1.5){\makebox(0,0){$\times$}}
\put(36,1.5){\makebox(0,0){$\bullet$}}
\put(4,1.5){\makebox(0,0){$\times$}}
\put(4,1.5){\line(1,0){32}}
\put(4,8){\makebox(0,0){$\scriptstyle -4$}}
\put(20,8){\makebox(0,0){$\scriptstyle 0$}}
\put(36,8){\makebox(0,0){$\scriptstyle 0$}}
\end{picture}
\end{array}\!\!\to\!
\begin{array}{c}
\begin{picture}(42,5)
\put(20,1.5){\makebox(0,0){$\times$}}
\put(36,1.5){\makebox(0,0){$\bullet$}}
\put(4,1.5){\makebox(0,0){$\times$}}
\put(4,1.5){\line(1,0){32}}
\put(4,8){\makebox(0,0){$\scriptstyle 0$}}
\put(20,8){\makebox(0,0){$\scriptstyle -4$}}
\put(36,8){\makebox(0,0){$\scriptstyle 0$}}
\end{picture}\\
\oplus\\[3pt]
\begin{picture}(42,5)
\put(20,1.5){\makebox(0,0){$\times$}}
\put(36,1.5){\makebox(0,0){$\bullet$}}
\put(4,1.5){\makebox(0,0){$\times$}}
\put(4,1.5){\line(1,0){32}}
\put(4,8){\makebox(0,0){$\scriptstyle -3$}}
\put(20,8){\makebox(0,0){$\scriptstyle -2$}}
\put(36,8){\makebox(0,0){$\scriptstyle 1$}}
\end{picture}
\end{array}\!\to
\begin{picture}(42,5)
\put(20,1.5){\makebox(0,0){$\times$}}
\put(36,1.5){\makebox(0,0){$\bullet$}}
\put(4,1.5){\makebox(0,0){$\times$}}
\put(4,1.5){\line(1,0){32}}
\put(4,8){\makebox(0,0){$\scriptstyle -2$}}
\put(20,8){\makebox(0,0){$\scriptstyle -3$}}
\put(36,8){\makebox(0,0){$\scriptstyle 0$}}
\end{picture}\to 0.\end{array}$$
\end{theorem}
As in \S\ref{threeinsix}, one can make a further normalisation in order to
split the two bundles that have arisen from the spectral
sequence, or from the equivalent diagram chasing, only as filtered bundles. For
the first of these we note that the freedom in its splitting lies in
$${\mathrm{Hom}}\big(\begin{picture}(42,5)
\put(20,1.5){\makebox(0,0){$\times$}}
\put(36,1.5){\makebox(0,0){$\bullet$}}
\put(4,1.5){\makebox(0,0){$\times$}}
\put(4,1.5){\line(1,0){32}}
\put(4,8){\makebox(0,0){$\scriptstyle 2$}}
\put(20,8){\makebox(0,0){$\scriptstyle -3$}}
\put(36,8){\makebox(0,0){$\scriptstyle 0$}}
\end{picture},\begin{picture}(42,5)
\put(20,1.5){\makebox(0,0){$\times$}}
\put(36,1.5){\makebox(0,0){$\bullet$}}
\put(4,1.5){\makebox(0,0){$\times$}}
\put(4,1.5){\line(1,0){32}}
\put(4,8){\makebox(0,0){$\scriptstyle 0$}}
\put(20,8){\makebox(0,0){$\scriptstyle -3$}}
\put(36,8){\makebox(0,0){$\scriptstyle 2$}}
\end{picture}\oplus\begin{picture}(42,5)
\put(20,1.5){\makebox(0,0){$\times$}}
\put(36,1.5){\makebox(0,0){$\bullet$}}
\put(4,1.5){\makebox(0,0){$\times$}}
\put(4,1.5){\line(1,0){32}}
\put(4,8){\makebox(0,0){$\scriptstyle -3$}}
\put(20,8){\makebox(0,0){$\scriptstyle 0$}}
\put(36,8){\makebox(0,0){$\scriptstyle 1$}}
\end{picture}\big)=
\begin{picture}(42,5)
\put(20,1.5){\makebox(0,0){$\times$}}
\put(36,1.5){\makebox(0,0){$\bullet$}}
\put(4,1.5){\makebox(0,0){$\times$}}
\put(4,1.5){\line(1,0){32}}
\put(4,8){\makebox(0,0){$\scriptstyle -2$}}
\put(20,8){\makebox(0,0){$\scriptstyle 0$}}
\put(36,8){\makebox(0,0){$\scriptstyle 2$}}
\end{picture}\oplus\begin{picture}(42,5)
\put(20,1.5){\makebox(0,0){$\times$}}
\put(36,1.5){\makebox(0,0){$\bullet$}}
\put(4,1.5){\makebox(0,0){$\times$}}
\put(4,1.5){\line(1,0){32}}
\put(4,8){\makebox(0,0){$\scriptstyle -1$}}
\put(20,8){\makebox(0,0){$\scriptstyle 3$}}
\put(36,8){\makebox(0,0){$\scriptstyle 0$}}
\end{picture}$$
whereas, from~(\ref{oneforms}), the freedom in splitting 
$\Lambda^1$ lies in
$$\begin{array}{l}
{\mathrm{Hom}}\big(\begin{picture}(42,5)
\put(20,1.5){\makebox(0,0){$\times$}}
\put(36,1.5){\makebox(0,0){$\bullet$}}
\put(4,1.5){\makebox(0,0){$\times$}}
\put(4,1.5){\line(1,0){32}}
\put(4,8){\makebox(0,0){$\scriptstyle 1$}}
\put(20,8){\makebox(0,0){$\scriptstyle -2$}}
\put(36,8){\makebox(0,0){$\scriptstyle 1$}}
\end{picture}\oplus\begin{picture}(42,5)
\put(20,1.5){\makebox(0,0){$\times$}}
\put(36,1.5){\makebox(0,0){$\bullet$}}
\put(4,1.5){\makebox(0,0){$\times$}}
\put(4,1.5){\line(1,0){32}}
\put(4,8){\makebox(0,0){$\scriptstyle -2$}}
\put(20,8){\makebox(0,0){$\scriptstyle 1$}}
\put(36,8){\makebox(0,0){$\scriptstyle 0$}}
\end{picture},\begin{picture}(42,5)
\put(20,1.5){\makebox(0,0){$\times$}}
\put(36,1.5){\makebox(0,0){$\bullet$}}
\put(4,1.5){\makebox(0,0){$\times$}}
\put(4,1.5){\line(1,0){32}}
\put(4,8){\makebox(0,0){$\scriptstyle -1$}}
\put(20,8){\makebox(0,0){$\scriptstyle -1$}}
\put(36,8){\makebox(0,0){$\scriptstyle 1$}}
\end{picture}\big)\\[6pt]
\qquad{}=\big(\begin{picture}(42,5)
\put(20,1.5){\makebox(0,0){$\times$}}
\put(36,1.5){\makebox(0,0){$\bullet$}}
\put(4,1.5){\makebox(0,0){$\times$}}
\put(4,1.5){\line(1,0){32}}
\put(4,8){\makebox(0,0){$\scriptstyle -1$}}
\put(20,8){\makebox(0,0){$\scriptstyle 1$}}
\put(36,8){\makebox(0,0){$\scriptstyle 1$}}
\end{picture}\oplus\begin{picture}(42,5)
\put(20,1.5){\makebox(0,0){$\times$}}
\put(36,1.5){\makebox(0,0){$\bullet$}}
\put(4,1.5){\makebox(0,0){$\times$}}
\put(4,1.5){\line(1,0){32}}
\put(4,8){\makebox(0,0){$\scriptstyle 2$}}
\put(20,8){\makebox(0,0){$\scriptstyle -1$}}
\put(36,8){\makebox(0,0){$\scriptstyle 0$}}
\end{picture}\big)\otimes\begin{picture}(42,5)
\put(20,1.5){\makebox(0,0){$\times$}}
\put(36,1.5){\makebox(0,0){$\bullet$}}
\put(4,1.5){\makebox(0,0){$\times$}}
\put(4,1.5){\line(1,0){32}}
\put(4,8){\makebox(0,0){$\scriptstyle -1$}}
\put(20,8){\makebox(0,0){$\scriptstyle -1$}}
\put(36,8){\makebox(0,0){$\scriptstyle 1$}}
\end{picture}\\[6pt]
\qquad\qquad{}=\begin{picture}(42,5)
\put(20,1.5){\makebox(0,0){$\times$}}
\put(36,1.5){\makebox(0,0){$\bullet$}}
\put(4,1.5){\makebox(0,0){$\times$}}
\put(4,1.5){\line(1,0){32}}
\put(4,8){\makebox(0,0){$\scriptstyle -2$}}
\put(20,8){\makebox(0,0){$\scriptstyle 0$}}
\put(36,8){\makebox(0,0){$\scriptstyle 2$}}
\end{picture}\oplus\begin{picture}(42,5)
\put(20,1.5){\makebox(0,0){$\times$}}
\put(36,1.5){\makebox(0,0){$\bullet$}}
\put(4,1.5){\makebox(0,0){$\times$}}
\put(4,1.5){\line(1,0){32}}
\put(4,8){\makebox(0,0){$\scriptstyle -2$}}
\put(20,8){\makebox(0,0){$\scriptstyle 1$}}
\put(36,8){\makebox(0,0){$\scriptstyle 0$}}
\end{picture}\oplus\begin{picture}(42,5)
\put(20,1.5){\makebox(0,0){$\times$}}
\put(36,1.5){\makebox(0,0){$\bullet$}}
\put(4,1.5){\makebox(0,0){$\times$}}
\put(4,1.5){\line(1,0){32}}
\put(4,8){\makebox(0,0){$\scriptstyle 1$}}
\put(20,8){\makebox(0,0){$\scriptstyle -2$}}
\put(36,8){\makebox(0,0){$\scriptstyle 1$}}
\end{picture}.\end{array}$$
We see that only $\begin{picture}(42,5)
\put(20,1.5){\makebox(0,0){$\times$}}
\put(36,1.5){\makebox(0,0){$\bullet$}}
\put(4,1.5){\makebox(0,0){$\times$}}
\put(4,1.5){\line(1,0){32}}
\put(4,8){\makebox(0,0){$\scriptstyle -2$}}
\put(20,8){\makebox(0,0){$\scriptstyle 0$}}
\put(36,8){\makebox(0,0){$\scriptstyle 2$}}
\end{picture}$ 
is common to both. Therefore, it is only this freedom that need be eliminated
from the freedom to split $\Lambda^1$. In fact, once this freedom is
eliminated, then Theorem~\ref{protoBGGcomplex} is improved as follows.
\begin{theorem}\label{fullBGGcomplex}
On any $5$-dimensional manifold equipped with a regular
infinitesimal flag structure defined by 
\begin{picture}(42,5)
\put(20,1.5){\makebox(0,0){$\times$}}
\put(36,1.5){\makebox(0,0){$\bullet$}}
\put(4,1.5){\makebox(0,0){$\times$}}
\put(4,1.5){\line(1,0){32}}
\end{picture}, there is a canonically
defined locally exact differential complex
$$\begin{array}{r}
0\to{\mathbb{R}}\to\!
\begin{array}{c}
\begin{picture}(42,5)
\put(20,1.5){\makebox(0,0){$\times$}}
\put(36,1.5){\makebox(0,0){$\bullet$}}
\put(4,1.5){\makebox(0,0){$\times$}}
\put(4,1.5){\line(1,0){32}}
\put(4,8){\makebox(0,0){$\scriptstyle 1$}}
\put(20,8){\makebox(0,0){$\scriptstyle -2$}}
\put(36,8){\makebox(0,0){$\scriptstyle 1$}}
\end{picture}\\
\oplus\\[3pt]
\begin{picture}(42,5)
\put(20,1.5){\makebox(0,0){$\times$}}
\put(36,1.5){\makebox(0,0){$\bullet$}}
\put(4,1.5){\makebox(0,0){$\times$}}
\put(4,1.5){\line(1,0){32}}
\put(4,8){\makebox(0,0){$\scriptstyle -2$}}
\put(20,8){\makebox(0,0){$\scriptstyle 1$}}
\put(36,8){\makebox(0,0){$\scriptstyle 0$}}
\end{picture}
\end{array}\!\to\!
\begin{array}{c}
\begin{picture}(42,5)
\put(20,1.5){\makebox(0,0){$\times$}}
\put(36,1.5){\makebox(0,0){$\bullet$}}
\put(4,1.5){\makebox(0,0){$\times$}}
\put(4,1.5){\line(1,0){32}}
\put(4,8){\makebox(0,0){$\scriptstyle 2$}}
\put(20,8){\makebox(0,0){$\scriptstyle -3$}}
\put(36,8){\makebox(0,0){$\scriptstyle 0$}}
\end{picture}\\
\oplus\\[3pt]
\begin{picture}(42,5)
\put(20,1.5){\makebox(0,0){$\times$}}
\put(36,1.5){\makebox(0,0){$\bullet$}}
\put(4,1.5){\makebox(0,0){$\times$}}
\put(4,1.5){\line(1,0){32}}
\put(4,8){\makebox(0,0){$\scriptstyle 0$}}
\put(20,8){\makebox(0,0){$\scriptstyle -3$}}
\put(36,8){\makebox(0,0){$\scriptstyle 2$}}
\end{picture}\\
\oplus\\[3pt]
\begin{picture}(42,5)
\put(20,1.5){\makebox(0,0){$\times$}}
\put(36,1.5){\makebox(0,0){$\bullet$}}
\put(4,1.5){\makebox(0,0){$\times$}}
\put(4,1.5){\line(1,0){32}}
\put(4,8){\makebox(0,0){$\scriptstyle -3$}}
\put(20,8){\makebox(0,0){$\scriptstyle 0$}}
\put(36,8){\makebox(0,0){$\scriptstyle 1$}}
\end{picture}
\end{array}\!\to\!
\begin{array}{c}
\begin{picture}(42,5)
\put(20,1.5){\makebox(0,0){$\times$}}
\put(36,1.5){\makebox(0,0){$\bullet$}}
\put(4,1.5){\makebox(0,0){$\times$}}
\put(4,1.5){\line(1,0){32}}
\put(4,8){\makebox(0,0){$\scriptstyle 1$}}
\put(20,8){\makebox(0,0){$\scriptstyle -4$}}
\put(36,8){\makebox(0,0){$\scriptstyle 1$}}
\end{picture}\\
\oplus\\[3pt]
\begin{picture}(42,5)
\put(20,1.5){\makebox(0,0){$\times$}}
\put(36,1.5){\makebox(0,0){$\bullet$}}
\put(4,1.5){\makebox(0,0){$\times$}}
\put(4,1.5){\line(1,0){32}}
\put(4,8){\makebox(0,0){$\scriptstyle -2$}}
\put(20,8){\makebox(0,0){$\scriptstyle -2$}}
\put(36,8){\makebox(0,0){$\scriptstyle 2$}}
\end{picture}\\
\oplus\\[3pt]
\begin{picture}(42,5)
\put(20,1.5){\makebox(0,0){$\times$}}
\put(36,1.5){\makebox(0,0){$\bullet$}}
\put(4,1.5){\makebox(0,0){$\times$}}
\put(4,1.5){\line(1,0){32}}
\put(4,8){\makebox(0,0){$\scriptstyle -4$}}
\put(20,8){\makebox(0,0){$\scriptstyle 0$}}
\put(36,8){\makebox(0,0){$\scriptstyle 0$}}
\end{picture}
\end{array}\!\to\!
\begin{array}{c}
\begin{picture}(42,5)
\put(20,1.5){\makebox(0,0){$\times$}}
\put(36,1.5){\makebox(0,0){$\bullet$}}
\put(4,1.5){\makebox(0,0){$\times$}}
\put(4,1.5){\line(1,0){32}}
\put(4,8){\makebox(0,0){$\scriptstyle 0$}}
\put(20,8){\makebox(0,0){$\scriptstyle -4$}}
\put(36,8){\makebox(0,0){$\scriptstyle 0$}}
\end{picture}\\
\oplus\\[3pt]
\begin{picture}(42,5)
\put(20,1.5){\makebox(0,0){$\times$}}
\put(36,1.5){\makebox(0,0){$\bullet$}}
\put(4,1.5){\makebox(0,0){$\times$}}
\put(4,1.5){\line(1,0){32}}
\put(4,8){\makebox(0,0){$\scriptstyle -3$}}
\put(20,8){\makebox(0,0){$\scriptstyle -2$}}
\put(36,8){\makebox(0,0){$\scriptstyle 1$}}
\end{picture}
\end{array}\enskip\quad\\
{}\to\begin{picture}(42,5)
\put(20,1.5){\makebox(0,0){$\times$}}
\put(36,1.5){\makebox(0,0){$\bullet$}}
\put(4,1.5){\makebox(0,0){$\times$}}
\put(4,1.5){\line(1,0){32}}
\put(4,8){\makebox(0,0){$\scriptstyle -2$}}
\put(20,8){\makebox(0,0){$\scriptstyle -3$}}
\put(36,8){\makebox(0,0){$\scriptstyle 0$}}
\end{picture}\to 0.\end{array}$$
\end{theorem}
\begin{proof}
As already remarked, to complete the proof we should find a preferred
class of splittings of the $1$-forms (\ref{oneforms}) so that the
$\begin{picture}(42,5)
\put(20,1.5){\makebox(0,0){$\times$}}
\put(36,1.5){\makebox(0,0){$\bullet$}}
\put(4,1.5){\makebox(0,0){$\times$}}
\put(4,1.5){\line(1,0){32}}
\put(4,8){\makebox(0,0){$\scriptstyle -2$}}
\put(20,8){\makebox(0,0){$\scriptstyle 0$}}
\put(36,8){\makebox(0,0){$\scriptstyle 2$}}
\end{picture}$-freedom present in the general
splitting is eliminated. As in \S\ref{threeinsix}, this can be achieved by
restricting a particular component of the exterior derivative
$d:\Lambda^2\to\Lambda^3$ defined via an arbitrary splitting. In this case, we
may consider the composition
$$\begin{picture}(42,5)
\put(20,1.5){\makebox(0,0){$\times$}}
\put(36,1.5){\makebox(0,0){$\bullet$}}
\put(4,1.5){\makebox(0,0){$\times$}}
\put(4,1.5){\line(1,0){32}}
\put(4,8){\makebox(0,0){$\scriptstyle 0$}}
\put(20,8){\makebox(0,0){$\scriptstyle -2$}}
\put(36,8){\makebox(0,0){$\scriptstyle 0$}}
\end{picture}\to\Lambda^2\xrightarrow{\,d\,}\Lambda^3\to
\begin{picture}(42,5)
\put(20,1.5){\makebox(0,0){$\times$}}
\put(36,1.5){\makebox(0,0){$\bullet$}}
\put(4,1.5){\makebox(0,0){$\times$}}
\put(4,1.5){\line(1,0){32}}
\put(4,8){\makebox(0,0){$\scriptstyle -2$}}
\put(20,8){\makebox(0,0){$\scriptstyle -2$}}
\put(36,8){\makebox(0,0){$\scriptstyle 2$}}
\end{picture}.$$
Using an adapted co-framing (\ref{yetanothercoframe}), one may readily verify
that
\begin{itemize}
\item this is actually a homomorphism of vector bundles,
\item insisting that it vanish reduces the freedom in splitting $\Lambda^1$\\
exactly as desired,
\item this also eliminates the freedom in splitting the filtered\\
occurring bundles in Theorem~\ref{protoBGGcomplex},
\end{itemize}
which completes the proof.\end{proof} The differential complex in
Theorem~\ref{fullBGGcomplex} is our BGG complex for this parabolic geometry.
Furthermore, one can easily check that in case
$[D,D]=\ell\oplus D$, equivalently if the
\begin{picture}(42,5)
\put(20,1.5){\makebox(0,0){$\times$}}
\put(36,1.5){\makebox(0,0){$\bullet$}}
\put(4,1.5){\makebox(0,0){$\times$}}
\put(4,1.5){\line(1,0){32}}
\put(4,8){\makebox(0,0){$\scriptstyle 4$}}
\put(20,8){\makebox(0,0){$\scriptstyle -4$}}
\put(36,8){\makebox(0,0){$\scriptstyle 0$}}
\end{picture}-component of curvature does not vanish, then further
cancellations may be effected and one reduces to the BGG complex for the five
variables geometry previously discussed in~\S\ref{five}.

\section{Pfaffian systems of rank three in seven variables}\label{threeinseven}
Let $M$ be a $7$-manifold endowed with a generic distribution $H\subset TM$ of
rank $4$. Equivalently, let $I\subset\Lambda^1$ be a Pfaffian system of rank
$3$ that is generic in Cartan's sense, meaning that its first derived system
$I'$ is zero. We write the corresponding filtration of the cotangent bundle as
$$\Lambda^1=\Lambda^1_H+I.$$
Genericity says that the Levi form, defined as the composition
$$I\rightarrow\Lambda^1\rightarrow\Lambda^2\rightarrow\Lambda^2_H,$$ 
is injective.

It turns out that there exactly two types of generic rank $4$ distributions in
dimension $7$, corresponding to the two open orbits of the action of
${\mathrm{GL}}(4,\mathbb R)\times{\mathrm{GL}}(3,\mathbb R)$ on the space of
linear maps ${\mathrm{Hom}}(\Lambda^2\mathbb R^4,\mathbb R^3)$ called elliptic,
respectively hyperbolic; see~\cite{M2}. We shall treat these two cases
simultaneously.

The Lie algebra $\mathfrak s\mathfrak p(6,\mathbb C)$ admits a grading of the 
form
$$\begin{array}{l}{\mathfrak{g}}_{-2}\oplus{\mathfrak{g}}_{-1}\oplus
\underbrace{{\mathfrak{g}}_0
\oplus{\mathfrak{g}}_1\oplus{\mathfrak{g}}_2}_{\mathfrak{p}}=\\[-5pt]
\qquad\begin{picture}(42,5)
\put(21.4,.1){\line(1,0){14.3}}
\put(21.4,3.2){\line(1,0){14.3}}
\put(20,1.5){\makebox(0,0){$\times$}}
\put(36,1.6){\makebox(0,0){$\bullet$}}
\put(28,1.5){\makebox(0,0){$\langle$}}
\put(4,1.5){\makebox(0,0){$\bullet$}}
\put(4,1.5){\line(1,0){16}}
\put(4,8){\makebox(0,0){$\scriptstyle 2$}}
\put(20,8){\makebox(0,0){$\scriptstyle 0$}}
\put(36,8){\makebox(0,0){$\scriptstyle 0$}}
\end{picture}\oplus
\begin{picture}(42,5)
\put(21.4,.1){\line(1,0){14.3}}
\put(21.4,3.2){\line(1,0){14.3}}
\put(20,1.5){\makebox(0,0){$\times$}}
\put(36,1.6){\makebox(0,0){$\bullet$}}
\put(28,1.5){\makebox(0,0){$\langle$}}
\put(4,1.5){\makebox(0,0){$\bullet$}}
\put(4,1.5){\line(1,0){16}}
\put(4,8){\makebox(0,0){$\scriptstyle 1$}}
\put(20,8){\makebox(0,0){$\scriptstyle -1$}}
\put(36,8){\makebox(0,0){$\scriptstyle 1$}}
\end{picture}\oplus
\begin{array}{c}
\begin{picture}(42,5)
\put(21.4,.1){\line(1,0){14.3}}
\put(21.4,3.2){\line(1,0){14.3}}
\put(20,1.5){\makebox(0,0){$\times$}}
\put(36,1.6){\makebox(0,0){$\bullet$}}
\put(28,1.5){\makebox(0,0){$\langle$}}
\put(4,1.5){\makebox(0,0){$\bullet$}}
\put(4,1.5){\line(1,0){16}}
\put(4,8){\makebox(0,0){$\scriptstyle 2$}}
\put(20,8){\makebox(0,0){$\scriptstyle -1$}}
\put(36,8){\makebox(0,0){$\scriptstyle 0$}}
\end{picture}\\[-1pt]
\oplus\\[2pt]
\begin{picture}(42,5)
\put(21.4,.1){\line(1,0){14.3}}
\put(21.4,3.2){\line(1,0){14.3}}
\put(20,1.5){\makebox(0,0){$\times$}}
\put(36,1.6){\makebox(0,0){$\bullet$}}
\put(28,1.5){\makebox(0,0){$\langle$}}
\put(4,1.5){\makebox(0,0){$\bullet$}}
\put(4,1.5){\line(1,0){16}}
\put(4,8){\makebox(0,0){$\scriptstyle 0$}}
\put(20,8){\makebox(0,0){$\scriptstyle -2$}}
\put(36,8){\makebox(0,0){$\scriptstyle 2$}}
\end{picture}\\[-1pt]
\oplus\\[2pt]
\begin{picture}(42,5)
\put(21.4,.1){\line(1,0){14.3}}
\put(21.4,3.2){\line(1,0){14.3}}
\put(20,1.5){\makebox(0,0){$\times$}}
\put(36,1.6){\makebox(0,0){$\bullet$}}
\put(28,1.5){\makebox(0,0){$\langle$}}
\put(4,1.5){\makebox(0,0){$\bullet$}}
\put(4,1.5){\line(1,0){16}}
\put(4,8){\makebox(0,0){$\scriptstyle 0$}}
\put(20,8){\makebox(0,0){$\scriptstyle 0$}}
\put(36,8){\makebox(0,0){$\scriptstyle 0$}}
\end{picture}
\end{array}\oplus
\begin{picture}(42,5)
\put(21.4,.1){\line(1,0){14.3}}
\put(21.4,3.2){\line(1,0){14.3}}
\put(20,1.5){\makebox(0,0){$\times$}}
\put(36,1.6){\makebox(0,0){$\bullet$}}
\put(28,1.5){\makebox(0,0){$\langle$}}
\put(4,1.5){\makebox(0,0){$\bullet$}}
\put(4,1.5){\line(1,0){16}}
\put(4,8){\makebox(0,0){$\scriptstyle 1$}}
\put(20,8){\makebox(0,0){$\scriptstyle -2$}}
\put(36,8){\makebox(0,0){$\scriptstyle 1$}}
\end{picture}\oplus
\begin{picture}(42,5)
\put(21.4,.1){\line(1,0){14.3}}
\put(21.4,3.2){\line(1,0){14.3}}
\put(20,1.5){\makebox(0,0){$\times$}}
\put(36,1.6){\makebox(0,0){$\bullet$}}
\put(28,1.5){\makebox(0,0){$\langle$}}
\put(4,1.5){\makebox(0,0){$\bullet$}}
\put(4,1.5){\line(1,0){16}}
\put(4,8){\makebox(0,0){$\scriptstyle 2$}}
\put(20,8){\makebox(0,0){$\scriptstyle -2$}}
\put(36,8){\makebox(0,0){$\scriptstyle 0$}}
\end{picture}.\end{array}$$
There are two real forms of this grading, namely $\mathfrak s\mathfrak p(2,1)$
and the split real form $\mathfrak s\mathfrak p(6,\mathbb R)$. One can see that
these gradings give rise to an elliptic generic rank $4$ distribution on the
corresponding $7$-dimensional homogeneous space $G/P$ in the first case and to
a hyperbolic generic rank $4$ distribution on the corresponding homogeneous
space $G/P$ in the second case. The parabolic geometries based on these
particular $G/P$ are known as {\em quaternionic contact\/} \cite[\S4.3.3]{CSl}
and {\em split quaternionic contact\/} \cite[\S4.3.4]{CSl}, respectively.
Regular infinitesimal flag structures of these types correspond exactly to
generic rank $4$ distributions on $7$-manifolds and the irreducible bundles of
these geometries can be written as
$$
\begin{picture}(42,5)
\put(21.4,.1){\line(1,0){14.3}}
\put(21.4,3.2){\line(1,0){14.3}}
\put(20,1.5){\makebox(0,0){$\times$}}
\put(36,1.6){\makebox(0,0){$\bullet$}}
\put(28,1.5){\makebox(0,0){$\langle$}}
\put(4,1.5){\makebox(0,0){$\bullet$}}
\put(4,1.5){\line(1,0){16}}
\put(4,8){\makebox(0,0){$\scriptstyle a$}}
\put(20,8){\makebox(0,0){$\scriptstyle b$}}
\put(36,8){\makebox(0,0){$\scriptstyle c$}}
\end{picture},
$$
where $a,c$ are non-negative integers. Accordingly, we can write the filtration
of the cotangent bundle of a generic rank $4$ distribution as
$$\Lambda^1=
\begin{picture}(42,5)
\put(21.4,.1){\line(1,0){14.3}}
\put(21.4,3.2){\line(1,0){14.3}}
\put(20,1.5){\makebox(0,0){$\times$}}
\put(36,1.6){\makebox(0,0){$\bullet$}}
\put(28,1.5){\makebox(0,0){$\langle$}}
\put(4,1.5){\makebox(0,0){$\bullet$}}
\put(4,1.5){\line(1,0){16}}
\put(4,8){\makebox(0,0){$\scriptstyle 1$}}
\put(20,8){\makebox(0,0){$\scriptstyle -2$}}
\put(36,8){\makebox(0,0){$\scriptstyle 1$}}
\end{picture}+
\begin{picture}(42,5)
\put(21.4,.1){\line(1,0){14.3}}
\put(21.4,3.2){\line(1,0){14.3}}
\put(20,1.5){\makebox(0,0){$\times$}}
\put(36,1.6){\makebox(0,0){$\bullet$}}
\put(28,1.5){\makebox(0,0){$\langle$}}
\put(4,1.5){\makebox(0,0){$\bullet$}}
\put(4,1.5){\line(1,0){16}}
\put(4,8){\makebox(0,0){$\scriptstyle 2$}}
\put(20,8){\makebox(0,0){$\scriptstyle -2$}}
\put(36,8){\makebox(0,0){$\scriptstyle 0$}}
\end{picture}
$$ 
and the filtration of the higher forms as
$$\Lambda^2=
\begin{array}{c}
\begin{picture}(42,5)
\put(21.4,.1){\line(1,0){14.3}}
\put(21.4,3.2){\line(1,0){14.3}}
\put(20,1.5){\makebox(0,0){$\times$}}
\put(36,1.6){\makebox(0,0){$\bullet$}}
\put(28,1.5){\makebox(0,0){$\langle$}}
\put(4,1.5){\makebox(0,0){$\bullet$}}
\put(4,1.5){\line(1,0){16}}
\put(4,8){\makebox(0,0){$\scriptstyle 2$}}
\put(20,8){\makebox(0,0){$\scriptstyle -2$}}
\put(36,8){\makebox(0,0){$\scriptstyle 0$}}
\end{picture}\\[-1pt]
\oplus\\[2pt]
\begin{picture}(42,5)
\put(21.4,.1){\line(1,0){14.3}}
\put(21.4,3.2){\line(1,0){14.3}}
\put(20,1.5){\makebox(0,0){$\times$}}
\put(36,1.6){\makebox(0,0){$\bullet$}}
\put(28,1.5){\makebox(0,0){$\langle$}}
\put(4,1.5){\makebox(0,0){$\bullet$}}
\put(4,1.5){\line(1,0){16}}
\put(4,8){\makebox(0,0){$\scriptstyle 0$}}
\put(20,8){\makebox(0,0){$\scriptstyle -3$}}
\put(36,8){\makebox(0,0){$\scriptstyle 2$}}
\end{picture}
\end{array}
+
\begin{array}{c}
\begin{picture}(42,5)
\put(21.4,.1){\line(1,0){14.3}}
\put(21.4,3.2){\line(1,0){14.3}}
\put(20,1.5){\makebox(0,0){$\times$}}
\put(36,1.6){\makebox(0,0){$\bullet$}}
\put(28,1.5){\makebox(0,0){$\langle$}}
\put(4,1.5){\makebox(0,0){$\bullet$}}
\put(4,1.5){\line(1,0){16}}
\put(4,8){\makebox(0,0){$\scriptstyle 3$}}
\put(20,8){\makebox(0,0){$\scriptstyle -4$}}
\put(36,8){\makebox(0,0){$\scriptstyle 1$}}
\end{picture}\\[-1pt]
\oplus\\[2pt]
\begin{picture}(42,5)
\put(21.4,.1){\line(1,0){14.3}}
\put(21.4,3.2){\line(1,0){14.3}}
\put(20,1.5){\makebox(0,0){$\times$}}
\put(36,1.6){\makebox(0,0){$\bullet$}}
\put(28,1.5){\makebox(0,0){$\langle$}}
\put(4,1.5){\makebox(0,0){$\bullet$}}
\put(4,1.5){\line(1,0){16}}
\put(4,8){\makebox(0,0){$\scriptstyle 1$}}
\put(20,8){\makebox(0,0){$\scriptstyle -3$}}
\put(36,8){\makebox(0,0){$\scriptstyle 1$}}
\end{picture}
\end{array}
+
\begin{picture}(42,5)
\put(21.4,.1){\line(1,0){14.3}}
\put(21.4,3.2){\line(1,0){14.3}}
\put(20,1.5){\makebox(0,0){$\times$}}
\put(36,1.6){\makebox(0,0){$\bullet$}}
\put(28,1.5){\makebox(0,0){$\langle$}}
\put(4,1.5){\makebox(0,0){$\bullet$}}
\put(4,1.5){\line(1,0){16}}
\put(4,8){\makebox(0,0){$\scriptstyle 2$}}
\put(20,8){\makebox(0,0){$\scriptstyle -3$}}
\put(36,8){\makebox(0,0){$\scriptstyle 0$}}
\end{picture}
$$ 
\mbox{ }
$$\Lambda^3 =
\begin{picture}(42,5)
\put(21.4,.1){\line(1,0){14.3}}
\put(21.4,3.2){\line(1,0){14.3}}
\put(20,1.5){\makebox(0,0){$\times$}}
\put(36,1.6){\makebox(0,0){$\bullet$}}
\put(28,1.5){\makebox(0,0){$\langle$}}
\put(4,1.5){\makebox(0,0){$\bullet$}}
\put(4,1.5){\line(1,0){16}}
\put(4,8){\makebox(0,0){$\scriptstyle 1$}}
\put(20,8){\makebox(0,0){$\scriptstyle -3$}}
\put(36,8){\makebox(0,0){$\scriptstyle 1$}}
\end{picture}
+
\begin{array}{c}
\begin{picture}(42,5)
\put(21.4,.1){\line(1,0){14.3}}
\put(21.4,3.2){\line(1,0){14.3}}
\put(20,1.5){\makebox(0,0){$\times$}}
\put(36,1.6){\makebox(0,0){$\bullet$}}
\put(28,1.5){\makebox(0,0){$\langle$}}
\put(4,1.5){\makebox(0,0){$\bullet$}}
\put(4,1.5){\line(1,0){16}}
\put(4,8){\makebox(0,0){$\scriptstyle 4$}}
\put(20,8){\makebox(0,0){$\scriptstyle -4$}}
\put(36,8){\makebox(0,0){$\scriptstyle 0$}}
\end{picture}\\[-1pt]
\oplus\\[2pt]
\begin{picture}(42,5)
\put(21.4,.1){\line(1,0){14.3}}
\put(21.4,3.2){\line(1,0){14.3}}
\put(20,1.5){\makebox(0,0){$\times$}}
\put(36,1.6){\makebox(0,0){$\bullet$}}
\put(28,1.5){\makebox(0,0){$\langle$}}
\put(4,1.5){\makebox(0,0){$\bullet$}}
\put(4,1.5){\line(1,0){16}}
\put(4,8){\makebox(0,0){$\scriptstyle 2$}}
\put(20,8){\makebox(0,0){$\scriptstyle -3$}}
\put(36,8){\makebox(0,0){$\scriptstyle 0$}}
\end{picture}
\\[-1pt]
\oplus\\[2pt]
\begin{picture}(42,5)
\put(21.4,.1){\line(1,0){14.3}}
\put(21.4,3.2){\line(1,0){14.3}}
\put(20,1.5){\makebox(0,0){$\times$}}
\put(36,1.6){\makebox(0,0){$\bullet$}}
\put(28,1.5){\makebox(0,0){$\langle$}}
\put(4,1.5){\makebox(0,0){$\bullet$}}
\put(4,1.5){\line(1,0){16}}
\put(4,8){\makebox(0,0){$\scriptstyle 0$}}
\put(20,8){\makebox(0,0){$\scriptstyle -2$}}
\put(36,8){\makebox(0,0){$\scriptstyle 0$}}
\end{picture}
\\[-1pt]
\oplus\\[2pt]
\begin{picture}(42,5)
\put(21.4,.1){\line(1,0){14.3}}
\put(21.4,3.2){\line(1,0){14.3}}
\put(20,1.5){\makebox(0,0){$\times$}}
\put(36,1.6){\makebox(0,0){$\bullet$}}
\put(28,1.5){\makebox(0,0){$\langle$}}
\put(4,1.5){\makebox(0,0){$\bullet$}}
\put(4,1.5){\line(1,0){16}}
\put(4,8){\makebox(0,0){$\scriptstyle 2$}}
\put(20,8){\makebox(0,0){$\scriptstyle -5$}}
\put(36,8){\makebox(0,0){$\scriptstyle 2$}}
\end{picture}
\end{array}
+
\begin{array}{c}
\begin{picture}(42,5)
\put(21.4,.1){\line(1,0){14.3}}
\put(21.4,3.2){\line(1,0){14.3}}
\put(20,1.5){\makebox(0,0){$\times$}}
\put(36,1.6){\makebox(0,0){$\bullet$}}
\put(28,1.5){\makebox(0,0){$\langle$}}
\put(4,1.5){\makebox(0,0){$\bullet$}}
\put(4,1.5){\line(1,0){16}}
\put(4,8){\makebox(0,0){$\scriptstyle 3$}}
\put(20,8){\makebox(0,0){$\scriptstyle -5$}}
\put(36,8){\makebox(0,0){$\scriptstyle 1$}}
\end{picture}\\[-1pt]
\oplus\\[2pt]
\begin{picture}(42,5)
\put(21.4,.1){\line(1,0){14.3}}
\put(21.4,3.2){\line(1,0){14.3}}
\put(20,1.5){\makebox(0,0){$\times$}}
\put(36,1.6){\makebox(0,0){$\bullet$}}
\put(28,1.5){\makebox(0,0){$\langle$}}
\put(4,1.5){\makebox(0,0){$\bullet$}}
\put(4,1.5){\line(1,0){16}}
\put(4,8){\makebox(0,0){$\scriptstyle 1$}}
\put(20,8){\makebox(0,0){$\scriptstyle -4$}}
\put(36,8){\makebox(0,0){$\scriptstyle 1$}}
\end{picture}
\end{array}
+
\begin{picture}(42,5)
\put(21.4,.1){\line(1,0){14.3}}
\put(21.4,3.2){\line(1,0){14.3}}
\put(20,1.5){\makebox(0,0){$\times$}}
\put(36,1.6){\makebox(0,0){$\bullet$}}
\put(28,1.5){\makebox(0,0){$\langle$}}
\put(4,1.5){\makebox(0,0){$\bullet$}}
\put(4,1.5){\line(1,0){16}}
\put(4,8){\makebox(0,0){$\scriptstyle 0$}}
\put(20,8){\makebox(0,0){$\scriptstyle -3$}}
\put(36,8){\makebox(0,0){$\scriptstyle 0$}}
\end{picture}
$$ 
\mbox{ }
$$
\Lambda^4 =
\begin{picture}(42,5)
\put(21.4,.1){\line(1,0){14.3}}
\put(21.4,3.2){\line(1,0){14.3}}
\put(20,1.5){\makebox(0,0){$\times$}}
\put(36,1.6){\makebox(0,0){$\bullet$}}
\put(28,1.5){\makebox(0,0){$\langle$}}
\put(4,1.5){\makebox(0,0){$\bullet$}}
\put(4,1.5){\line(1,0){16}}
\put(4,8){\makebox(0,0){$\scriptstyle 0$}}
\put(20,8){\makebox(0,0){$\scriptstyle -2$}}
\put(36,8){\makebox(0,0){$\scriptstyle 0$}}
\end{picture}
+
\begin{array}{c}
\begin{picture}(42,5)
\put(21.4,.1){\line(1,0){14.3}}
\put(21.4,3.2){\line(1,0){14.3}}
\put(20,1.5){\makebox(0,0){$\times$}}
\put(36,1.6){\makebox(0,0){$\bullet$}}
\put(28,1.5){\makebox(0,0){$\langle$}}
\put(4,1.5){\makebox(0,0){$\bullet$}}
\put(4,1.5){\line(1,0){16}}
\put(4,8){\makebox(0,0){$\scriptstyle 3$}}
\put(20,8){\makebox(0,0){$\scriptstyle -5$}}
\put(36,8){\makebox(0,0){$\scriptstyle 1$}}
\end{picture}\\[-1pt]
\oplus\\[2pt]
\begin{picture}(42,5)
\put(21.4,.1){\line(1,0){14.3}}
\put(21.4,3.2){\line(1,0){14.3}}
\put(20,1.5){\makebox(0,0){$\times$}}
\put(36,1.6){\makebox(0,0){$\bullet$}}
\put(28,1.5){\makebox(0,0){$\langle$}}
\put(4,1.5){\makebox(0,0){$\bullet$}}
\put(4,1.5){\line(1,0){16}}
\put(4,8){\makebox(0,0){$\scriptstyle 1$}}
\put(20,8){\makebox(0,0){$\scriptstyle -4$}}
\put(36,8){\makebox(0,0){$\scriptstyle 1$}}
\end{picture}
\end{array}
+
\begin{array}{c}
\begin{picture}(42,5)
\put(21.4,.1){\line(1,0){14.3}}
\put(21.4,3.2){\line(1,0){14.3}}
\put(20,1.5){\makebox(0,0){$\times$}}
\put(36,1.6){\makebox(0,0){$\bullet$}}
\put(28,1.5){\makebox(0,0){$\langle$}}
\put(4,1.5){\makebox(0,0){$\bullet$}}
\put(4,1.5){\line(1,0){16}}
\put(4,8){\makebox(0,0){$\scriptstyle 4$}}
\put(20,8){\makebox(0,0){$\scriptstyle -5$}}
\put(36,8){\makebox(0,0){$\scriptstyle 0$}}
\end{picture}\\[-1pt]
\oplus\\[2pt]
\begin{picture}(42,5)
\put(21.4,.1){\line(1,0){14.3}}
\put(21.4,3.2){\line(1,0){14.3}}
\put(20,1.5){\makebox(0,0){$\times$}}
\put(36,1.6){\makebox(0,0){$\bullet$}}
\put(28,1.5){\makebox(0,0){$\langle$}}
\put(4,1.5){\makebox(0,0){$\bullet$}}
\put(4,1.5){\line(1,0){16}}
\put(4,8){\makebox(0,0){$\scriptstyle 2$}}
\put(20,8){\makebox(0,0){$\scriptstyle -4$}}
\put(36,8){\makebox(0,0){$\scriptstyle 0$}}
\end{picture}\\[-1pt]
\oplus\\[2pt]
\begin{picture}(42,5)
\put(21.4,.1){\line(1,0){14.3}}
\put(21.4,3.2){\line(1,0){14.3}}
\put(20,1.5){\makebox(0,0){$\times$}}
\put(36,1.6){\makebox(0,0){$\bullet$}}
\put(28,1.5){\makebox(0,0){$\langle$}}
\put(4,1.5){\makebox(0,0){$\bullet$}}
\put(4,1.5){\line(1,0){16}}
\put(4,8){\makebox(0,0){$\scriptstyle 0$}}
\put(20,8){\makebox(0,0){$\scriptstyle -3$}}
\put(36,8){\makebox(0,0){$\scriptstyle 0$}}
\end{picture}\\[-1pt]
\oplus\\[2pt]
\begin{picture}(42,5)
\put(21.4,.1){\line(1,0){14.3}}
\put(21.4,3.2){\line(1,0){14.3}}
\put(20,1.5){\makebox(0,0){$\times$}}
\put(36,1.6){\makebox(0,0){$\bullet$}}
\put(28,1.5){\makebox(0,0){$\langle$}}
\put(4,1.5){\makebox(0,0){$\bullet$}}
\put(4,1.5){\line(1,0){16}}
\put(4,8){\makebox(0,0){$\scriptstyle 2$}}
\put(20,8){\makebox(0,0){$\scriptstyle -6$}}
\put(36,8){\makebox(0,0){$\scriptstyle 2$}}
\end{picture}
\end{array}
+
\begin{picture}(42,5)
\put(21.4,.1){\line(1,0){14.3}}
\put(21.4,3.2){\line(1,0){14.3}}
\put(20,1.5){\makebox(0,0){$\times$}}
\put(36,1.6){\makebox(0,0){$\bullet$}}
\put(28,1.5){\makebox(0,0){$\langle$}}
\put(4,1.5){\makebox(0,0){$\bullet$}}
\put(4,1.5){\line(1,0){16}}
\put(4,8){\makebox(0,0){$\scriptstyle 1$}}
\put(20,8){\makebox(0,0){$\scriptstyle -5$}}
\put(36,8){\makebox(0,0){$\scriptstyle 1$}}
\end{picture}
$$
\mbox{ }
$$
\Lambda^5 =
\begin{picture}(42,5)
\put(21.4,.1){\line(1,0){14.3}}
\put(21.4,3.2){\line(1,0){14.3}}
\put(20,1.5){\makebox(0,0){$\times$}}
\put(36,1.6){\makebox(0,0){$\bullet$}}
\put(28,1.5){\makebox(0,0){$\langle$}}
\put(4,1.5){\makebox(0,0){$\bullet$}}
\put(4,1.5){\line(1,0){16}}
\put(4,8){\makebox(0,0){$\scriptstyle 2$}}
\put(20,8){\makebox(0,0){$\scriptstyle -4$}}
\put(36,8){\makebox(0,0){$\scriptstyle 0$}}
\end{picture}
+
\begin{array}{c}
\begin{picture}(42,5)
\put(21.4,.1){\line(1,0){14.3}}
\put(21.4,3.2){\line(1,0){14.3}}
\put(20,1.5){\makebox(0,0){$\times$}}
\put(36,1.6){\makebox(0,0){$\bullet$}}
\put(28,1.5){\makebox(0,0){$\langle$}}
\put(4,1.5){\makebox(0,0){$\bullet$}}
\put(4,1.5){\line(1,0){16}}
\put(4,8){\makebox(0,0){$\scriptstyle 3$}}
\put(20,8){\makebox(0,0){$\scriptstyle -6$}}
\put(36,8){\makebox(0,0){$\scriptstyle 1$}}
\end{picture}\\[-1pt]
\oplus\\[2pt]
\begin{picture}(42,5)
\put(21.4,.1){\line(1,0){14.3}}
\put(21.4,3.2){\line(1,0){14.3}}
\put(20,1.5){\makebox(0,0){$\times$}}
\put(36,1.6){\makebox(0,0){$\bullet$}}
\put(28,1.5){\makebox(0,0){$\langle$}}
\put(4,1.5){\makebox(0,0){$\bullet$}}
\put(4,1.5){\line(1,0){16}}
\put(4,8){\makebox(0,0){$\scriptstyle 1$}}
\put(20,8){\makebox(0,0){$\scriptstyle -5$}}
\put(36,8){\makebox(0,0){$\scriptstyle 1$}}
\end{picture}
\end{array}
+
\begin{array}{c}
\begin{picture}(42,5)
\put(21.4,.1){\line(1,0){14.3}}
\put(21.4,3.2){\line(1,0){14.3}}
\put(20,1.5){\makebox(0,0){$\times$}}
\put(36,1.6){\makebox(0,0){$\bullet$}}
\put(28,1.5){\makebox(0,0){$\langle$}}
\put(4,1.5){\makebox(0,0){$\bullet$}}
\put(4,1.5){\line(1,0){16}}
\put(4,8){\makebox(0,0){$\scriptstyle 2$}}
\put(20,8){\makebox(0,0){$\scriptstyle -5$}}
\put(36,8){\makebox(0,0){$\scriptstyle 0$}}
\end{picture}\\[-1pt]
\oplus\\[2pt]
\begin{picture}(42,5)
\put(21.4,.1){\line(1,0){14.3}}
\put(21.4,3.2){\line(1,0){14.3}}
\put(20,1.5){\makebox(0,0){$\times$}}
\put(36,1.6){\makebox(0,0){$\bullet$}}
\put(28,1.5){\makebox(0,0){$\langle$}}
\put(4,1.5){\makebox(0,0){$\bullet$}}
\put(4,1.5){\line(1,0){16}}
\put(4,8){\makebox(0,0){$\scriptstyle 0$}}
\put(20,8){\makebox(0,0){$\scriptstyle -6$}}
\put(36,8){\makebox(0,0){$\scriptstyle 2$}}
\end{picture}
\end{array}
$$
\mbox{ }
$$
\Lambda^6=
\begin{picture}(42,5)
\put(21.4,.1){\line(1,0){14.3}}
\put(21.4,3.2){\line(1,0){14.3}}
\put(20,1.5){\makebox(0,0){$\times$}}
\put(36,1.6){\makebox(0,0){$\bullet$}}
\put(28,1.5){\makebox(0,0){$\langle$}}
\put(4,1.5){\makebox(0,0){$\bullet$}}
\put(4,1.5){\line(1,0){16}}
\put(4,8){\makebox(0,0){$\scriptstyle 2$}}
\put(20,8){\makebox(0,0){$\scriptstyle -5$}}
\put(36,8){\makebox(0,0){$\scriptstyle 0$}}
\end{picture}
+
\begin{picture}(42,5)
\put(21.4,.1){\line(1,0){14.3}}
\put(21.4,3.2){\line(1,0){14.3}}
\put(20,1.5){\makebox(0,0){$\times$}}
\put(36,1.6){\makebox(0,0){$\bullet$}}
\put(28,1.5){\makebox(0,0){$\langle$}}
\put(4,1.5){\makebox(0,0){$\bullet$}}
\put(4,1.5){\line(1,0){16}}
\put(4,8){\makebox(0,0){$\scriptstyle 1$}}
\put(20,8){\makebox(0,0){$\scriptstyle -6$}}
\put(36,8){\makebox(0,0){$\scriptstyle 1$}}
\end{picture}
\quad\quad\quad\quad\quad
\Lambda^7= 
\begin{picture}(42,5)
\put(21.4,.1){\line(1,0){14.3}}
\put(21.4,3.2){\line(1,0){14.3}}
\put(20,1.5){\makebox(0,0){$\times$}}
\put(36,1.6){\makebox(0,0){$\bullet$}}
\put(28,1.5){\makebox(0,0){$\langle$}}
\put(4,1.5){\makebox(0,0){$\bullet$}}
\put(4,1.5){\line(1,0){16}}
\put(4,8){\makebox(0,0){$\scriptstyle 0$}}
\put(20,8){\makebox(0,0){$\scriptstyle -5$}}
\put(36,8){\makebox(0,0){$\scriptstyle 0$}}
\end{picture}
$$
Choosing an adapted co-framing of the Pfaffian system, one can explicitly
verify that all the expected cancellations at the $E_0$-level of the associated
spectral sequence take place and, therefore, one obtains the following.
\begin{theorem}
There is a canonically defined locally exact differential complex
$$\begin{array}{l}
0\rightarrow{\mathbb{R}}\rightarrow 
\begin{picture}(42,5)
\put(21.4,.1){\line(1,0){14.3}}
\put(21.4,3.2){\line(1,0){14.3}}
\put(20,1.5){\makebox(0,0){$\times$}}
\put(36,1.6){\makebox(0,0){$\bullet$}}
\put(28,1.5){\makebox(0,0){$\langle$}}
\put(4,1.5){\makebox(0,0){$\bullet$}}
\put(4,1.5){\line(1,0){16}}
\put(4,8){\makebox(0,0){$\scriptstyle 1$}}
\put(20,8){\makebox(0,0){$\scriptstyle -2$}}
\put(36,8){\makebox(0,0){$\scriptstyle 1$}}
\end{picture}
\rightarrow\!\!
\begin{array}{c}
\begin{picture}(42,5)
\put(21.4,.1){\line(1,0){14.3}}
\put(21.4,3.2){\line(1,0){14.3}}
\put(20,1.5){\makebox(0,0){$\times$}}
\put(36,1.6){\makebox(0,0){$\bullet$}}
\put(28,1.5){\makebox(0,0){$\langle$}}
\put(4,1.5){\makebox(0,0){$\bullet$}}
\put(4,1.5){\line(1,0){16}}
\put(4,8){\makebox(0,0){$\scriptstyle 0$}}
\put(20,8){\makebox(0,0){$\scriptstyle -3$}}
\put(36,8){\makebox(0,0){$\scriptstyle 2$}}
\end{picture}\\
+\\[3pt]
\begin{picture}(42,5)
\put(21.4,.1){\line(1,0){14.3}}
\put(21.4,3.2){\line(1,0){14.3}}
\put(20,1.5){\makebox(0,0){$\times$}}
\put(36,1.6){\makebox(0,0){$\bullet$}}
\put(28,1.5){\makebox(0,0){$\langle$}}
\put(4,1.5){\makebox(0,0){$\bullet$}}
\put(4,1.5){\line(1,0){16}}
\put(4,8){\makebox(0,0){$\scriptstyle 3$}}
\put(20,8){\makebox(0,0){$\scriptstyle -4$}}
\put(36,8){\makebox(0,0){$\scriptstyle 1$}}
\end{picture}
\end{array}
\!\!\rightarrow\!\!
\begin{array}{c}
\begin{picture}(42,5)
\put(21.4,.1){\line(1,0){14.3}}
\put(21.4,3.2){\line(1,0){14.3}}
\put(20,1.5){\makebox(0,0){$\times$}}
\put(36,1.6){\makebox(0,0){$\bullet$}}
\put(28,1.5){\makebox(0,0){$\langle$}}
\put(4,1.5){\makebox(0,0){$\bullet$}}
\put(4,1.5){\line(1,0){16}}
\put(4,8){\makebox(0,0){$\scriptstyle 2$}}
\put(20,8){\makebox(0,0){$\scriptstyle -5$}}
\put(36,8){\makebox(0,0){$\scriptstyle 2$}}
\end{picture}\\
\oplus\\[3pt]
\begin{picture}(42,5)
\put(21.4,.1){\line(1,0){14.3}}
\put(21.4,3.2){\line(1,0){14.3}}
\put(20,1.5){\makebox(0,0){$\times$}}
\put(36,1.6){\makebox(0,0){$\bullet$}}
\put(28,1.5){\makebox(0,0){$\langle$}}
\put(4,1.5){\makebox(0,0){$\bullet$}}
\put(4,1.5){\line(1,0){16}}
\put(4,8){\makebox(0,0){$\scriptstyle 4$}}
\put(20,8){\makebox(0,0){$\scriptstyle -4$}}
\put(36,8){\makebox(0,0){$\scriptstyle 0$}}
\end{picture}
\end{array}
\!\!\rightarrow\!\!
\begin{array}{c}
\begin{picture}(42,5)
\put(21.4,.1){\line(1,0){14.3}}
\put(21.4,3.2){\line(1,0){14.3}}
\put(20,1.5){\makebox(0,0){$\times$}}
\put(36,1.6){\makebox(0,0){$\bullet$}}
\put(28,1.5){\makebox(0,0){$\langle$}}
\put(4,1.5){\makebox(0,0){$\bullet$}}
\put(4,1.5){\line(1,0){16}}
\put(4,8){\makebox(0,0){$\scriptstyle 4$}}
\put(20,8){\makebox(0,0){$\scriptstyle -5$}}
\put(36,8){\makebox(0,0){$\scriptstyle 0$}}
\end{picture}\\
\oplus\\[3pt]
\begin{picture}(42,5)
\put(21.4,.1){\line(1,0){14.3}}
\put(21.4,3.2){\line(1,0){14.3}}
\put(20,1.5){\makebox(0,0){$\times$}}
\put(36,1.6){\makebox(0,0){$\bullet$}}
\put(28,1.5){\makebox(0,0){$\langle$}}
\put(4,1.5){\makebox(0,0){$\bullet$}}
\put(4,1.5){\line(1,0){16}}
\put(4,8){\makebox(0,0){$\scriptstyle 2$}}
\put(20,8){\makebox(0,0){$\scriptstyle -6$}}
\put(36,8){\makebox(0,0){$\scriptstyle 2$}}
\end{picture}
\end{array}
\!\!\rightarrow\\
\qquad\begin{array}{c}
\begin{picture}(42,5)
\put(21.4,.1){\line(1,0){14.3}}
\put(21.4,3.2){\line(1,0){14.3}}
\put(20,1.5){\makebox(0,0){$\times$}}
\put(36,1.6){\makebox(0,0){$\bullet$}}
\put(28,1.5){\makebox(0,0){$\langle$}}
\put(4,1.5){\makebox(0,0){$\bullet$}}
\put(4,1.5){\line(1,0){16}}
\put(4,8){\makebox(0,0){$\scriptstyle 3$}}
\put(20,8){\makebox(0,0){$\scriptstyle -6$}}
\put(36,8){\makebox(0,0){$\scriptstyle 1$}}
\end{picture}\\
+\\[3pt]
\begin{picture}(42,5)
\put(21.4,.1){\line(1,0){14.3}}
\put(21.4,3.2){\line(1,0){14.3}}
\put(20,1.5){\makebox(0,0){$\times$}}
\put(36,1.6){\makebox(0,0){$\bullet$}}
\put(28,1.5){\makebox(0,0){$\langle$}}
\put(4,1.5){\makebox(0,0){$\bullet$}}
\put(4,1.5){\line(1,0){16}}
\put(4,8){\makebox(0,0){$\scriptstyle 0$}}
\put(20,8){\makebox(0,0){$\scriptstyle -6$}}
\put(36,8){\makebox(0,0){$\scriptstyle 2$}}
\end{picture}
\end{array}
\!\!\rightarrow
\begin{picture}(42,5)
\put(21.4,.1){\line(1,0){14.3}}
\put(21.4,3.2){\line(1,0){14.3}}
\put(20,1.5){\makebox(0,0){$\times$}}
\put(36,1.6){\makebox(0,0){$\bullet$}}
\put(28,1.5){\makebox(0,0){$\langle$}}
\put(4,1.5){\makebox(0,0){$\bullet$}}
\put(4,1.5){\line(1,0){16}}
\put(4,8){\makebox(0,0){$\scriptstyle 1$}}
\put(20,8){\makebox(0,0){$\scriptstyle -6$}}
\put(36,8){\makebox(0,0){$\scriptstyle 1$}}
\end{picture}
\rightarrow
\begin{picture}(42,5)
\put(21.4,.1){\line(1,0){14.3}}
\put(21.4,3.2){\line(1,0){14.3}}
\put(20,1.5){\makebox(0,0){$\times$}}
\put(36,1.6){\makebox(0,0){$\bullet$}}
\put(28,1.5){\makebox(0,0){$\langle$}}
\put(4,1.5){\makebox(0,0){$\bullet$}}
\put(4,1.5){\line(1,0){16}}
\put(4,8){\makebox(0,0){$\scriptstyle 0$}}
\put(20,8){\makebox(0,0){$\scriptstyle -5$}}
\put(36,8){\makebox(0,0){$\scriptstyle 0$}}
\end{picture}
\rightarrow 0\end{array}$$
on any smooth $7$-manifold equipped with a generic distribution of rank~$4$. 
\end{theorem}
The bundles 
$\begin{picture}(42,5)
\put(21.4,.1){\line(1,0){14.3}}
\put(21.4,3.2){\line(1,0){14.3}}
\put(20,1.5){\makebox(0,0){$\times$}}
\put(36,1.6){\makebox(0,0){$\bullet$}}
\put(28,1.5){\makebox(0,0){$\langle$}}
\put(4,1.5){\makebox(0,0){$\bullet$}}
\put(4,1.5){\line(1,0){16}}
\put(4,8){\makebox(0,0){$\scriptstyle 0$}}
\put(20,8){\makebox(0,0){$\scriptstyle -3$}}
\put(36,8){\makebox(0,0){$\scriptstyle 2$}}
\end{picture}
+
\begin{picture}(42,5)
\put(21.4,.1){\line(1,0){14.3}}
\put(21.4,3.2){\line(1,0){14.3}}
\put(20,1.5){\makebox(0,0){$\times$}}
\put(36,1.6){\makebox(0,0){$\bullet$}}
\put(28,1.5){\makebox(0,0){$\langle$}}
\put(4,1.5){\makebox(0,0){$\bullet$}}
\put(4,1.5){\line(1,0){16}}
\put(4,8){\makebox(0,0){$\scriptstyle 3$}}
\put(20,8){\makebox(0,0){$\scriptstyle -4$}}
\put(36,8){\makebox(0,0){$\scriptstyle 1$}}
\end{picture}$ and
$\begin{picture}(42,5)
\put(21.4,.1){\line(1,0){14.3}}
\put(21.4,3.2){\line(1,0){14.3}}
\put(20,1.5){\makebox(0,0){$\times$}}
\put(36,1.6){\makebox(0,0){$\bullet$}}
\put(28,1.5){\makebox(0,0){$\langle$}}
\put(4,1.5){\makebox(0,0){$\bullet$}}
\put(4,1.5){\line(1,0){16}}
\put(4,8){\makebox(0,0){$\scriptstyle 3$}}
\put(20,8){\makebox(0,0){$\scriptstyle -6$}}
\put(36,8){\makebox(0,0){$\scriptstyle 1$}}
\end{picture}
+
\begin{picture}(42,5)
\put(21.4,.1){\line(1,0){14.3}}
\put(21.4,3.2){\line(1,0){14.3}}
\put(20,1.5){\makebox(0,0){$\times$}}
\put(36,1.6){\makebox(0,0){$\bullet$}}
\put(28,1.5){\makebox(0,0){$\langle$}}
\put(4,1.5){\makebox(0,0){$\bullet$}}
\put(4,1.5){\line(1,0){16}}
\put(4,8){\makebox(0,0){$\scriptstyle 0$}}
\put(20,8){\makebox(0,0){$\scriptstyle -6$}}
\put(36,8){\makebox(0,0){$\scriptstyle 2$}}
\end{picture}$
are sub-quotients of~$\Lambda^2$, respectively $\Lambda^5$. However, there is a
preferred class of splittings of the filtration of $\Lambda^1$ that canonically
splits these bundles as the following result shows.

\begin{theorem}\label{anotherpreferredsplitting}
The splittings of the short exact sequence 
$$0\rightarrow \begin{picture}(42,5)
\put(21.4,.1){\line(1,0){14.3}}
\put(21.4,3.2){\line(1,0){14.3}}
\put(20,1.5){\makebox(0,0){$\times$}}
\put(36,1.6){\makebox(0,0){$\bullet$}}
\put(28,1.5){\makebox(0,0){$\langle$}}
\put(4,1.5){\makebox(0,0){$\bullet$}}
\put(4,1.5){\line(1,0){16}}
\put(4,8){\makebox(0,0){$\scriptstyle 2$}}
\put(20,8){\makebox(0,0){$\scriptstyle -2$}}
\put(36,8){\makebox(0,0){$\scriptstyle 0$}}
\end{picture}\rightarrow
\Lambda^1
\rightarrow
\begin{picture}(42,5)
\put(21.4,.1){\line(1,0){14.3}}
\put(21.4,3.2){\line(1,0){14.3}}
\put(20,1.5){\makebox(0,0){$\times$}}
\put(36,1.6){\makebox(0,0){$\bullet$}}
\put(28,1.5){\makebox(0,0){$\langle$}}
\put(4,1.5){\makebox(0,0){$\bullet$}}
\put(4,1.5){\line(1,0){16}}
\put(4,8){\makebox(0,0){$\scriptstyle 1$}}
\put(20,8){\makebox(0,0){$\scriptstyle -2$}}
\put(36,8){\makebox(0,0){$\scriptstyle 1$}}
\end{picture}
\rightarrow 0$$
are acted freely upon by  
$${\mathrm{Hom}}\big( \begin{picture}(42,5)
\put(21.4,.1){\line(1,0){14.3}}
\put(21.4,3.2){\line(1,0){14.3}}
\put(20,1.5){\makebox(0,0){$\times$}}
\put(36,1.6){\makebox(0,0){$\bullet$}}
\put(28,1.5){\makebox(0,0){$\langle$}}
\put(4,1.5){\makebox(0,0){$\bullet$}}
\put(4,1.5){\line(1,0){16}}
\put(4,8){\makebox(0,0){$\scriptstyle 1$}}
\put(20,8){\makebox(0,0){$\scriptstyle -2$}}
\put(36,8){\makebox(0,0){$\scriptstyle 1$}}
\end{picture},\begin{picture}(42,5)
\put(21.4,.1){\line(1,0){14.3}}
\put(21.4,3.2){\line(1,0){14.3}}
\put(20,1.5){\makebox(0,0){$\times$}}
\put(36,1.6){\makebox(0,0){$\bullet$}}
\put(28,1.5){\makebox(0,0){$\langle$}}
\put(4,1.5){\makebox(0,0){$\bullet$}}
\put(4,1.5){\line(1,0){16}}
\put(4,8){\makebox(0,0){$\scriptstyle 2$}}
\put(20,8){\makebox(0,0){$\scriptstyle -2$}}
\put(36,8){\makebox(0,0){$\scriptstyle 0$}}
\end{picture}\big)= 
\begin{picture}(42,5)
\put(21.4,.1){\line(1,0){14.3}}
\put(21.4,3.2){\line(1,0){14.3}}
\put(20,1.5){\makebox(0,0){$\times$}}
\put(36,1.6){\makebox(0,0){$\bullet$}}
\put(28,1.5){\makebox(0,0){$\langle$}}
\put(4,1.5){\makebox(0,0){$\bullet$}}
\put(4,1.5){\line(1,0){16}}
\put(4,8){\makebox(0,0){$\scriptstyle 3$}}
\put(20,8){\makebox(0,0){$\scriptstyle -3$}}
\put(36,8){\makebox(0,0){$\scriptstyle 1$}}
\end{picture}\oplus \begin{picture}(42,5)
\put(21.4,.1){\line(1,0){14.3}}
\put(21.4,3.2){\line(1,0){14.3}}
\put(20,1.5){\makebox(0,0){$\times$}}
\put(36,1.6){\makebox(0,0){$\bullet$}}
\put(28,1.5){\makebox(0,0){$\langle$}}
\put(4,1.5){\makebox(0,0){$\bullet$}}
\put(4,1.5){\line(1,0){16}}
\put(4,8){\makebox(0,0){$\scriptstyle 1$}}
\put(20,8){\makebox(0,0){$\scriptstyle -2$}}
\put(36,8){\makebox(0,0){$\scriptstyle 1$}}
\end{picture}.$$ 
There is a preferred class of splittings in which the
\begin{picture}(42,5)
\put(21.4,.1){\line(1,0){14.3}}
\put(21.4,3.2){\line(1,0){14.3}}
\put(20,1.5){\makebox(0,0){$\times$}}
\put(36,1.6){\makebox(0,0){$\bullet$}}
\put(28,1.5){\makebox(0,0){$\langle$}}
\put(4,1.5){\makebox(0,0){$\bullet$}}
\put(4,1.5){\line(1,0){16}}
\put(4,8){\makebox(0,0){$\scriptstyle 3$}}
\put(20,8){\makebox(0,0){$\scriptstyle -3$}}
\put(36,8){\makebox(0,0){$\scriptstyle 1$}}
\end{picture}-freedom is eliminated. This restricted choice of splittings 
canonically splits the bundles 
\rule{0pt}{14pt}
$\begin{picture}(42,5)
\put(21.4,.1){\line(1,0){14.3}}
\put(21.4,3.2){\line(1,0){14.3}}
\put(20,1.5){\makebox(0,0){$\times$}}
\put(36,1.6){\makebox(0,0){$\bullet$}}
\put(28,1.5){\makebox(0,0){$\langle$}}
\put(4,1.5){\makebox(0,0){$\bullet$}}
\put(4,1.5){\line(1,0){16}}
\put(4,8){\makebox(0,0){$\scriptstyle 0$}}
\put(20,8){\makebox(0,0){$\scriptstyle -3$}}
\put(36,8){\makebox(0,0){$\scriptstyle 2$}}
\end{picture}
+
\begin{picture}(42,5)
\put(21.4,.1){\line(1,0){14.3}}
\put(21.4,3.2){\line(1,0){14.3}}
\put(20,1.5){\makebox(0,0){$\times$}}
\put(36,1.6){\makebox(0,0){$\bullet$}}
\put(28,1.5){\makebox(0,0){$\langle$}}
\put(4,1.5){\makebox(0,0){$\bullet$}}
\put(4,1.5){\line(1,0){16}}
\put(4,8){\makebox(0,0){$\scriptstyle 3$}}
\put(20,8){\makebox(0,0){$\scriptstyle -4$}}
\put(36,8){\makebox(0,0){$\scriptstyle 1$}}
\end{picture}$ and
$\begin{picture}(42,5)
\put(21.4,.1){\line(1,0){14.3}}
\put(21.4,3.2){\line(1,0){14.3}}
\put(20,1.5){\makebox(0,0){$\times$}}
\put(36,1.6){\makebox(0,0){$\bullet$}}
\put(28,1.5){\makebox(0,0){$\langle$}}
\put(4,1.5){\makebox(0,0){$\bullet$}}
\put(4,1.5){\line(1,0){16}}
\put(4,8){\makebox(0,0){$\scriptstyle 3$}}
\put(20,8){\makebox(0,0){$\scriptstyle -6$}}
\put(36,8){\makebox(0,0){$\scriptstyle 1$}}
\end{picture}
+
\begin{picture}(42,5)
\put(21.4,.1){\line(1,0){14.3}}
\put(21.4,3.2){\line(1,0){14.3}}
\put(20,1.5){\makebox(0,0){$\times$}}
\put(36,1.6){\makebox(0,0){$\bullet$}}
\put(28,1.5){\makebox(0,0){$\langle$}}
\put(4,1.5){\makebox(0,0){$\bullet$}}
\put(4,1.5){\line(1,0){16}}
\put(4,8){\makebox(0,0){$\scriptstyle 0$}}
\put(20,8){\makebox(0,0){$\scriptstyle -6$}}
\put(36,8){\makebox(0,0){$\scriptstyle 2$}}
\end{picture}$.
\end{theorem}
\begin{proof} The only difficulty is in restricting the class of splittings
and, as usual, one looks to the exterior derivative $d:\Lambda^2\to\Lambda^3$
in the presence of a chosen splitting. More specifically, one checks (e.g., in
an adapted co-frame) that, having chosen a splitting of~$\Lambda^1$, the
resulting component of $d:\Lambda^2\to\Lambda^3$ mapping
$\begin{picture}(42,5)
\put(21.4,.1){\line(1,0){14.3}}
\put(21.4,3.2){\line(1,0){14.3}}
\put(20,1.5){\makebox(0,0){$\times$}}
\put(36,1.6){\makebox(0,0){$\bullet$}}
\put(28,1.5){\makebox(0,0){$\langle$}}
\put(4,1.5){\makebox(0,0){$\bullet$}}
\put(4,1.5){\line(1,0){16}}
\put(4,8){\makebox(0,0){$\scriptstyle 1$}}
\put(20,8){\makebox(0,0){$\scriptstyle -3$}}
\put(36,8){\makebox(0,0){$\scriptstyle 1$}}
\end{picture}$ to 
$\begin{picture}(42,5)
\put(21.4,.1){\line(1,0){14.3}}
\put(21.4,3.2){\line(1,0){14.3}}
\put(20,1.5){\makebox(0,0){$\times$}}
\put(36,1.6){\makebox(0,0){$\bullet$}}
\put(28,1.5){\makebox(0,0){$\langle$}}
\put(4,1.5){\makebox(0,0){$\bullet$}}
\put(4,1.5){\line(1,0){16}}
\put(4,8){\makebox(0,0){$\scriptstyle 4$}}
\put(20,8){\makebox(0,0){$\scriptstyle -4$}}
\put(36,8){\makebox(0,0){$\scriptstyle 0$}}
\end{picture}$ is actually a homomorphism. However,
\begin{equation}\label{torsionishere}
\begin{array}{rcc}{\mathrm{Hom}}\big(\begin{picture}(42,5)
\put(21.4,.1){\line(1,0){14.3}}
\put(21.4,3.2){\line(1,0){14.3}}
\put(20,1.5){\makebox(0,0){$\times$}}
\put(36,1.6){\makebox(0,0){$\bullet$}}
\put(28,1.5){\makebox(0,0){$\langle$}}
\put(4,1.5){\makebox(0,0){$\bullet$}}
\put(4,1.5){\line(1,0){16}}
\put(4,8){\makebox(0,0){$\scriptstyle 1$}}
\put(20,8){\makebox(0,0){$\scriptstyle -3$}}
\put(36,8){\makebox(0,0){$\scriptstyle 1$}}
\end{picture},\begin{picture}(42,5)
\put(21.4,.1){\line(1,0){14.3}}
\put(21.4,3.2){\line(1,0){14.3}}
\put(20,1.5){\makebox(0,0){$\times$}}
\put(36,1.6){\makebox(0,0){$\bullet$}}
\put(28,1.5){\makebox(0,0){$\langle$}}
\put(4,1.5){\makebox(0,0){$\bullet$}}
\put(4,1.5){\line(1,0){16}}
\put(4,8){\makebox(0,0){$\scriptstyle 4$}}
\put(20,8){\makebox(0,0){$\scriptstyle -4$}}
\put(36,8){\makebox(0,0){$\scriptstyle 0$}}
\end{picture}\big)&=&\begin{picture}(42,5)
\put(21.4,.1){\line(1,0){14.3}}
\put(21.4,3.2){\line(1,0){14.3}}
\put(20,1.5){\makebox(0,0){$\times$}}
\put(36,1.6){\makebox(0,0){$\bullet$}}
\put(28,1.5){\makebox(0,0){$\langle$}}
\put(4,1.5){\makebox(0,0){$\bullet$}}
\put(4,1.5){\line(1,0){16}}
\put(4,8){\makebox(0,0){$\scriptstyle 1$}}
\put(20,8){\makebox(0,0){$\scriptstyle 0$}}
\put(36,8){\makebox(0,0){$\scriptstyle 1$}}
\end{picture}\otimes\begin{picture}(42,5)
\put(21.4,.1){\line(1,0){14.3}}
\put(21.4,3.2){\line(1,0){14.3}}
\put(20,1.5){\makebox(0,0){$\times$}}
\put(36,1.6){\makebox(0,0){$\bullet$}}
\put(28,1.5){\makebox(0,0){$\langle$}}
\put(4,1.5){\makebox(0,0){$\bullet$}}
\put(4,1.5){\line(1,0){16}}
\put(4,8){\makebox(0,0){$\scriptstyle 4$}}
\put(20,8){\makebox(0,0){$\scriptstyle -4$}}
\put(36,8){\makebox(0,0){$\scriptstyle 0$}}
\end{picture}\\[5pt]
&&\|\\[5pt]
&&\begin{picture}(42,5)
\put(21.4,.1){\line(1,0){14.3}}
\put(21.4,3.2){\line(1,0){14.3}}
\put(20,1.5){\makebox(0,0){$\times$}}
\put(36,1.6){\makebox(0,0){$\bullet$}}
\put(28,1.5){\makebox(0,0){$\langle$}}
\put(4,1.5){\makebox(0,0){$\bullet$}}
\put(4,1.5){\line(1,0){16}}
\put(4,8){\makebox(0,0){$\scriptstyle 5$}}
\put(20,8){\makebox(0,0){$\scriptstyle -4$}}
\put(36,8){\makebox(0,0){$\scriptstyle 1$}}
\end{picture}\oplus\begin{picture}(42,5)
\put(21.4,.1){\line(1,0){14.3}}
\put(21.4,3.2){\line(1,0){14.3}}
\put(20,1.5){\makebox(0,0){$\times$}}
\put(36,1.6){\makebox(0,0){$\bullet$}}
\put(28,1.5){\makebox(0,0){$\langle$}}
\put(4,1.5){\makebox(0,0){$\bullet$}}
\put(4,1.5){\line(1,0){16}}
\put(4,8){\makebox(0,0){$\scriptstyle 3$}}
\put(20,8){\makebox(0,0){$\scriptstyle -3$}}
\put(36,8){\makebox(0,0){$\scriptstyle 1$}}
\end{picture}\end{array}\end{equation}
and one checks (again using an adapted co-frame or by arguing with irreducible
bundles and Schur's lemma) that the original freedom in splitting $\Lambda^1$
can be used to eliminate the
$\begin{picture}(42,5)
\put(21.4,.1){\line(1,0){14.3}}
\put(21.4,3.2){\line(1,0){14.3}}
\put(20,1.5){\makebox(0,0){$\times$}}
\put(36,1.6){\makebox(0,0){$\bullet$}}
\put(28,1.5){\makebox(0,0){$\langle$}}
\put(4,1.5){\makebox(0,0){$\bullet$}}
\put(4,1.5){\line(1,0){16}}
\put(4,8){\makebox(0,0){$\scriptstyle 3$}}
\put(20,8){\makebox(0,0){$\scriptstyle -3$}}
\put(36,8){\makebox(0,0){$\scriptstyle 1$}}
\end{picture}$-component. The remaining freedom in splitting 
$\Lambda^1$ is thereby restricted in the appropriate manner (with the stated
knock-on effect on the previously identified subquotients of $\Lambda^2$
and~$\Lambda^5$).
\end{proof}

Of course, one could rephrase Theorem~\ref{anotherpreferredsplitting} as
defining the notion of a $2$-adapted co-framing and observe that the effect is
to reduce the structure bundle of $\Lambda^1$ to~$P/\exp({\mathfrak{g}}_2)$. 
In any case, combining the two theorems above we immediately obtain the 
following improved complex.
\begin{theorem}
On any smooth $7$-manifold endowed with a generic $4$-distribution, there is a
canonically defined locally exact differential complex
$$\begin{array}{l}
0\rightarrow{\mathbb{R}}\rightarrow 
\begin{picture}(42,5)
\put(21.4,.1){\line(1,0){14.3}}
\put(21.4,3.2){\line(1,0){14.3}}
\put(20,1.5){\makebox(0,0){$\times$}}
\put(36,1.6){\makebox(0,0){$\bullet$}}
\put(28,1.5){\makebox(0,0){$\langle$}}
\put(4,1.5){\makebox(0,0){$\bullet$}}
\put(4,1.5){\line(1,0){16}}
\put(4,8){\makebox(0,0){$\scriptstyle 1$}}
\put(20,8){\makebox(0,0){$\scriptstyle -2$}}
\put(36,8){\makebox(0,0){$\scriptstyle 1$}}
\end{picture}
\rightarrow\!\!
\begin{array}{c}
\begin{picture}(42,5)
\put(21.4,.1){\line(1,0){14.3}}
\put(21.4,3.2){\line(1,0){14.3}}
\put(20,1.5){\makebox(0,0){$\times$}}
\put(36,1.6){\makebox(0,0){$\bullet$}}
\put(28,1.5){\makebox(0,0){$\langle$}}
\put(4,1.5){\makebox(0,0){$\bullet$}}
\put(4,1.5){\line(1,0){16}}
\put(4,8){\makebox(0,0){$\scriptstyle 0$}}
\put(20,8){\makebox(0,0){$\scriptstyle -3$}}
\put(36,8){\makebox(0,0){$\scriptstyle 2$}}
\end{picture}\\
\oplus\\[3pt]
\begin{picture}(42,5)
\put(21.4,.1){\line(1,0){14.3}}
\put(21.4,3.2){\line(1,0){14.3}}
\put(20,1.5){\makebox(0,0){$\times$}}
\put(36,1.6){\makebox(0,0){$\bullet$}}
\put(28,1.5){\makebox(0,0){$\langle$}}
\put(4,1.5){\makebox(0,0){$\bullet$}}
\put(4,1.5){\line(1,0){16}}
\put(4,8){\makebox(0,0){$\scriptstyle 3$}}
\put(20,8){\makebox(0,0){$\scriptstyle -4$}}
\put(36,8){\makebox(0,0){$\scriptstyle 1$}}
\end{picture}
\end{array}
\!\!\rightarrow\!\!
\begin{array}{c}
\begin{picture}(42,5)
\put(21.4,.1){\line(1,0){14.3}}
\put(21.4,3.2){\line(1,0){14.3}}
\put(20,1.5){\makebox(0,0){$\times$}}
\put(36,1.6){\makebox(0,0){$\bullet$}}
\put(28,1.5){\makebox(0,0){$\langle$}}
\put(4,1.5){\makebox(0,0){$\bullet$}}
\put(4,1.5){\line(1,0){16}}
\put(4,8){\makebox(0,0){$\scriptstyle 2$}}
\put(20,8){\makebox(0,0){$\scriptstyle -5$}}
\put(36,8){\makebox(0,0){$\scriptstyle 2$}}
\end{picture}\\
\oplus\\[3pt]
\begin{picture}(42,5)
\put(21.4,.1){\line(1,0){14.3}}
\put(21.4,3.2){\line(1,0){14.3}}
\put(20,1.5){\makebox(0,0){$\times$}}
\put(36,1.6){\makebox(0,0){$\bullet$}}
\put(28,1.5){\makebox(0,0){$\langle$}}
\put(4,1.5){\makebox(0,0){$\bullet$}}
\put(4,1.5){\line(1,0){16}}
\put(4,8){\makebox(0,0){$\scriptstyle 4$}}
\put(20,8){\makebox(0,0){$\scriptstyle -4$}}
\put(36,8){\makebox(0,0){$\scriptstyle 0$}}
\end{picture}
\end{array}
\!\!\rightarrow\!\!
\begin{array}{c}
\begin{picture}(42,5)
\put(21.4,.1){\line(1,0){14.3}}
\put(21.4,3.2){\line(1,0){14.3}}
\put(20,1.5){\makebox(0,0){$\times$}}
\put(36,1.6){\makebox(0,0){$\bullet$}}
\put(28,1.5){\makebox(0,0){$\langle$}}
\put(4,1.5){\makebox(0,0){$\bullet$}}
\put(4,1.5){\line(1,0){16}}
\put(4,8){\makebox(0,0){$\scriptstyle 4$}}
\put(20,8){\makebox(0,0){$\scriptstyle -5$}}
\put(36,8){\makebox(0,0){$\scriptstyle 0$}}
\end{picture}\\
\oplus\\[3pt]
\begin{picture}(42,5)
\put(21.4,.1){\line(1,0){14.3}}
\put(21.4,3.2){\line(1,0){14.3}}
\put(20,1.5){\makebox(0,0){$\times$}}
\put(36,1.6){\makebox(0,0){$\bullet$}}
\put(28,1.5){\makebox(0,0){$\langle$}}
\put(4,1.5){\makebox(0,0){$\bullet$}}
\put(4,1.5){\line(1,0){16}}
\put(4,8){\makebox(0,0){$\scriptstyle 2$}}
\put(20,8){\makebox(0,0){$\scriptstyle -6$}}
\put(36,8){\makebox(0,0){$\scriptstyle 2$}}
\end{picture}
\end{array}
\!\!\rightarrow\\
\qquad\begin{array}{c}
\begin{picture}(42,5)
\put(21.4,.1){\line(1,0){14.3}}
\put(21.4,3.2){\line(1,0){14.3}}
\put(20,1.5){\makebox(0,0){$\times$}}
\put(36,1.6){\makebox(0,0){$\bullet$}}
\put(28,1.5){\makebox(0,0){$\langle$}}
\put(4,1.5){\makebox(0,0){$\bullet$}}
\put(4,1.5){\line(1,0){16}}
\put(4,8){\makebox(0,0){$\scriptstyle 3$}}
\put(20,8){\makebox(0,0){$\scriptstyle -6$}}
\put(36,8){\makebox(0,0){$\scriptstyle 1$}}
\end{picture}\\
\oplus\\[3pt]
\begin{picture}(42,5)
\put(21.4,.1){\line(1,0){14.3}}
\put(21.4,3.2){\line(1,0){14.3}}
\put(20,1.5){\makebox(0,0){$\times$}}
\put(36,1.6){\makebox(0,0){$\bullet$}}
\put(28,1.5){\makebox(0,0){$\langle$}}
\put(4,1.5){\makebox(0,0){$\bullet$}}
\put(4,1.5){\line(1,0){16}}
\put(4,8){\makebox(0,0){$\scriptstyle 0$}}
\put(20,8){\makebox(0,0){$\scriptstyle -6$}}
\put(36,8){\makebox(0,0){$\scriptstyle 2$}}
\end{picture}
\end{array}
\!\!\rightarrow
\begin{picture}(42,5)
\put(21.4,.1){\line(1,0){14.3}}
\put(21.4,3.2){\line(1,0){14.3}}
\put(20,1.5){\makebox(0,0){$\times$}}
\put(36,1.6){\makebox(0,0){$\bullet$}}
\put(28,1.5){\makebox(0,0){$\langle$}}
\put(4,1.5){\makebox(0,0){$\bullet$}}
\put(4,1.5){\line(1,0){16}}
\put(4,8){\makebox(0,0){$\scriptstyle 1$}}
\put(20,8){\makebox(0,0){$\scriptstyle -6$}}
\put(36,8){\makebox(0,0){$\scriptstyle 1$}}
\end{picture}
\rightarrow
\begin{picture}(42,5)
\put(21.4,.1){\line(1,0){14.3}}
\put(21.4,3.2){\line(1,0){14.3}}
\put(20,1.5){\makebox(0,0){$\times$}}
\put(36,1.6){\makebox(0,0){$\bullet$}}
\put(28,1.5){\makebox(0,0){$\langle$}}
\put(4,1.5){\makebox(0,0){$\bullet$}}
\put(4,1.5){\line(1,0){16}}
\put(4,8){\makebox(0,0){$\scriptstyle 0$}}
\put(20,8){\makebox(0,0){$\scriptstyle -5$}}
\put(36,8){\makebox(0,0){$\scriptstyle 0$}}
\end{picture}
\rightarrow 0.\end{array}$$
\end{theorem}
Of course, on the homogeneous model $G/P$ this is the standard BGG complex.
Finally, as mentioned already in the introduction, let us reiterate that our
construction does not see the torsion of this parabolic geometry. In fact, the
torsion lies in
$\begin{picture}(42,5)
\put(21.4,.1){\line(1,0){14.3}}
\put(21.4,3.2){\line(1,0){14.3}}
\put(20,1.5){\makebox(0,0){$\times$}}
\put(36,1.6){\makebox(0,0){$\bullet$}}
\put(28,1.5){\makebox(0,0){$\langle$}}
\put(4,1.5){\makebox(0,0){$\bullet$}}
\put(4,1.5){\line(1,0){16}}
\put(4,8){\makebox(0,0){$\scriptstyle 5$}}
\put(20,8){\makebox(0,0){$\scriptstyle -4$}}
\put(36,8){\makebox(0,0){$\scriptstyle 1$}}
\end{picture}$,
which is exactly the component of (\ref{torsionishere}) that we have not
eliminated.

\section{The Rumin-Seshadri complex}
Although not a replacement for the de~Rham complex in resolving the constants,
we take the opportunity here to describe another natural differential complex,
the Rumin-Seshadri complex~\cite{S}, the construction of which follows the same
general technique. This complex is defined on any symplectic $2n$-manifold $M$
as follows. Denoting by $J$ the symplectic $2$-form, let us consider the
filtered differential complex $E^\bullet$ defined by
$$E^p=\Lambda^p\oplus\Lambda^{p-1}\enskip\mbox{for }
p=0,1,\ldots,{2n+1}$$
with differentials
$$(\omega,\mu)\mapsto(d\omega+(-1)^pJ\wedge\mu,d\mu).$$
Notice that this complex has local cohomology at both $p=0$ and $p=1$.
Specifically, the kernel of $E^0\to E^1$ is 
$\{(f,0)\mbox{ s.t.\ }f\mbox{ is constant}\}$ and the cohomology at $p=1$ is
generated by $(\alpha,-1)$, where $\alpha$ is any local potential for the
symplectic form $J$, meaning that $d\alpha=J$. Evidently, this is a filtered
complex: the de~Rham complex is a sub-complex. The associated spectral
sequence immediately gives rise to the Rumin-Seshadri differential complex 
on~$M$. It has the form
\begin{equation}
\label{RScomplex}\begin{array}{ccccccccccc}
\Lambda^0&\stackrel{d}{\longrightarrow}&\Lambda^1
&\stackrel{d_\perp}{\longrightarrow}&\Lambda_\perp^2
&\stackrel{d_\perp}{\longrightarrow}&\Lambda_\perp^3
&\stackrel{d_\perp}{\longrightarrow}&\cdots
&\stackrel{d_\perp}{\longrightarrow}&\Lambda_\perp^{n}\\[2pt]
&&&&&&&&&&\big\downarrow\makebox[0pt][l]{\scriptsize$d_\perp^{(2)}$}\\
\Lambda^0&\stackrel{d_\perp}{\longleftarrow}&\Lambda^1
&\stackrel{d_\perp}{\longleftarrow}&\Lambda_\perp^2
&\stackrel{d_\perp}{\longleftarrow}&\Lambda_\perp^3
&\stackrel{d_\perp}{\longleftarrow}&\cdots
&\stackrel{d_\perp}{\longleftarrow}&\Lambda_\perp^{n}
\end{array}\end{equation}
where $\Lambda_\perp^p$ denotes the $p$-forms that are trace-free with
respect to~$J$.  The conclusion is as follows.
\begin{theorem}
On any symplectic manifold, there is a differential 
complex\/ {\rm (\ref{RScomplex})} with local cohomology in 
degrees\/ $0$ and~$\,1$. 
On the level of sheaves, 
in both these degrees the cohomology is the locally constant 
sheaf~$\,{\mathbb{R}}$. In all other degrees it is 
locally exact. On a compact symplectic manifold of 
dimension~$\geq 4$,
\begin{equation}\label{firstcohomology}
\frac{{\mathrm{ker}}\,d_\perp:\Gamma(M,\Lambda^1)\to
\Gamma(M,\Lambda_\perp^2)}
{{\mathrm{im}}\,d:\Gamma(M,\Lambda^0)\to
\Gamma(M,\Lambda^1)}\cong H^1(M,{\mathbb{R}}).\end{equation}
\end{theorem} 
\begin{proof} The construction of the complex and the identification of
its local cohomology are immediate form the spectral sequence. To see
(\ref{firstcohomology}), note that for a $1$-form $\omega$ to be in the
kernel of $d_\perp$ is to say that $d\omega=fJ$ for some smooth 
function~$f$ but then
$$\begin{array}{ccccl}
0=d^2\omega=df\wedge J&\implies&df=0&\implies&
f\mbox{ is constant}\\
&&&\implies&f=0
\enskip\mbox{or}\enskip J=d(\omega/f).\end{array}$$
However, the symplectic form cannot be exact for $M$ compact so 
$f=0$ and thus $d\omega=0$. 
\end{proof}
In four dimensions, the complex (\ref{RScomplex}) is due to
R.T.~Smith~\cite{Sm}. In higher dimensions, it was also found by L.-S.~Tseng
and S.-T.~Yau~\cite{SY} who show that it is elliptic and go on to study its
cohomology on compact manifolds. The complex of first order operators after the
second-order operator in the middle, was introduced by T.~Bouche~\cite{Bo} and
who dubbed it the {\em coeffective\/} complex (he regarded it as a subcomplex
of the second half of the de~Rham complex $\Lambda^{n}\to\cdots\Lambda^{2n}$).
The coeffective cohomology was further studied by M.~Fern\'andez,
R.~Ib\'a\~nez, and M. de~Le\'on (see, for example,~\cite{fil}).

\section*{Acknowledgements}
\noindent RLB gratefully acknowledges NSF support from grants DMS-8352009 and
DMS-8905207 (from the 1980s when some of this work was done), and current NSF
support from DMS-1105868\@. MGE, ARG, and KN would like to thank the Erwin
Schr\"odinger Institute for hospitality in July 2011 during which this work was
crucially advanced. MGE gratefully acknowledges support from the Australian
Research Council. ARG gratefully acknowledges support from the Royal Society of
New Zealand via Marsden Grant 10-UOA-113.

We thank Boris Doubrov for pointing out to us the homogeneous Engel manifold 
of~\S\ref{engelrevisited}, Jean-Pierre Demailly for drawing our 
attention to~\cite{fil}, and Li-Sheng Tseng for drawing our attention 
to~\cite{Sm}.

\newpage


\vspace{-10pt}\begin{thebibliography}{99}

\bibitem{BE} R.J. Baston and M.G. Eastwood,
The Penrose Transform: its Interaction with Representation Theory,
Oxford University Press, 1989.

\bibitem{BGG} I.N. Bernstein, I.M. Gelfand, and S.I. Gelfand,
Differential operators on the base affine space and a study 
of ${\mathfrak{g}}$-modules,
In: Lie Groups and their Representations,
Halsted, New York, 1975, pp.~21--64.

\bibitem{Bo} T. Bouche,
La cohomologie coeffective d'une vari\'et\'e symplectique,
Bull. Sci. Math., {\bf 114} (1990), 115--122. 

\bibitem{B} R.L. Bryant,
Conformal geometry and $3$-plane fields on $6$-manifolds,
In: Proceedings of the RIMS symposium 
`Developments of Cartan Geometry and Related Mathematical Problems' 
(24--27 October 2005),
Kokyuroku, {\bf 1502} (2006), pp.~1--15. 

\bibitem{CD} D.M.J. Calderbank and T. Diemer,
Differential invariants and curved Bernstein-Gelfand-Gelfand sequences,
Jour. Reine Angew. Math., {\bf 537} (2001), 67--103.

\bibitem{CSc} A. \v{C}ap and H. Schichl, 
Parabolic geometries and canonical Cartan connections, 
Hokkaido Math. Jour., {\bf 29} (2001), 453--505.

\bibitem{CSl} A. \v{C}ap and J. Slov\'ak, 
Parabolic Geometries I: Background and General Theory,
Math. Surv. and Monographs, {\bf 154},
Amer. Math. Soc., Providence, RI, 2009.

\bibitem{CSS} A. \v{C}ap, J. Slov\'ak, and V. Sou\v{c}ek, 
Bernstein-Gelfand-Gelfand sequences,
Ann. Math., {\bf 154} (2001), 97--113.

\bibitem{Ca} E. Cartan,
Les syst\`emes de Pfaff \'a cinq variables et les \'equations aux
deriv\'ees partielles du second ordre,
Ann. Ecole Norm. Sup., {\bf 27} (1910), 109--192.

\bibitem{Ch} T.Y. Chow,
You could have invented spectral sequences,
Notices Amer. Math. Soc., {\bf 53} (2006), 15--19.

\bibitem{fil} M.~Fern\'andez, R.~Ib\'a\~nez, and M. de~Le\'on,
Coeffective and de Rham cohomologies of symplectic manifolds,
Jour. Geom. Phys., {\bf 27} (1998), 281--296. 

\bibitem{M} R.G. Montgomery,
Engel deformations and contact structures,
In: Northern California Symplectic Geometry Seminar,
Amer. Math. Soc. Transl., {\bf 196}, 
Amer. Math. Soc., Providence, RI, 1999, pp.~103--117.

\bibitem{M2} R.G. Montgomery,
A Tour of Subriemannian Geometries, their Geodesics and Applications,
Math. Surv. and Monographs, {\bf 91},
Amer. Math. Soc., Providence, RI, 2002.

\bibitem{P} J. Peetre,
R\'ectification \`a l'article `une caract\'erisation abstraite des 
op\'era\-teurs diff\'erentiels,'
Math. Scand., {\bf 8} (1960), 116--120. 

\bibitem{R} M. Rumin,
Un complexe de formes diff\'erentielles sur les vari\'et\'es de contact,
Comptes Rendus Acad. Sci. Paris Math., {\bf 310} (1990), 401--404.

\bibitem{S} N. Seshadri,
Private communication, 
September 2007.

\bibitem{Sm} R.T. Smith,
Examples of elliptic complexes,
Bull. Amer. Math. Soc., {\bf 82} (1976), 297--299. 

\bibitem{SY} L.-S. Tseng and S.-T. Yau,
Cohomology and Hodge theory on symplectic manifolds: I and II,
arXiv:0909.5418 and arXiv:1011.1250.

\end{thebibliography}
\end{document}